\documentclass[reqno, 11pt]{amsart}
\usepackage[left=2cm, right=2cm, top=2cm, bottom=2cm]{geometry}
\usepackage[T1]{fontenc}
\usepackage[utf8]{inputenc}	
\usepackage[english]{babel}
\usepackage[babel]{csquotes}
\usepackage{bbm}
\usepackage{amssymb}
\usepackage{amsmath}
\usepackage{amsthm}
\usepackage{latexsym}
\usepackage{xcolor}
\usepackage{mathrsfs}
\usepackage{bbm}
\usepackage{verbatim}
\usepackage[colorlinks=true, pdfstartview=FitV, linkcolor=blue, citecolor=blue, urlcolor=blue]{hyperref}
\usepackage[style=numeric, sorting=nty, sortcites=true, backend=biber, maxbibnames=10, giveninits=true]{biblatex}
\usepackage{faktor}
\usepackage[shortlabels]{enumitem}

\numberwithin{equation}{section}
\newtheorem{thm}{Theorem}[section]
\newtheorem{cor}[thm]{Corollary}
\newtheorem{lem}[thm]{Lemma}

\theoremstyle{definition}
\newtheorem{defin}[thm]{Definition}
\newtheorem{remark}[thm]{Remark}

\renewcommand{\d}{{\mathrm d}} 
\newcommand{\norm}[1]{\left\|#1\right\|} 
\newcommand{\ip}[2]{\left\langle #1,#2 \right\rangle} 
\def\laweq{\stackrel{\mathcal L}{=}}
\def\tom{\widetilde{\Omega}}
\def\tF{\widetilde{\mathscr{F}}}
\def\tP{\widetilde{\mathbb{P}}}
\def\tE{\widetilde{\mathbb{E}}}

\def\hom{\widehat{\Omega}}
\def\hF{\widehat{\mathscr{F}}}
\def\hP{\widehat{\mathbb{P}}}
\def\hE{\widehat{\mathbb{E}}}
\def\hW{\widehat W}

\def\hphi{\widehat\varphi}
\def\hmu{\widehat\mu}
\def\pom{{\Omega}'}
\def\pE{\E\:\hspace{-1.4mm}'}

\def\ki{\varepsilon_{k_i}}
\def\ji{\varepsilon_{j_i}}

\def\bn{{\boldsymbol n}}

\def\bH{{\boldsymbol H}}

\def\b #1{{\boldsymbol #1}}
\def\enne{\mathbb{N}}

\def\erre{\mathbb{R}}
\def\R{\mathbb{R}}

\def\P{\mathbb{P}}

\def\E{\mathop{{}\mathbb{E}}}

\def\cL{\mathscr{L}}
\def\cF{\mathscr{F}}

\def\eps{\varepsilon}
\def\cP{\mathscr{P}}

\def\OO{\mathcal{O}}
\def\embed{\hookrightarrow}

\def\t{{\min\{t, \zeta_n\}}}

\def\sups{\sup_{s \in [0,\,\t]}}
\def\supt{\sup_{\tau \in [0,\,t]}}

\def\christoph #1{{\color{black} #1}}
\def\review #1{{\color{black}#1}}

\makeatletter
\DeclareFontFamily{OMX}{MnSymbolE}{}
\DeclareSymbolFont{MnLargeSymbols}{OMX}{MnSymbolE}{m}{n}
\SetSymbolFont{MnLargeSymbols}{bold}{OMX}{MnSymbolE}{b}{n}
\DeclareFontShape{OMX}{MnSymbolE}{m}{n}{
	<-6>  MnSymbolE5
	<6-7>  MnSymbolE6
	<7-8>  MnSymbolE7
	<8-9>  MnSymbolE8
	<9-10> MnSymbolE9
	<10-12> MnSymbolE10
	<12->   MnSymbolE12
}{}
\DeclareFontShape{OMX}{MnSymbolE}{b}{n}{
	<-6>  MnSymbolE-Bold5
	<6-7>  MnSymbolE-Bold6
	<7-8>  MnSymbolE-Bold7
	<8-9>  MnSymbolE-Bold8
	<9-10> MnSymbolE-Bold9
	<10-12> MnSymbolE-Bold10
	<12->   MnSymbolE-Bold12
}{}

\let\llangle\@undefined
\let\rrangle\@undefined
\DeclareMathDelimiter{\llangle}{\mathopen}%
{MnLargeSymbols}{'164}{MnLargeSymbols}{'164}
\DeclareMathDelimiter{\rrangle}{\mathclose}%
{MnLargeSymbols}{'171}{MnLargeSymbols}{'171}
\makeatother
\AtEveryBibitem{%
\clearfield{doi}%
\clearfield{url}%
\clearfield{issn}%
\clearfield{isbn}
} 
\renewbibmacro*{publisher+location+date}{%
\ifentrytype{book}
{\printlist{publisher}%
\iflistundef{location}
{\setunit*{\addcomma\space}}
{\setunit*{\addcomma\space}}%
\printlist{location}}
{\printlist{location}%
\iflistundef{publisher}
{\setunit*{\addcomma\space}}
{\setunit*{\addcomma\space}}%
\printlist{publisher}}
\setunit*{\addcomma\space}%
\usebibmacro{date}%
\newunit} 

\usepackage{xcolor}
\begin{document}
    \title[The stochastic nonlocal CH equation with regular potential and multiplicative noise]{The stochastic nonlocal Cahn--Hilliard equation \\ with regular potential and multiplicative noise}
    
    \author[Andrea Di Primio]{Andrea Di Primio$^*$}	
    \address{Classe di Scienze, Scuola Normale Superiore, 56126 Pisa, Italy}
    \email{andrea.diprimio@sns.it}
    
    \author{Christoph Hurm}
    \address{Fakult\"{a}t f\"{u}r	 Mathematik, Universit\"{a}t Regensburg, 93040 Regensburg, Germany}
    \email{christoph.hurm@ur.de}

    \thanks{${}^*$Corresponding author.}

    \begin{abstract}
	In this work, we deal with the stochastic counterpart of the nonlocal Cahn--Hilliard equation with regular potential in a smooth bounded one-, two- or three-dimensional domain. The problem is endowed with homogeneous Neumann boundary conditions and random initial data. Furthermore, the system is driven by cylindrical noise of multiplicative type. For the resulting system, we are able to show the existence of probabilistically-weak (or martingale) solutions in two and three dimensions, that are unique and probabilistically-strong under suitable assumptions on the stochastic diffusion. Moreover, we investigate the nonlocal-to-local asymptotics toward solutions of the local stochastic Cahn--Hilliard equations, establishing, under regularity conditions, a precise rate of convergence as well.
    \end{abstract}
    \maketitle

    \section{Introduction}
    \label{sec:intro}
    Phase separation is a peculiar physical phenomenon arising in certain multicompoment systems when subject to thermodynamical stress. A canonical example of this behavior is given by a binary homogeneous mixture that is rapidly quenched under a critical temperature. Under these conditions, separation dynamics may occur: the two components of the system reorganize through a pattern formation process, often leading, at equilibrium, to the formation of subregions that are rich in a single phase. The dynamics of phase separation generally unfold through two successive regimes, each characterized by different temporal and spatial scales. The initial stage, called spinodal decomposition, is the driver of the early stages of phase separation: microstructures form at small spatial scales as a result of thermodynamical instabilities. On longer time scales, the system enters the coarsening stage: macroscopic patterns emerge as multiple microstructures merge into larger regions, eventually towards some equilibrium state. From the modeling standpoint, the first complete description of phase separation dates back to the pioneering contributions of J. W. Cahn and J. E. Hilliard, in the context of metallic binary alloys (see, e.g., \cite{CahnHilliard58}). Precisely, for $d \in \{1,2,3\}$, let $\OO \subset \mathbb{R}^d$ be a bounded smooth domain and fix a time horizon $T > 0$. The well-known Cahn--Hilliard equation for a binary system reads
    \begin{equation} \label{eq:CH1}
    \begin{cases}
        \partial_t \varphi - \Delta \mu =0 & \quad \text{in } (0,T) \times \OO, \\
         \mu = -\Delta\varphi+ F'(\varphi)& \quad \text{in } (0,T) \times \OO, \\
         \partial_{\bn} \mu = \partial_{\bn} \varphi = 0 & \quad \text{on } (0,T) \times \partial\OO, \\
         \varphi(0) = \varphi_0 & \quad \text{in } \OO,
    \end{cases}
    \end{equation}
    where $\varphi$ denotes the relative concentration difference between the two components, $\mu$ denotes the chemical potential and $F$ represents a potential energy density for the system. The choice of the potential $F$ plays a crucial role when evaluating how close \eqref{eq:CH1} is to an actual physical model. Indeed, since $\varphi$ denotes a difference of two concentrations, it is natural to expect the physically meaningful constraint
    \begin{equation} \label{eq:physical}
        \varphi(x,t) \in [-1,1] \qquad \forall \: (x,t) \in  \OO \times [0,T),
    \end{equation}
    so that, at a fixed time $t \geq 0$, the level sets $\{\varphi(\cdot,\, t) = 1\}$ and $\{\varphi(\cdot,\, t) = -1\}$ represent the subregions filled with a single component of the system, while the set $\{\varphi(\cdot,\, t) \in (-1,1)\}$ denotes the so-called diffuse interface, i.e., a layer of positive thickness where some mixing takes place. However, it is a priori unclear how \eqref{eq:CH1} should ensure this property. In this regard, choosing a function $F$ whose derivatives blow up at $\pm 1$ can yield to a physically consistent description of phase separation. The most appropriate choice of the potential $F$ can be derived from thermodynamical arguments and is the so-called Flory--Huggins (or singular) double-well potential (see \cite{Flory42, Huggins41}), i.e.,
    \[
    F_{\text{log}}:[-1,1] \to \mathbb R, \qquad F_{\text{log}}(s) = \dfrac{\theta}{2}\left[(1+s)\log(1+s)+(1-s)\log(1-s)\right]- \dfrac{\theta_0}{2}s^2,
    \]
    extended by continuity at the endpoints. Here, $\theta$ denotes the absolute temperature of the system (which is assumed to be constant) and $\theta_0$, that satisfies $0 < \theta < \theta_0$, is the critical temperature under which phase separation takes place. Nonetheless, from both the modeling and the analytical viewpoints, alternative choices of $F$ are also of significant interest. A second, widely employed example for the potential $F$ is actually a regular approximation of $F_{\text{log}}$, for instance, 
    \[
    F_{\text{reg}}: \mathbb R \to \mathbb R, \qquad F_{\text{reg}}(s) = \dfrac{(s^2-1)^2}{4},
    \]
    which retains a double-well structure, but shifts the two global minimum points at $s = \pm 1$ and is globally smooth over $\mathbb R$ (i.e., no derivative of $F_{\text{reg}}$ exhibits singularities). Although regular potentials as $F_{\text{reg}}$ do not strictly enforce the physical constraint \eqref{eq:physical}, they provide a flexible and robust modeling framework. For instance, they allow for a clearer mathematical analysis and facilitate numerical approximations. Moreover, the investigation of the Cahn--Hilliard system \eqref{eq:CH1} with regular potentials is a necessary preliminary step to understand the singular case, while also retaining independent interest in several models (see, for instance, \cite[Section 1]{DPGW} for some remarks on the topic). Regardless of the specific form of the potential energy density $F$, formally, the Cahn--Hilliard system \eqref{eq:CH1} has a $H^{-1}$-gradient flow structure: indeed, it can be derived starting from the energy functional
    \begin{equation} \label{eq:localenergy}
        \mathcal{E}_\text{L}(\varphi) = \int_\OO \dfrac{1}{2}|\nabla \varphi(x)|^2 + F(\varphi(x)) \: \d x
    \end{equation}
    and computing the chemical potential $\mu$ as its first variational derivative. In this work, we shall actually be interested in a variant of system \eqref{eq:CH1}, which accounts for long-range interactions. Indeed, the energy functional \eqref{eq:localenergy} only takes into consideration contributions due to local interactions (so that system \eqref{eq:CH1} is also named the local Cahn--Hilliard equation), penalizing the presence of spatial variations across the interface. At the end of the last century, G. Giacomin and J. L. Lebowitz (see \cite{GL1}, \cite{GL2} and \cite{GL3}) addressed the fundamental problem of deriving a model of phase separation through a stochastic hydrodynamic limit. The resulting model is similar to \eqref{eq:CH1} in structure, but captures the energy contributions due to nonlocal interactions. Retaining the notation in \eqref{eq:CH1}, the so-called nonlocal Cahn--Hilliard equation reads
    \begin{equation} \label{eq:CH2}
    \begin{cases}
        \partial_t \varphi - \Delta \mu =0 & \quad \text{in } (0,T) \times \OO, \\
         \mu = (K*1)\varphi-K * \varphi + F'(\varphi)& \quad \text{in } (0,T) \times \OO, \\
         \partial_{\bn} \mu = 0 & \quad \text{on } (0,T) \times \partial\OO, \\
         \varphi(0) = \varphi_0 & \quad \text{in } \OO,
    \end{cases}
    \end{equation}
    where $K: \mathbb R^d \to \mathbb R$ is a suitable interaction kernel and $F: \mathbb R \to \mathbb R$ is still a potential energy density function. The peculiarity of this model is that an analogous derivation for the local equation \eqref{eq:CH1} is, to our knowledge, still not available to this day. Of course, the energy functional \eqref{eq:energy} is not the appropriate one when considering the nonlocal variant \eqref{eq:CH2}: indeed, the gradient squared term has to be replaced by a suitable nonlocal counterpart, that is
    \begin{equation} \label{eq:energy}
        \mathcal{E}_\text{NL}(\varphi) = \dfrac{1}{4}\int_\OO\int_\OO K(x-y)(\varphi(x)-\varphi(y))^2 \: \d 
x \, \d y + \int_\OO F(\varphi(x)) \: \d x.
    \end{equation}
    Equation \eqref{eq:CH2} has been studied in a number of recent contributions, both as a standalone problem or coupled with other equations from mathematical physics (e.g., the Navier--Stokes equations), refer for instance to \cite{GGG, FG2, FG3, FGG} and the references therein.\\
    On the modeling side, however, it should be noted that the purely deterministic description of phase separation (either local or nonlocal) fails in capturing thermodynamical fluctuations or environmental noise, that are characterized by unpredictable oscillations at the microscopic level and play are particularly relevant during spinodal decomoposition. The first effort towards a refinement of the local Cahn--Hilliard equation in this direction was proposed in 1970 (see \cite{cook}), where a stochastic version of equation \eqref{eq:CH1} was introduced (known as the stochastic Cahn--Hilliard equation or the Cahn--Hilliard--Cook equation). Since then, this stochastic counterpart has been analyzed in a number of works (see, for instance, \cite{daprato-deb,deb-zamb, deb-goud, EM1991,goud} and, more recently, \cite{scar-SCH, scar-SVCH,scarpa21, DPGS2}, see also \cite{DPGS} for a coupled model). All of the previously mentioned works, however, deal with the local version of the model. Indeed, the only references about nonlocal diffuse interface models concern more complicated models including a drift, either under additive noise (see \cite{corn}), or without a direct stochastic forcing (see \cite{medjo1}, for instance).
    In this work, we shall consider the stochastic version of the nonlocal Cahn--Hilliard equation, that is, 
    \begin{equation} \label{eq:strongNCH}
    \begin{cases}
        \d \varphi - \Delta \mu \: \d t = G(\varphi) \: \d W & \quad \text{in } (0,T) \times \OO, \\
         \mu = (K*1)\varphi - K * \varphi + F'(\varphi)& \quad \text{in } (0,T) \times \OO, \\
         \partial_{\bn} \mu = 0 & \quad \text{on } (0,T) \times \partial\OO, \\
         \varphi(0) = \varphi_0 & \quad \text{in } \OO,
    \end{cases}
    \end{equation}
    where $F: \erre \to \erre$ is a given regular potential density function, $K: \erre^d \to \erre$ is a given interaction kernel and $G$ is a suitable stochastically integrable diffusion process of multiplicative type.
    After the necessary preliminaries, we employ a Yosida regularization scheme and a stochastic compactness argument to show the existence of martingale solutions. If the noise enables conservation of mass (see also \cite{DPGS2} for further observations on the topic), then, by establishing pathwise uniqueness, we are able to show that even probabilistically-strong solutions exist. Finally, we estalish the nonlocal-to-local asymptotics of the problem, showing that along a suitable sequence of nonlocal kernels, solutions to the nonlocal problem \eqref{eq:strongNCH} converge to solutions of the corresponding local problem. The content of this work is structured as follows. In Section \ref{sec:main}, we introduce the notation used throughout the work and state the main results. Sections \ref{sec:proof1} and \ref{sec:proof2} are devoted to the proof of existence of a martingale solutions and of a probabilistically-strong solution, respectively. Finally, in Sections \ref{sec:convergence_unique} and \ref{sec:proof3} we investigate the nonlocal-to-local asymptotics of the model.
    \section{Preliminaries and main results} \label{sec:main}
    \subsection{Functional setting and notation}
    The notation and both the analytical and stochastic frameworks that are needed to analyze the problem are illustrated in the following paragraphs.
    \paragraph{\textit{Stochastic setting.}} Let $(\Omega,\cF,(\cF_t)_{t\in[0,T]},\P)$ be any filtered probability space satisfying the usual conditions (namely the filtration is saturated and right-continuous), with $T>0$ being a prescribed final time. Let us denote the progressive $\sigma$-algebra on $\Omega \times [0,T]$ by $\cP$. Let further $(E,\, \mathcal M,\, m)$ be any probability space. If $X: \Omega \to E$ and $Y: \Omega \to E$ are two random variables, then the notation
    \[
    X \laweq Y
    \]
    denotes the identity of their laws. When $E$ is also a Banach space (in these cases, we always understand that $\mathcal M$ is the Borel $\sigma$-algebra of $E$ unless otherwise stated),
	for any positive real quantity $p \geq 1$ the symbol $L^p(\Omega; E)$ denotes the set of strongly measurable $E$-valued random variables on $\Omega$ with finite moments up to order $p$. If the space $E$ is allowed to depend on time, as in classical Bochner spaces, we may denote the resulting space with $L^p_\cP(\Omega; E)$	to stress that measurability is intended with respect to the progressive
	$\sigma$-algebra $\cP$.	Throughout the work, the symbol $W$ denotes a cylindrical Wiener process on some separable, fixed Hilbert space $U$. We shall denote by $\{u_j\}_{j\in\enne} \subset U$ a fixed arbitrary orthonormal system for $U$. In order to avoid any ambiguities, we precise the rigorous interpretation of the stochastic perturbation appearing in \eqref{eq:strongNCH} as an It\^{o} integral.
		Being a cylindrical process on $U$, the Wiener process $W$ admits the representation
		\begin{equation} \label{eq:representation}
			W = \sum_{k=0}^{+\infty} \beta_k u_k,
		\end{equation}
		where $\{\beta_k\}_{k \in \enne}$ is a family of real and independent Brownian motions. However, the series \eqref{eq:representation} can not be expected to always converge in $U$. In general, there exists some larger Hilbert space $U_0$ such that $U \embed U_0$ with Hilbert-Schmidt embedding $\iota$ and such that we can identify $W$ as a $Q^0$-Wiener process on $U_0$, for some trace-class operator $Q^0$ (see \cite[Subsection 2.5.1]{LiuRo}). Actually, it holds that $Q^0 = \iota \circ \iota^*$. In the following, we may implicitly assume this extension by simply saying that $W$ is a cylindrical process on $U$. This holds also for stochastic integration with respect to $W$. Indeed, the symbol
		\[
		\int_0^\cdot B(s)\,\d W(s) := \int_0^\cdot B(s) \circ \iota^{-1}(s)\,\d W(s),
		\]
		for every process $B \in L^2_\cP(\Omega;L^2(0,T;\cL^2(U,K)))$, where $K$ is any real Hilbert space and the space $\cL^2(U,K)$ is defined later. It is well known that such a  definition is well posed and does not depend on the choice of $U_0$ or $\iota$ (see \cite[Subsection 2.5.2]{LiuRo}).
    \paragraph{\textit{Analytical setting}.} \label{subsub:2}
	Let $E$ be any Banach space. The symbol $\boldsymbol E$ denotes the product space $E^d$ (or even $E^{d\times d}$, if no ambiguities arise), while its topological dual is denoted by $E^*$ and the corresponding duality pairing is denoted by $\ip{\cdot}{\cdot}_{E^*,E}$.
	If $E$ is a Hilbert space, then the scalar product of $E$ is denoted by $(\cdot,\cdot)_E$. Given two Banach spaces $E_1$ and $E_2$ and an extended real number $s \in [1,+\infty]$, the symbol $L^s(E_1; E_2)$ indicates the usual spaces of strongly measurable, Bochner-integrable functions
	defined on the Banach space $E_1$ and with values in the Banach spaces $E_2$. If $E_2$ is omitted, it is understood that $E_2 = \mathbb{R}$.
	When $E_1$ and $E_2$ are also Hilbert spaces, we denote the space of Hilbert-Schmidt operators
	from $E_1$ to $E_2$ by the symbol $\cL^2(E_1,\,E_2)$, and we endow it with its canonical norm $\norm{\,\cdot\,}_{\cL^2(E_1,\,E_2)}$.
    For all $s\in(1,+\infty)$ and for every separable and reflexive Banach space $E$, we also define
	\[
	L^s_w(\Omega; L^\infty(0,T; E^*)):=
	\left\{v:\Omega\to L^\infty(0,T; E^*) \text{ weakly*-measurable and }
	\norm{v}_{L^\infty(0,T; E^*)}\in L^s(\Omega)
	\right\}\,,
	\]
	which yields by
	\cite[Theorem 8.20.3]{edwards} the identification
	\[
	L^s_w(\Omega; L^\infty(0,T; E^*))=
	\left(L^{\frac{s}{s-1}}(\Omega; L^1(0,T; E))\right)^*\,.
	\]
	Fixed a smooth and bounded domain $\OO \subset \mathbb R^d$, with $d \in \{1,2,3\}$, we denote by $W^{s,p}(\OO)$ the classical Sobolev spaces of order $s \in \mathbb R$ and $p \in [1,+\infty]$, and we denote by $\norm{\,\cdot\,}_{W^{s,p}(\OO))}$ their canonical norms (with the usual understanding that $W^{0,p}(\OO) = L^p(\OO)$). Accordingly, we define the Hilbert space $H^s(\OO):=W^{s,2}(\OO)$ for all $s\in\erre$,
	endowed with its canonical norm $\norm{\,\cdot\,}_{H^s(\OO)}$. Moreover, we set
	\[
	H:=L^2(\OO)\,, \qquad V:=H^1(\OO)\,
	\]
	endowed with their standard norms $\norm{\cdot}_H$,
	$\norm{\cdot}_{V}$, respectively. As usual, we identify the Hilbert space $H$ with its dual through
	the corresponding Riesz isomorphism, so that we have the variational structure
	\[
	V\embed H \embed V^*,
	\]
	with dense and compact embeddings (both in the cases $d = 2$ and $d = 3$). Finally, we set the zero-mean spaces
	\[
	H_0 := \left\{ u \in H : \overline{u} := \dfrac{1}{|\OO|}\int_{\OO} u \: \d x = 0\right\}, \qquad V_0 :=V \cap H_0.
	\]
	The space $H_0$ is endowed with the structure induced by $H$, hence we still carry over the same notation. Instead, owing to the Poincaré inequality, we set the $H^1$-seminorm structure on $V_0$, namely
	\[
	( u,\,v)_{V_0} := (\nabla u, \,\nabla v)_{ \bH}, \qquad \|u\|_{V_0} := \|\nabla u\|_ \bH, \qquad u,\,v \in V_0.
	\]
	This generates the zero-mean variational triplet
	\[
	V_0 \embed H_0 \embed V_0^*,
	\]
	once again with compact and dense embeddings in two and three dimensions. It is well-known that the operator
    \begin{align*}
        \mathcal{A}: V_0\rightarrow V_0^*,\quad \langle\mathcal{A}u,v\rangle_{V_0} := (\nabla u, \,\nabla v)_{ \bH}\quad \text{ for all }v\in V_0
    \end{align*}
    is a linear isomorphism. The operator $\mathcal{A}$ is then called the negative Neumann Laplace operator. We denote its inverse, which is a bounded linear operator, by
    \begin{align*}
        \mathcal{N} := \mathcal{A}^{-1}: V_0^*\rightarrow V_0.
    \end{align*}
	\subsection{Structural assumptions}
    The following assumptions will be needed throughout the whole work. Here and in the following, we refer to the stochastic nonlocal Cahn--Hilliard problem \eqref{eq:strongNCH}.
    \begin{enumerate}[label = \textbf{(A\arabic*)}, ref = \textbf{(A\arabic*)}] \itemsep0.3em
    \item \label{hyp:K} The nonlocal kernel $K \in W^{1,1}_{\text{loc}}(\mathbb{R}^d)$ satisfies $K(x) = K(-x)$ for all $x \in \mathbb{R}^d$. Throughout the entire work, we stress once and for all that the norms of $K$ shall be computed over a ball $B_R \supset \review{\overline{\OO}}$ with a sufficiently large radius $R > 0$. Moreover, we assume that
    \[a(x) := (K*1)(x) = \int_\OO K(x-y) \: \d y \geq 0\]
    for almost every $x\in \OO$.
	\item\label{hyp:potential} The potential $F:\erre \to \erre$ is such that $F \in C^2(\mathbb R)$ and satisfies the following properties:
    \begin{enumerate}[(i)]
      \item \label{hyp:sign} it holds $F'(0) = 0$ and $F(r) \geq 0$ for all $r \in \erre$,
      \item \label{hyp:growth} there exists a constant $C_F > 0$ such that
      \[
      |F'(r)|+|F''(r)| \leq C_F(1 + F(r))
      \]
      for all $r \in \erre$,
      \item \label{hyp:coercivity} \review{the second derivative $F''$ is bounded from below and} there exists a constant $C_0>0$ such that 
      \begin{align*}
          \review{\inf_{s \in \mathbb R} F^{\prime\prime}(s)+ \underset{x \in \OO}{\operatorname{ess\inf}} \: a(x)} \geq C_0.
      \end{align*}
  \end{enumerate}
		\item\label{hyp:G} The measurable operator $G: H \to \cL^2(U,H)$ satisfies the following properties:
        \begin{enumerate}[(i)]
            \item there exists a constant $L_G > 0$ such that
            \[
            \|G(\psi_1)-G(\psi_2)\|_{\cL^2(U,H)} \leq L_G\|\psi_1-\psi_2\|_H \qquad \forall \: \psi_1,\,\psi_2 \in H,
            \]
            i.e., $G$ is $L_G$-Lipschitz-continuous;
            \item \label{hyp:Gii} there exists a constant, that we still denote with $L_G > 0$ without loss of generality, such that
            \review{
            \[
            \sum_{k=0}^{+\infty} \|G(\psi)[u_k]\|_{L^\infty(\OO)}^2 \leq L_G
            \]
            }
            for all $\psi \in H$;
            \item there exists a constant, that we still denote with $L_G > 0$ without loss of generality, such that
            \[
            \|G(\psi)\|_{\cL^2(U,V)} \leq L_G\left( 1 + \|\psi\|_V\right) \qquad \forall \: \psi \in V,
            \]
            i.e., $G$ is linearly bounded in $V$.
        \end{enumerate}
    \end{enumerate}
    \begin{remark} \label{rem:assumptions}
        Let us briefly comment further on the previous assumptions. First of all, observe that for a sufficiently large $R > 0$, so that $\review{\OO - \OO := \left\{ x-y: x \in \OO, \, y \in \OO \right\}} \subset B_R$, we have
        \[
        \review{
        |a(x)| =\left| \int_{\OO} K(x-y) \: \d y \right| = \left| \int_{x-\OO} K(z) \: \d z \right| \leq \|K\|_{L^1(\OO - \OO)} \leq \|K\|_{L^1(B_R)}
        }
        \]
        and therefore $a \in L^\infty(\OO)$, owing to Assumption \ref{hyp:K}. Assumptions \ref{hyp:potential}-\ref{hyp:sign} and \ref{hyp:growth} allow for a wide class of regular potentials encompassing the usual quartic $F_{\text{reg}}$, as seen in Section \ref{sec:intro}, but also including functions of exponential growth. However, Assumption \ref{hyp:potential}-\ref{hyp:coercivity} poses a compatibility condition between the nonlocal kernel $K$ and the potential $F$, in terms of a coercivity condition. Under Assumption \ref{hyp:K}, this is trivially satisfied whenever $F$ is strongly convex, but is in general crucial in the nonlocal case, as no control on $V$-norms is otherwise possible. Indeed, this hypothesis has two important consequences. \review{First of all, we clearly have
        \[
        F''(s) + a(x) \geq  \inf_{s \in \mathbb R} F''(s) + \underset{x \in \OO}{\operatorname{ess\inf}} \: a(x) \geq C_0
        \]
        for all $s \in \mathbb R$ and almost all $x \in \OO$. Observe that, in general, the above may hold even if $F$ is not convex. Secondly,} it is clear that $F$ can be written as
        \[
        F(s) = \Psi(s) + R(s) \qquad s \in \mathbb R,
        \]
        where $\Psi$ is a strongly convex function and $R$ is a concave second-order perturbation. Letting indeed
        \[
        \inf_{s \in \mathbb R} F''(s) = -\alpha,
        \]
        and assuming $\alpha \geq 0$ (otherwise one can take $F = \Psi$ and $R \equiv 0$) one can set
        \[
        R(s) = -\dfrac{\gamma}{2} s^2
        \]
        for some $\gamma > \alpha$ and then set $\Psi(s) = F(s) - R(s)$ for all $s \in \mathbb R$. Observe further that the function $a$ defined in Assumption \ref{hyp:K}, in order to be compatible with the potential $F$, must satisfy
        \[
        a(x) \geq C_0 + \alpha
        \]
        for almost all $x \in \OO$. Finally, Assumption \ref{hyp:G} is often assumed in the context of stochastic diffuse interface models (see, for instance, \cite{scarpa21, DPGS, DPGS2, feir-petc}), also in the singular case. \review{In most instances, the noise is assumed to be of superposition type, i.e., for every $k \in \mathbb N$ and $\psi \in H$
        \[
        G(\psi)[u_k] = g_k(\psi)
        \]
        with $\{g_k\}_{k \in \mathbb N} \subset W^{1,\infty}(\mathbb R)$ and such that
        \[
        \sum_{k=0}^{+\infty}\|g_k\|_{W^{1,\infty}(\mathbb R)}^2 \leq C.
        \]
        When pathwise mass conservation is required, it is also necessary to require that the range of $G$ lies in $\cL^2(U,H_0)$ instead, so for example one can define
        \[
        G(\psi)[u_k] = g_k(\psi) - \overline{g_k(\psi)}
        \]
        for every $k \in \mathbb N$ and $\psi \in H$, or express $G(\psi)[u_k]$ as the divergence of suitably regular vector fields (see also \cite{DPGS2} for further details).}
        \end{remark}
        \subsection{Main results} Before stating the main results of the work, let us first precise the concept of martingale and probabilistically-strong solutions for system \eqref{eq:strongNCH}, starting with the former.
    \begin{defin}
        \label{def:mart_sol}
		Let $p \geq 2$ and let $\varphi_0$ satisfy
  \[
  \varphi_0 \in L^p(\Omega, \mathscr{F}_0; H), \qquad F(\varphi_0) \in  L^{\frac p2}(\Omega, \mathscr{F}_0; L^1(\OO)).
  \]
		A martingale solution to the stochastic nonlocal Cahn--Hilliard problem \eqref{eq:strongNCH} with respect to the initial datum $\varphi_0$ is a family
        \[
        \left\{ \left( \hom, \hF, (\hF)_{t \in [0,T]}, \hP\right), \hW, \hphi \right\}
        \]
        such that $\left( \hom, \hF, (\hF)_{t \in [0,T]}, \hP\right)$ is a filtered probability space satisfying the usual conditions, $\hW$ is a cylindrical Wiener process on a separable Hilbert space $U$, and $\hphi$ is an $\hF_t$-adapted process such that
		\begin{align*}
			&\hphi \in L^p_w(\hom; L^\infty(0,T;H)) \cap L^p_\cP(\hom;L^2(0,T;V)) \\
			&\hmu:=a\hphi-K*\hphi+F'(\hphi) \in L^{\frac p2}_\cP(\hom; L^2(0,T; V))\,,\\
			&\hphi(0) \laweq \varphi_0,
		\end{align*}
		and
		\begin{equation*} 
		  (\hphi(t),\psi)_H +
			\int_0^t\!\int_\OO \nabla  \hmu(s)\cdot \nabla \psi\,\review{\d x\,}\d s		= ( \hphi(0),\,\psi)_{H} +
			\left(\int_0^t G( \hphi(s))\,\d \hW(s),\, \psi\right)_{H}
        \end{equation*}
		for all $t \in [0,T]$, for every $\psi\in V_1$ and $\hP$-almost surely.
    \end{defin}
	\begin{defin}
        \label{def:strong_sol}
		Let $p \geq 2$ and let $\varphi_0$ satisfy
  \[
  \varphi_0 \in L^p(\Omega, \mathscr{F}_0; H), \qquad F(\varphi_0) \in  L^{\frac p2}(\Omega, \mathscr{F}_0; L^1(\OO)).
  \]
		Given a stochastic basis
        \[
        \left( \Omega,\, \cF,\,(\cF_t)_{t \in [0,T]},\, \P \right)
        \]
        satisfying the usual conditions and a cylindrical Wiener process $W$ on a separable Hilbert space $U$, a probabilistically-strong solution to the stochastic nonlocal Cahn--Hilliard problem \eqref{eq:strongNCH} with respect to the initial datum $\varphi_0$ is an $\cF_t$-adapted process $\varphi$
        such that
		\begin{align*}
			&\varphi \in  L^p_w(\Omega; L^\infty(0,T;H)) \cap L^p_\cP(\Omega;L^2(0,T;V)) \\
			&\mu:=a\varphi-K*\varphi+F'(\varphi) \in L^{\frac p2}_\cP(\Omega; L^2(0,T; V))\,,\\
			&\varphi(0) = \varphi_0,
		\end{align*}
		and
		\begin{equation*} 
		  (\varphi(t),\psi)_H +
			\int_0^t\!\int_\OO \nabla  \mu(s)\cdot \nabla \psi\,\review{\d x}\,\d s		= ( \varphi(0),\,\psi)_{H} +
			\left(\int_0^t G( \varphi(s))\,\d W(s),\, \psi\right)_{H}
        \end{equation*}
		for all $t \in [0,T]$, for every $\psi\in V$ and $\P$-almost surely.
    \end{defin} \noindent
    We are now in a position to state the main results of this work. \review{The first result concerns the existence of martingale solutions and will be proven in Section \ref{sec:proof1}.}
    \begin{thm} \label{thm:existenceM}
         Let Assumptions \ref{hyp:K}-\ref{hyp:G} hold. Let $p > 2$ and let $\varphi_0$ satisfy
  \[
  \varphi_0 \in L^p(\Omega, \mathscr{F}_0; H), \qquad F(\varphi_0) \in  L^{\frac p2}(\Omega, \mathscr{F}_0; L^1(\OO)).
  \]Then, there exists at least one martingale solution to the stochastic nonlocal Cahn--Hilliard equation \eqref{eq:strongNCH}.
    \end{thm} \noindent
    Under stronger assumptions on the stochastic diffusion $G$, we are also able to prove a continuous dependence estimate on the full paths, yielding pathwise uniqueness of martingale solutions.
    \begin{enumerate}[label = \textbf{(A\arabic*)}, ref = \textbf{(A\arabic*)}, resume]
        \item \label{hyp:additional} Fixed a complete orthonormal system $\{u_k\}_{k \in \mathbb{N}}\subset U$, the stochastic diffusion $G:H\to \cL^2(U, H)$ satisfies
        \[
        \overline{G(\psi)[u_k]} = 0
        \]
        for all $k \in \mathbb N$, and there exists a positive constant, that we still denote with $L_G > 0$, such that
            \[
        \|G(\psi_1)-G(\psi_2)\|_{\cL^2(U,V^*)} \leq L_G\|\psi_1-\psi_2\|_{V^*}
            \]
            for all $\psi_1,\,\psi_2 \in H$, i.e., $G$ is $L_G$-Lipschitz-continuous with respect to the $V^*$-topology. 
    \end{enumerate}
    \begin{remark}
        As we shall better see later, the assumption of Lipschitz continuity in the $V^*$-topology could be replaced by a smallness assumption on the Lipschitz constant with respect to the $H$-topology. Nonetheless, in order to avoid ambiguities in the following, we stick with this more classical assumption. 
    \end{remark}
    \begin{remark}
        Assumption \ref{hyp:additional} implies that the dynamics of problem is mass-conservative in a pathwise sense, i.e., that any probabilistically-strong solution satisfies
        \[
        \overline{\varphi(t)} = \overline{\varphi_0}
        \]
        for all $t > 0$, $\P$-almost surely. This follows by taking $\psi = 1$ in the variational formulation of the problem and recovers a key feature of the deterministic version of the system. Moreover, Assumption \ref{hyp:additional} accounts for the fact that a continuous dependence estimate in the $H$-norm for the Cahn--Hilliard equation (both local and nonlocal) is generally out of reach due to the nonlinear potential term.
    \end{remark}
    \begin{thm} \label{thm:uniqueness}
        Let Assumptions \ref{hyp:K}-\ref{hyp:additional} hold. Let $p > 2$ and, for $i \in \{1, 2\}$, let $\varphi_{0i}$ comply with Theorem \ref{thm:existenceM}. \review{Assume that $\hphi_1$ and $\hphi_2$ are defined on the same filtered probability space $(\hom, \hF, (\hF_t)_{t \in [0,T]},\hP)$ and solve \eqref{eq:strongNCH} driven by the same $\hF_t$-adapted Wiener process $\hW$ and with respect to the initial data $\varphi_{01}$ and $\varphi_{02}$, respectively.} If $\overline{\varphi_{01}} = \overline{\varphi_{02}}$ $\P$-almost surely, then, the continuous dependence estimate
        \begin{equation*}
        \|\hphi_1-\hphi_2\|_{L^p(\hom; L^\infty(0,T;V^*))} + \|\hphi_1-\hphi_2\|_{L^p(\hom; L^2(0,T;H))} \leq C\|\hphi_{01}-\hphi_{02}\|_{L^p(\hom;V^*)}
    \end{equation*}
    holds, \review{where $\hphi_{0i} := \hphi_{i}(0) \laweq \varphi_{0i}$ as $H$-valued random variables for $i \in \{1,2\}$.} In particular, the martingale solution to \eqref{eq:strongNCH} is pathwise unique.
    \end{thm}
    In any case, \review{as a consequence of the previous result and the Yamada--Watanabe theorem (see, for instance, \cite[Theorem E.0.8]{LiuRo}),} we have the straightforward corollary
    \begin{cor} \label{cor:existenceS}
       Let Assumptions \ref{hyp:K}-\ref{hyp:additional} hold. Let $p > 2$ and let $\varphi_0$ satisfy
  \[
  \varphi_0 \in L^p(\Omega, \mathscr{F}_0; H), \qquad F(\varphi_0) \in  L^{\frac p2}(\Omega, \mathscr{F}_0; L^1(\OO)).
  \]Then, there exists a unique probabilistically-strong solution to the stochastic nonlocal Cahn--Hilliard equation \eqref{eq:strongNCH}.
    \end{cor} \noindent
    \review{The proofs of Theorem \ref{thm:uniqueness} and Corollary \ref{cor:existenceS} will be discussed in Section \ref{sec:proof2}. }
    The last results of the work concerns the nonlocal-to-local asymptotics for the stochastic Cahn--Hilliard equation, i.e., establishing whether and in which sense the solutions of the nonlocal equations converge to the ones of the local problem, that reads
    \begin{equation} \label{eq:localCHstrong}
        \begin{cases} 
        \d \varphi - \Delta \mu \: \d t = G(\varphi) \: \d W & \quad \text{in } (0,T) \times \OO, \\
         \mu = -\Delta\varphi+ F'(\varphi)& \quad \text{in } (0,T) \times \OO, \\
         \partial_{\bn} \mu = \partial_\bn \varphi=0 & \quad \text{on } (0,T) \times \partial\OO, \\
         \varphi(0) = \varphi_0 & \quad \text{in } \OO.
    \end{cases}
    \end{equation}
    More precisely, we shall make use of a specific sequence of nonlocal kernels. In the following, we denote by $\{\rho_\varepsilon: \mathbb R \to [0,+\infty)\}_{\varepsilon>0}$ a family of nonnegative mollifiers satisfying the following properties:
    \begin{enumerate}[(i)]
        \item the function $\rho_\varepsilon$ is integrable for all $\eps > 0$;
        \item the function $\rho_\varepsilon$ is even, i.e., $\rho_\eps(r) = \rho_\eps(-r)$ for all $r \in \erre$ and $\varepsilon >0$;
        \item the relation
        \[
        \int_{0}^{+\infty}\rho_\varepsilon(r)\:r^{n-1}\:\text{d}r = \frac{2}{C_n}
        \]
        holds for all integers $n \geq 2$, where
        \[
        C_n := \int_{\mathbb{S}^{n-1}}|e_1\cdot\sigma|^2\:\text{d}\mathcal{H}^{n-1}(\sigma),
        \]
        with $e_1$ being the first element of the canonical basis in $\mathbb R^{n}$, $\mathbb S^{n-1}$ denoting the unit sphere of $\mathbb R^n$ and $\mathcal H^{n-1}$ being the $(n-1)$-dimensional Hausdorff measure;
        \item it holds 
        \[
        \lim_{\eps \to 0^+}\int_{\delta}^\infty\rho_\varepsilon(r)\:r^{n-1}\:\text{d}r = 0
        \]
        for all $\delta > 0$.
    \end{enumerate}
    The nonlocal kernels are of the following form.
    \begin{enumerate}[label = \textbf{(A\arabic*)}, ref = \textbf{(A\arabic*)}, resume] 
    \item \label{hyp:asymptotics} For any $\varepsilon > 0$, we define a nonlocal kernel $K_\varepsilon: \mathbb{R}^n\rightarrow[0,+\infty)$ by
    \[
        K_\varepsilon(x) = \frac{\rho_\varepsilon(x)}{|x|^2}
    \]
    for any $x \in \mathbb R^n$.
    \end{enumerate}
    \begin{remark}
        Assumption \ref{hyp:asymptotics} is not new in the context of nonlocal-to-local asymptotics for the Cahn--Hilliard equation: see for instance \cite{Davoli1, Davoli2, abels2023strong} for further motivation. Observe that by \cite[Lemma 3.1]{Davoli1}, we have $K_\eps \in W^{1,1}(\mathbb R^n)$.
    \end{remark}
    For the sake of convenience, for any $\eps>0$ and any sufficiently regular function $u:\OO\rightarrow\R$, we define the linear nonlocal operator $\mathcal{L}_\eps u$ as
    \begin{align}\label{Def:NonlocOp}
        \mathcal{L}_\eps u(x) := \int_\OO K_\eps(x-y)(u(x)-u(y)) \: \d y
    \end{align}
    for almost all $x\in\OO$. In the limit $\eps\rightarrow0^+$, it is well-known that the nonlocal operator as in \eqref{Def:NonlocOp} converges to the negative Laplacian, cf., e.g., \cite{abels2023strong}.
    In this perspective, we analyze the nonlocal problems
    \begin{equation} \label{eq:NCH:eps}
    \begin{cases}
        \d \varphi_\eps - \Delta \mu_\eps \: \d t = G(\varphi_\eps) \: \d W & \quad \text{in } (0,T) \times \OO, \\
         \mu_\eps = \mathcal{L}_\eps\varphi_\eps + F'(\varphi_\eps)& \quad \text{in } (0,T) \times \OO, \\
         \partial_{\bn} \mu_\eps = 0 & \quad \text{on } (0,T) \times \partial\OO, \\
         \varphi_\eps(0) = \varphi_{0,\eps} & \quad \text{in } \OO.
    \end{cases}
    \end{equation}
    Since $K_\varepsilon$ is nonnegative for all $\eps > 0$, it is immediate to check that $K_\varepsilon$ satisfies Assumption \ref{hyp:K} for all $\varepsilon>0$. Therefore, under Assumptions \ref{hyp:K}-\ref{hyp:G} and \ref{hyp:asymptotics}, we infer from Theorem \ref{thm:existenceM} that \eqref{eq:NCH:eps} admits at least one martingale solution. On the remaining Assumption \ref{hyp:additional}, we observe the following:
    \begin{enumerate} \itemsep0.5em
        \item under Assumption \ref{hyp:additional}, Corollary \ref{cor:existenceS} gives unique probabilistically-strong solutions for the nonlocal problems \eqref{eq:NCH:eps}. The same actually holds for the local problem \eqref{eq:localCHstrong} as well, as addressed in \cite{scar-SCH}. In this setting, all the needed solutions can be framed with respect to the same stochastic basis. In particular, it makes sense to compute differences between solutions and estimate a rate of nonlocal-to-local convergence, extending the work previously done in \cite{abels2023strong};
        \item if we do not assume Assumption \ref{hyp:additional}, we lack uniqueness both on the nonlocal and on the local problems, and we can only work with martingale solutions. A meaningful convergence result can then only be achieved in law and for subsequences by means of a compactness argument similar to the one that is necessary to prove existence of martingale solutions.
    \end{enumerate}
    Therefore, in order to give a precise nonlocal-to-local convergence result, Assumption \ref{hyp:additional} can not be dropped. In general, we expect that the second case could be improved if a uniqueness result is available at least for the limit problem, for instance by adopting the technique introduced by Gy\"{o}ngy and Krylov (see \cite[Lemma 1.1]{gyo-kry}, see also \cite{DPGS2} for an application in the context of diffuse interface models). 
    \begin{remark}
        Actually, it turns out that the structure of the kernels $K_\varepsilon$ would in principle enable some relaxation of Assumption \ref{hyp:additional}, as the Lipschitz continuity in the $V^*$-topology turns out not to be necessary to investigate the nonlocal-to-local asymptotics of the system. As will be clear later, this is due to the validity of the so-called nonlocal Ehrling interpolation inequality (see \cite[Lemma 3.4]{Davoli1}). However, in the present setting Assumption \ref{hyp:additional} is still needed as-is to grant pathwise uniqueness to the limit local problem (and, coincidentally, to the nonlocal problems as well).
    \end{remark} \noindent
    Nonetheless, both the statements of the final results of this work shall need an additional assumption. In light of the results in \cite{abels2023strong}, establishing a rate of convergence strongly depends on some degree of higher-order regularity of the limit solution, which has been proven in \cite{scarpa21} only in the case of sub-quartic potentials (in three dimensions). Therefore, we shall need the following assumption.
    \begin{enumerate}[label = \textbf{(A\arabic*)}, ref = \textbf{(A\arabic*)}, start = 6] 
    \item \label{hyp:quarticpotential} There exists a positive constant, that without loss of generality we still denote with $C_F >0$, such that
    \[
    F''(s) \leq C_F(1 + |s|^{q_F}) \qquad \forall \: s \in \mathbb R,
    \]
    where $q_F \geq 2$ if $d = 2$ and $q_F = 2$ if $d = 3$.
    \end{enumerate}
    \review{Under the additional Assumption \ref{hyp:quarticpotential}, it has been shown that solutions of the local problem enjoy more regularity. More precisely, confining ourselves to what is needed in this work, the following lemma holds.
    \begin{lem}[{\cite[Theorem 2.2]{scarpa21}} and {\cite[Theorem 2.10]{scar-SCH}}]\label{lem:higherreg}
        Let Assumptions \ref{hyp:K}-\ref{hyp:G} and \ref{hyp:quarticpotential} hold. Let $p \geq 4$ and let $\varphi_0 \in L^p(\Omega, \cF_0; V)$ satisfy $F(\varphi_0) \in L^\frac p2(\Omega, \cF_0; L^1(\OO))$. Then, the local problem \eqref{eq:localCHstrong} admits at least a martingale solution $\{ ( \hom, \hF, (\hF_t)_{t \in [0,T]}, \hP), \hW, \hphi \}$ such that
        \[
        \hphi \in L^\frac p4(\hom; L^2(0;T;H^3(\OO))).
        \]
        Furthermore, if also Assumption \ref{hyp:additional} holds, then the solution is unique and probabilistically-strong.
    \end{lem} \noindent
    Having Lemma \ref{lem:higherreg} at disposal, we can use the results of \cite{abels2023strong} to derive concrete rates of convergence. This will be done in Section \ref{sec:convergence_unique}. In particular, we show the following result.}
    \begin{thm} \label{thm:convergence_unique}
        Let Assumptions \ref{hyp:K}-\ref{hyp:quarticpotential} hold. Let $p \geq 4$ and, for $\eps > 0$, let $\varphi_{0,\varepsilon}$ comply with Theorem \ref{thm:existenceM}, and let $\varphi_\varepsilon$ denote the probabilistically-strong solution to \eqref{eq:NCH:eps} originating from $\varphi_{0,\varepsilon}$. Let further $\varphi_0 \in L^p(\Omega, \cF_0; V)$ and $F(\varphi_0) \in L^\frac p2(\Omega, \cF_0;L^1(\OO))$ and let $\varphi$ denote the probabilistically-strong solution to \eqref{eq:localCHstrong} originating from $\varphi_{0}$. Assuming that the integral average of $\varphi_{0,\varepsilon}$ does not depend on $\varepsilon$, if
        \[
        \|\varphi_{0,\eps}-\varphi_0\|_{L^\frac p4(\Omega; V^*)} \leq C\sqrt \varepsilon,
        \]
        then the estimate
        \begin{equation*}
        \|\varphi_\varepsilon-\varphi\|_{L^\frac p4(\Omega; L^\infty(0,T;V^*))} + \|\varphi_\eps-\varphi\|_{L^\frac p4(\Omega; L^2(0,T;H))} \leq C\sqrt \varepsilon\left( 1 + \|\varphi\|_{L^\frac p4(\Omega; L^2(0,T;H^3(\OO)))} \right)
    \end{equation*}
    holds. In particular, we have $\varphi_\varepsilon \to \varphi$ strongly in $L^\frac p4(\Omega; L^\infty(0,T;V^*) \cap L^2(0,T;H))$.
    \end{thm} \noindent
    The assumption $p \geq 4$ is not restrictive, as the initial conditions are in any case tailored for convergence (in fact, they could also be taken deterministic without much loss of generality). However, the main drawback of Theorem \ref{thm:convergence_unique} is restricting the family of potentials to a subset of polynomial functions. Dropping Assumption \ref{hyp:quarticpotential} yields the need of establishing $\varepsilon$-uniform estimates in order to retrieve a solution. To this end, an additional assumption on $G$ is needed.
    \begin{enumerate}[label = \textbf{(A\arabic*)}, ref = \textbf{(A\arabic*)}, start = 7] 
    \item \label{hyp:G2} There exists a positive constant, that without loss of generality we denote with $L_G >0$, such that
    \[
    \|G(\psi)\|_{\cL^2(U,V)} \leq L_G\left(1+\|\psi\|_H \right) \qquad \forall \: \psi \in V.
    \]
    \end{enumerate}
    \begin{remark}
        Assumption \ref{hyp:G2} is a technical requirement to make use of a nonlocal Poincaré inequality requiring a higher order norm of $G$. Refinements or alternative formulations of this assumption may be possible, also taking Assumption \ref{hyp:G} into account. For example, it would be sufficient to ask that
        \[
        \operatorname{Range}(G) \subset \cL^2(U,W^{1,p}(\OO))
        \]
        for some $p > d$. However, this property would appear much more natural on the restriction of $G$ over $W^{1,p}(\OO)$, but proving additional spatial regularity (beyond $H^1(\OO)$) on the solutions of the nonlocal Cahn--Hilliard equation appears challenging.
    \end{remark}
    The last result of this work\review{, which we prove in Section \ref{sec:proof3},} reads as follows.
    \begin{thm} \label{thm:convergence}
    \review{Let Assumptions \ref{hyp:K}-\ref{hyp:asymptotics} and \ref{hyp:G2} hold and let $p > 2$. Fixed $\eps_0>0$ sufficiently small, define a family of random variables $\{\varphi_{0,\eps}\}_{\eps \in (0,\eps_0)} \subset L^p(\Omega;H)$ such that $\varphi_{0,\eps}$ satisfies the assumptions of Theorem \ref{thm:existenceM} for all $\eps \in (0,\eps_0)$ and 
    \[
    \sup_{\eps\in(0,\eps_0)}\E \left| \dfrac{1}{4}\int_\OO\int_\OO K_\eps(x-y)(\varphi_{0,\eps}(x)-\varphi_{0,\eps}(y))^2 \: \d 
x \, \d y + \int_\OO F(\varphi_{0,\eps}(x)) \: \d x \right|^\frac p2 \leq C
    \]
    for a constant $C > 0$ independent of $\eps$. Let further $\varphi_0 \in L^p(\Omega, \cF_0; V)$ and $F(\varphi_0) \in L^\frac p2(\Omega, \cF_0;L^1(\OO))$ and assume that $\varphi_{0,\eps} \to \varphi_0$ in $L^p(\Omega;H)$. For any $\eps \in (0,\eps_0)$, let $\varphi_\eps$ denote the unique probabilistically-strong solution to \eqref{eq:strongNCH} with respect to the initial datum $\varphi_{0,\eps}$. Accordingly, let $\varphi$ denote the unique probabilistically-strong solution to \eqref{eq:localCHstrong} with respect to the initial datum $\varphi_{0}$. Then, for any $q < p$, it holds 
        \[
            \varphi_{\eps} \to \varphi \quad \text{in }L^q(\Omega;L^\infty(0,T;H))
        \]
        as $\eps \to 0^+$.
    }
    \end{thm}
     \section{\review{Existence of martingale solutions}} \label{sec:proof1}
    The proof of existence is carried out via a two-stage approximation scheme and a stochastic compactness argument, borrowing some ideas from \cite[Section 3]{scarpa21}. The argument is split in several steps for the sake of clarity.
    \subsection{The approximation scheme.} 
    As anticipated, the approximation scheme is carried out in two steps. Firstly, we recast problem \eqref{eq:strongNCH} exploiting Remark \ref{rem:assumptions}, hence arriving at
    \begin{equation} \label{eq:strongNCH2}
    \begin{cases}
        \d \varphi - \Delta \mu \: \d t = G(\varphi) \: \d W & \quad \text{in } (0,T) \times \OO, \\
         \mu = a\varphi - K * \varphi + \Psi'(\varphi) - \gamma\varphi& \quad \text{in } (0,T) \times \OO, \\
         \partial_{\bn} \mu = 0 & \quad \text{on } (0,T) \times \partial\OO, \\
         \varphi(0) = \varphi_0 & \quad \text{in } \OO,
    \end{cases}
    \end{equation}
    for some $\gamma > 0$ to be defined later and such that $\Psi$ is strongly convex. This enables us to  employ a Yosida regularization scheme (see, for reference, the classical textbooks \cite{brezis, barbu-monot}). Indeed, the function $\Psi$ is a proper, convex, continuously differentiable function on $\mathbb R$ (hence, in particular, lower semicontinuous), and therefore its derivative $\Psi'$ can be identified with a maximal monotone graph on $\mathbb{R}\times\mathbb{R}$. This allows us to define for all $\lambda\in(0,1)$ the resolvent operator $J_\lambda$ and the Yosida approximation $\Psi'_\lambda$ of $\Psi'$, that is
    \begin{align*}
        &J_\lambda: \mathbb{R}\rightarrow\mathbb{R},\qquad J_\lambda(s) := (I+\lambda\Psi')^{-1}(s)\qquad 
 \forall \: s\in\mathbb{R}, \\
        &\Psi'_\lambda: \mathbb{R}\rightarrow\mathbb{R},\qquad \Psi'_\lambda(s) := \lambda^{-1}(s-J_\lambda(s))\qquad \forall \: s\in\mathbb{R},
    \end{align*}
    where $I$ denotes the identity operator on $\mathbb R$.
    Now, for every $\lambda\in (0,1)$, we define the Yosida approximation of $F$ by setting
    \begin{align}\label{eq:F_lam}
        F_\lambda: \mathbb{R} \rightarrow [0,+\infty),\qquad F_\lambda(s) := F(0) + \int_0^s\Psi'_\lambda(x)\;\mathrm{d}x - \frac{\gamma}{2}s^2\qquad \forall \: s\in \mathbb{R}.
    \end{align}
    By construction, the function $F_\lambda$ is two times continuously differentiable and
    \begin{align}
        F^\prime_\lambda(s) = \Psi'_\lambda(s) - \gamma s \qquad \forall \: s\in \mathbb{R}
    \end{align}
    yielding also, by differentiating once more, that $F'_\lambda$ is Lipschitz-continuous (since $F''_\lambda$ is bounded).
    Owing to the considerations above, we are now in a position to introduce the Yosida approximation of problem \eqref{eq:CH2}, namely
    \begin{equation} \label{eq:NCH_lambda}
    \begin{cases} 
		\d\varphi_\lambda - \Delta\mu_\lambda\,\d t =
		G(\varphi_\lambda)\,\d W
		\qquad&\text{in } (0,T)\times\OO,\\
		\mu_\lambda = a\varphi_\lambda -K*\varphi_\lambda+ F_\lambda'(\varphi_\lambda)
		\qquad&\text{in } (0,T)\times\OO,\\
		\varphi_\lambda(0)=\varphi_0 \qquad&\text{in } \OO.
    \end{cases}
\end{equation}
    In order to solve \eqref{eq:NCH_lambda}, we employ a discretization scheme. Indeed, let $\{(\ell_j ,\, e_j)\}_{j\in\enne_+} \subset \mathbb R \times D(\christoph{\mathcal{A}})$ be countably many independent eigencouples of the negative Laplacian with homogeneous Neumann boundary condition, i.e., let the functions $e_j$  solve 
    \[
	\begin{cases}
		-\Delta e_j = \ell_j e_j, & \quad \text{in } \OO,\\
		\partial_\bn e_j = 0 & \quad \text{on } \partial\OO,
	\end{cases} \qquad \forall \: j \in \enne_+
	\]
    in the strong sense. Then, up to a possible renormalization, the set $\{e_j\}_{j\in\enne_+}$ is an orthonormal system in $H$ and an orthogonal system in $V$. For $n\in\enne_+$, we now define the finite dimensional space $H_n := \mathrm{span}\{e_1,\dots,e_n\} \subset V$, endowed with the $H$-norm. 
    This enables us to introduce the approximate stochastic diffusion
    \[
    G_{n}: H_n \rightarrow \cL^2(U,H_n)
    \]
    such that 
    \begin{align*}
        G_{n}(v)[u_k] := \sum_{j=1}^n(G(v)[u_k],e_j)_He_j
    \end{align*}
    for all $k\in\enne$ and $v\in H_n$. Note that the operator \review{$G_{n}$} is well-defined: indeed, for $v\in H_n$ and $n\in\enne_+$, it holds
    \begin{equation} \label{eq:calc_G}
    \begin{split} 
        \|G_{n}(v)\|_{\cL^2(U,H_n)}^2 &=  \|G_{n}(v)\|_{\cL^2(U,H)}^2 \\ & = \sum_{k=1}^{+\infty}\sum_{j=1}^n|(G(v)[u_k],e_j)_H|^2 \\
        &\leq \sum_{k=1}^{+\infty}\sum_{j=1}^{+\infty}|(G(v)[u_k],e_j)_H|^2 \\ & = \|G(v)\|_{\cL^2(U,H)}^2
    \end{split}
    \end{equation}
    and hence $G_{n}(v)\in \cL^2(U,H_n)$ for all $v\in H_n$ and $n\in\enne_+$. Owing to Assumption \ref{hyp:G}, a similar calculation then shows that $G_{n}$ is $L_G$-Lipschitz-continuous (hence uniformly in $n$) as an operator from $H_n$ to $\cL^2(U,H_n)$ for all $n \in \mathbb N_+$. Next, we define the orthogonal projection of the initial datum onto $H_n$, namely, we set
    \begin{align*}
        \varphi_{0,n} := \sum_{j=1}^n(\varphi_0,e_j)_H\,e_j
    \end{align*}
    for $n\in \enne_+$. Now, we can formulate the discretized problem projected onto $H_n$, which reads
    \begin{equation} \label{eq:NCH_galerkin}
    \begin{cases}
		\d\varphi_{\lambda,n} - \Delta\mu_{\lambda,n}\,\d t =
		G_{n}(\varphi_{\lambda,n})\,\d W
		\qquad&\text{in } (0,T)\times\OO,\\
		\mu_{\lambda,n} = a\varphi_{\lambda,n} - K*\varphi_{\lambda,n}+ F_\lambda'(\varphi_{\lambda,n})
		\qquad&\text{in } (0,T)\times\OO,\\
        \partial_{\b n}\mu_{\lambda,n} = 0 \qquad&\text{on }(0,T)\times\partial\OO, \\
		\varphi_{\lambda,n}(0)=\varphi_{0,n} \qquad&\text{in } \OO.
    \end{cases} 
    \end{equation}
    For any fixed $\lambda>0$ and $n\in\enne_+$, we look for a solution $(\varphi_{\lambda,n},\mu_{\lambda,n})$ to the discretized problem \eqref{eq:NCH_galerkin} of the form 
    \begin{align*}
        \varphi_{\lambda,n} = \sum_{j=1}^n {\alpha}_j^{\lambda,n}e_j,\qquad \mu_{\lambda,n} = \sum_{j=1}^n {\beta}_j^{\lambda,n}e_j,
    \end{align*}
    for some suitable vector-valued processes 
    \begin{align*}
        {\alpha}^{\lambda,n}:= ({\alpha}_1^{\lambda,n},\ldots,{\alpha}_n^{\lambda,n}): \Omega \times [0,T] \rightarrow \R^n, \qquad  {\beta}^{\lambda,n}:= ({\beta}_1^{\lambda,n},\ldots,{\beta}_n^{\lambda,n}): \Omega \times [0,T] \rightarrow \R^n.
    \end{align*}
    Using this ansatz on $(\varphi_{\lambda,n},\mu_{\lambda,n})$ and choosing $e_i$ as test functions for {all} $i\in\{1,\ldots,n\}$, it follows that the variational form of \eqref{eq:NCH_galerkin} is given by 
    \begin{equation} \label{eq:var_1}
        \int_{\OO}\varphi_{\lambda,n}(t)e_i\: \d x + \int_0^t \int_{\OO}\nabla\mu_{\lambda,n}(s)\cdot\nabla e_i\: \d x\,\d s 
        = \int_{\OO}\varphi_{0,n}e_i\: \d x + \int_{\OO}\left(\int_0^t G_{n}(\varphi_{\lambda,n}(s))\: \d W(s)\right)e_i\: \d x
    \end{equation}
    for all $i \in \{1,\ldots,n\}$ and every $t\in[0,T]$, $\mathbb{P}$-almost surely, as well as 
    \begin{align} \label{eq:var_2}
        \int_{\OO}\mu_{\lambda,n}(t)e_i\: \d x = \int_{\OO}{a\varphi_{\lambda,n}(t)e_i}\:\d x - \int_{\OO}{[K*\varphi_{\lambda,n}(t)]e_i}\:\d x + \int_{\OO}F^\prime_\lambda(\varphi_{\lambda,n}(t))e_i \:\d x 
    \end{align}
    for all $i \in \{1,\ldots,n\}$ and every $t\in[0,T]$, $\mathbb{P}$-almost surely. Due to the orthogonality properties of $\{e_i\}_{i\in\enne_+}$, we deduce that the processes ${\alpha}^{\lambda,n}$ and ${\beta}^{\lambda,n}$ satisfy the following ordinary stochastic differential equations,
    \begin{equation} \label{eq:NSCH_galerkin2}
    \begin{cases}
        \displaystyle \d {\alpha}_i^{\lambda,n} + {\ell_i\beta}_i^{\lambda,n} = \left(G_{n}\left(\sum_{l=1}^n {\alpha}_l^{\lambda,n}e_l\right)\d W,\,e_i\right)_H, 
        \\
        \displaystyle {\beta}_i^{\lambda,n} = \sum_{j=1}^n{\alpha}_j^{\lambda,n}\int_\OO ae_je_i \: \d x - \sum_{j=1}^n{\alpha_j^{\lambda,n} \int_\OO [K*e_j]e_i\: \d x} + \int_\OO F^\prime_{\lambda}\left(\sum_{l=1}^n{\alpha}_l^{\lambda,n}e_l\right)e_i \: \d x, 
        \\
        \displaystyle a_i^{\lambda,n}(0) = (\varphi_0,e_i)_H 
    \end{cases} 
    \end{equation}
    for all $i \in \{1,\ldots,n\}$. Note that the stochastic integrals on the right-hand side have to be regarded, for all $i \in \{1,\ldots,n\}$, as $G_i^{n}\,\d W$, where 
    \begin{align*}
        G_i^{n}: H_n \rightarrow \cL^2(U,\R),\qquad G_i^{n}\left(\sum_{l=1}^nc_le_l\right)[u_k] := \left(G_n\left(\sum_{l=1}^nc_le_l\right)\review{[u_k]},e_i\right)_H
    \end{align*}
    for all $k\in\enne$. Since all the nonlinearities appearing in {the discretized problem \eqref{eq:NSCH_galerkin2}} are Lipschitz continuous, the theory of abstract stochastic evolution equations applies. Hence, there exists a unique solution
    \begin{align*}
        {\alpha}^{\lambda,n},\, {\beta}^{\lambda,n} \in L^{p}(\Omega;C^0([0,T];\R^n)).
    \end{align*}
    Therefore, the approximated problem \eqref{eq:NSCH_galerkin2} admits a unique solution 
    \begin{align*}
        \varphi_{\lambda,n},\, \mu_{\lambda,n}\in  L^p(\Omega;C^0([0,T];H_n))
    \end{align*}
    where $p$ is given by the moment regularity of initial conditions.
    \subsection{Uniform estimates with respect to $n$} 
    \label{ssec:unif_n}
    Next, we prove that the approximated solutions $(\varphi_{\lambda,n},\mu_{\lambda,n})$ satisfy some basic regularity estimates uniformly in $n$, {while} keeping $\lambda>0$ fixed. Here, the symbol $C$ (possibly indexed) denotes a positive constant whose special dependencies are explicitly pointed out, if necessary. In particular, throughout this subsection, the dependence on $\lambda$ is always explicited by the symbol $C_\lambda$. Before continuing, we point out a key preliminary lemma.
    \begin{lem} \label{lem:coercivity}
        Assume Assumptions \ref{hyp:K} and \ref{hyp:potential}. There exists $\gamma > 0$ such that
        \[
        F_\lambda''(s) + a(x) = \Psi''_\lambda(s) - \gamma + a(x) \geq \dfrac{C_0}{2}
        \]
        for all $s \in \mathbb R$, almost every $x \in \OO$ and all $\lambda \in (0,1]$.
    \end{lem}
    \begin{proof}
        As established in Remark \ref{rem:assumptions}, for every $s \in \mathbb R$ we have
        \[
        \Psi''(s) \geq \gamma - \alpha > 0
        \]
        and invoking, for instance, \cite[Lemma 2.10]{GGG}, we get
        \begin{equation} \label{eq:coerc1}
            \Psi''_\lambda(s) \geq \dfrac{\gamma-\alpha}{1+\gamma-\alpha},
        \end{equation}
        for all $s \in \mathbb R$ and for all $\lambda \in (0,1]$. Therefore, recalling that $a(x) \geq C_0 + \alpha$ almost everywhere in $\OO$,
        \[
        \begin{split}
            \Psi_\lambda''(s) - \gamma + a(x) & \geq \dfrac{\gamma-\alpha}{1+\gamma-\alpha} - \gamma + C_0 + \alpha \\
            & = - \dfrac{(\gamma-\alpha)^2}{1+\gamma-\alpha} + C_0
        \end{split}
        \]
        and as 
        \[
        x \mapsto \dfrac{x^2}{1+x}
        \]
        vanishes as $x \to 0^+$, we can choose $\gamma-\alpha$ small enough in such a way that the claim holds.
    \end{proof} \noindent
    \review{In all of the following, }in light of Lemma \ref{lem:coercivity}, we shall implicitly assume, without loss of generality, that $\lambda \leq 1$. Moreover, we set $\gamma > 0$ in \eqref{eq:strongNCH2} in such a way that Lemma \ref{lem:coercivity} holds. \review{The first estimate reads as follows.}
    \begin{lem} \label{lem:first_estimate}
    \review{The sequence $\{\varphi_{\lambda, n}\}_{n \in \mathbb N}$ is uniformly bounded in $L^p(\Omega;C^0([0,T];H)) \cap  L^p(\Omega;L^2(0,T;V))$ for all $p \geq 2$.}
    \end{lem}
    \begin{proof}
    We apply the It\^{o} formula for the $H$-norm of $\varphi_{\lambda,n}$ (see, for instance, \cite[Theorem 4.32]{dapratozab} or \cite[Theorem 4.2.5]{LiuRo}). This yields
    \begin{multline}
        \|\varphi_{\lambda,n}(t)\|_H^2 = \|\varphi_{0,n}\|_H^2 + 2\int_0^t(\Delta\mu_{\lambda,n}(s),\varphi_{\lambda,n}(s))_{H} + \|G_{n}(\varphi_{\lambda,n}(s))\|_{\cL^2(U,H_n)}^2\: \d s \\
        + 2\int_0^t\big(\varphi_{\lambda,n}(s),G_n(\varphi_{\lambda,n}(s))\: \d W(s)\big)_H.
    \end{multline}
    Recalling the definition of the chemical potential and the properties of convolutions, as well as appealing to Lemma \ref{lem:coercivity}, we obtain 
    \begin{align*}
        (-\Delta\mu_{\lambda,n},\varphi_{\lambda,n})_{H} 
        &= (\nabla\mu_{\lambda,n},\nabla\varphi_{\lambda,n})_{\b H} \\
        &= (a\nabla\varphi_{\lambda,n} + \varphi_{\lambda,n}\nabla a  - \nabla {K}*\varphi_{\lambda,n}+ F''_\lambda(\varphi_{\lambda,n})\nabla\varphi_{\lambda,n},\,\nabla\varphi_{\lambda,n})_{\b H}  \\    
        &\geq \dfrac{C_0}{2}\|\nabla\varphi_{\lambda,n}\|_\b H^2 - 2\|\nabla {K}\|_{\b L^1({B_R})}\|\nabla\varphi_{\lambda,n}\|_\b H\|\varphi_{\lambda,n}\|_H \\
        &\geq\frac{C_0}{4} \|\nabla\varphi_{\lambda,n}\|_\b H^2 - \dfrac{4}{C_0}{\|\nabla {K}\|_{\b L^1({B_R})}^2}\|\varphi_{\lambda,n}\|_H^2.
    \end{align*}
    Moreover, thanks to \eqref{eq:calc_G} and Assumption \ref{hyp:G} it holds
    \begin{equation*}
    \begin{split}
        \|G_{n}(\varphi_{\lambda,n})\|_{\cL^2(U,H_n)}^2 
        &\leq \|G(\varphi_{\lambda,n})\|_{\cL^2(U,H)}^2 \\
        & \leq L_G^2
    \end{split}
    \end{equation*}
   for all $n \in \mathbb N$. Combining the above estimates, we have
    \begin{multline} \label{eq:1est1}
        \|\varphi_{\lambda,n}(t)\|_H^2 + \dfrac{C_0}{2}\int_0^t\|\nabla\varphi_{\lambda,n}(s)\|_\b H^2 \: \d s \\
        \leq \|\varphi_{0,n}\|_H^2 + 4\left( \dfrac{2}{C_0}\|\nabla K\|^2_{\b L^1(B_R)} + L_G^2\right)\int_0^t\|\varphi_{\lambda,n}(s)\|_H^2 \: \d s + 2L_G^2t \\+ 2\int_0^t\big(\varphi_{\lambda,n}(s),G_n(\varphi_{\lambda,n}(s)) \:\d W(s)\big)_H
    \end{multline}
    for all $n \in \mathbb N$. Note that by linearity and $1$-Lipschitz continuity of the projection $\Pi_n: H\rightarrow H_n$, it follows
    \begin{align*}
        \|\varphi_{0,n}\|_H^2 \leq \|\varphi_{0}\|_H^2.
    \end{align*}
    Next, raising \eqref{eq:1est1} to the power $\frac p2$, then taking the supremum on the interval $[0,t]$ and $\mathbb P$-expectations, we eventually arrive at
    \begin{multline*}
        \E\supt\|\varphi_{\lambda,n}(\tau)\|_H^p + \E\left|\int_0^t\|\nabla\varphi_{\lambda,n}(\tau)\|_\b H^2\: \d\tau\right|^{\frac{p}{2}} \\
        \leq C\left(1 + \E\|\varphi_{0}\|_H^p + \left|\E\int_0^t\|\varphi_{\lambda,n}(\tau)\|_H^2 \: \d \tau\right|^\frac p2 \: \d\tau + \E\supt\left|\int_0^\tau\big(\varphi_{\lambda,n}(s),G_n(\varphi_{\lambda,n}(s)) \: \d W(s)\big)_H\right|^{\frac{p}{2}}\right),
    \end{multline*}
    where the constant $C$ only depends on the structural parameters of the problem, including $p$ and $T$. Finally, the Burkholder--Davis--Gundy and Hölder inequalities imply 
    \begin{align}
        &\E\sup_{\christoph{\tau} \in [0,t]}\left|\int_0^{\christoph{\tau}}\big(\varphi_{\lambda,n}(\christoph{s}),G_n(\varphi_{\lambda,n}(\christoph{s}))\: \d W(\christoph{s})\big)_H\right|^{\frac{p}{2}} \nonumber\\
        &\qquad \leq C\E\left|\int_0^t\|\varphi_{\lambda,n}(\christoph{s})\|_H^2\|G_n(\varphi_{\lambda,n}(\christoph{s}))\|_{\cL^2(U,H_n)}^2 \:\d\christoph{s}\right|^{\frac{p}{4}} \nonumber\\
        &\qquad \leq C\E\left|\supt\|\varphi_{\lambda,n}(\tau)\|_H^2\int_0^t\|G_n(\varphi_{\lambda,n}(\christoph{s}))\|_{\cL^2(U,H_n)}^2\: \d\christoph{s}\right|^{\frac{p}{4}} \nonumber\\
        &\qquad \leq Ct^\frac p4\E\left|\supt\|\varphi_{\lambda,n}(\tau)\|_H^2\right|^{\frac{p}{4}} \nonumber\\
        &\qquad \leq \dfrac{1}{2}\E\supt\|\varphi_{\lambda,n}(\tau)\|_H^{p} + Ct^{\frac{p}{2}} \label{eq:est_Galerkin_1}
    \end{align}
    where the constant $C$ depends on the structural parameters of the problem, including $p$.
    Hence, thanks to \eqref{eq:est_Galerkin_1} and Young's inequality, we arrive at
    \begin{align*}
        \E\supt\|\varphi_{\lambda,n}(\tau)\|_H^p + \E\left|\int_0^t\|\nabla\varphi_{\lambda,n}(\tau)\|_\b H^2 \: \d\tau\right|^{\frac{p}{2}} 
        \leq C\left(1 + \E\|\varphi_{0}\|_H^p + \E\left| \int_0^t\|\varphi_{\lambda,n}(\tau)\|_H^2 \: \d\tau \right|^\frac p2\right).
    \end{align*}
    Finally, an application of Gronwall's lemma implies that there exists some constant $C_1$, independent of $n$ and $\lambda$, such that 
    \begin{align}
        \|\varphi_{\lambda,n}\|_{L^p(\Omega;C^0([0,T];H))} + \|\varphi_{\lambda,n}\|_{L^p(\Omega;L^2(0,T;V))} \leq C_1 \label{eq:unif_est_1}
    \end{align}
    and the first estimate is proved.
    \end{proof} \noindent
    A second useful estimate is an energy inequality. In particular, we want to apply the It\^{o} formula (see \cite[Theorem 4.32]{dapratozab}) to the approximated free energy functional 
    \begin{align*}
         \mathcal{E}_\lambda: H_n \to \mathbb R, \qquad \mathcal{E}_\lambda(v) := \frac{1}{4}\int_\OO\int_\OO K(x-y)|v(x) - v(y)|^2 \: \d y \: \d x + \int_\OO F_\lambda(v(x)) \: \d x,
    \end{align*}
    for all $v \in H_n$. In the following, for convenience, we treat the nonlocal part of $\mathcal{E}_\lambda$ as a bilinear nonlocal operator $\Lambda: H_n \times H_n \to \mathbb R$, so that
    \begin{equation} \label{eq:lambda}
        \begin{split}
        \Lambda(h, k) & = \int_\OO\int_\OO K(x-y)(h(x)-h(y))(k(x)-k(y))\: \d y\:\d x \\
        & = 2(1, K*hk)_H-2(h, K*k)_H
    \end{split}
    \end{equation}
    the equality holding by means of Assumption \ref{hyp:K}, ensuring that the convolution operator is a symmetric map on $H \times H$.
    \begin{lem} \label{lem:nonlocal}
        The bilinear form $\Lambda$ is well defined on $H \times H$. Moreover, $\Lambda$ is bounded, i.e., there exists a positive constant $C > 0$ such that
        \[
        |\Lambda(h, k)| \leq C\|h\|_H\|k\|_H
        \]
        for all $h$ and $k \in H$.
    \end{lem}
    \begin{proof}
        Let $h$ and $k$ to be arbitrary but fixed elements of $H$. The claim is a direct consequence of Young's convolutional inequality. Indeed,
        \[
        \begin{split}
             |\Lambda(h, k)| & = |2(1, K*hk)_H-2(h, K*k)_H| \\
             & \leq 2|(1, K*hk)_H|+ 2|(h, K*k)_H| \\
             & \leq 2\|1\|_{L^\infty(\OO)}\|K * hk\|_{L^1(\OO)} + 2\|h\|_H\|K*k\|_H \\
             & \leq 2\|K\|_{L^1(B_R)} \left( \|hk\|_{L^1(\OO)} + \|h\|_H\|k\|_H \right) \\
             & \leq 4\|K\|_{L^1(B_R)}\|h\|_H\|k\|_H,
        \end{split}
        \]
        and the proof is complete.
    \end{proof} \noindent
    Moreover, let us show a simple, yet crucial property of the approximated energy functional.
    \begin{lem} \label{lem:positivity}
        \review{For any constant $M > 0$, there exists $\eta = \eta(M) > 0$ and a constant $C_\eta > 0$ only depending on $\eta$ such that
        \[
        \mathcal E_\lambda(v) > M\|v\|^2_H- C_\eta
        \]}\noindent
        for all $v \in H$ and all $\lambda \in (0,\eta)$.
    \end{lem}
    \begin{proof}
        Fix $v \in H$ and $\eta > 0$, and let $\lambda \in (0, \eta)$. Applying \cite[Lemma 2.11]{GGG} and owing to the computations in the proof of Lemma \ref{lem:nonlocal}, we get
        \[
        \begin{split}
            \mathcal E_\lambda(v) & = \dfrac 14 \Lambda(v, v) + \int_\OO F_\lambda(v(x)) \: \d x \\
            & = \dfrac 14 \Lambda(v, v) + \int_\OO \Psi_\lambda(v(x)) - \review{\dfrac{\gamma}{2}}v^2(x) \: \d x + |\OO|F(0) \\
            & \geq \dfrac{1}{4\eta}\|v\|^2_H - \|K\|_{L^1(\OO)}\|v\|^2_H -\review{\dfrac{\gamma}{2}}\|v\|^2_H - C_\eta \\
            & = \left( \dfrac{1}{4\eta} - \|K\|_{L^1(\OO)} - \gamma \right)\|v\|^2_H - C_\eta, 
        \end{split}
        \]
        so that by choosing $\eta > 0$ such that
        \[
       \dfrac{1}{4\eta} - \|K\|_{L^1(\OO)} - \gamma  > \review{M} \Leftrightarrow \eta < \dfrac{1}{4(\review{M}+\|K\|_{L^1(\OO)}+\gamma)}
        \]
        the proof is complete.
    \end{proof} \noindent
    Lemma \ref{lem:positivity} implies that the approximated energy functional is almost coercive if $\lambda$ is sufficiently small. Without loss of generality, we shall then assume that, for some fixed $0<\christoph{\nu}< \eta$, 
    \[
    \lambda < \min\{\christoph{\nu}, 1\}
    \]
    and consider, for convenience and without relabeling, $\mathcal E_\lambda$ to be coercive by adding the constant $C_{\christoph{\nu}}$, i.e., 
    \[
    \mathcal E_\lambda: H_n \to \mathbb R \qquad \mathcal{E}_\lambda(v) := \frac{1}{4}\int_\OO\int_\OO K(x-y)|v(x) - v(y)|^2 \: \d y \: \d x + \int_\OO F_\lambda(v(x)) \: \d x + C_{\christoph{\nu}},\qquad \forall \: v\in H_n.
    \]In order to apply the It\^o lemma, we need to show that $\mathcal{E}_\lambda$ is twice Fréchet differentiable. Observe that the constant $C_{\christoph{\nu}}$ does not affect any derivative of $\mathcal E_\lambda$.
    \begin{lem}
        The functional $\mathcal E_\lambda$ is twice Fréchet differentiable with uniformly continuous derivatives on bounded subsets of $H_n$.
    \end{lem}
    \begin{proof}
        First, observe that $\mathcal{E}_\lambda$ is indeed Fréchet-differentiable with $D\mathcal{E}_\lambda: H_n\rightarrow H_n^*$ given by  
    \begin{align*}
        D\mathcal{E}_\lambda(v)[h] := \frac{1}{2}\Lambda(v, h) + \int_\OO F^\prime_\lambda(v(x))h(x) \: \d x,\qquad \forall \: v,\,h\in H_n.
    \end{align*}
    Indeed, for any $v,\,h\in H_n$, we have 
    \begin{align*}
        &\mathcal{E}_\lambda(v+h) - \mathcal{E}_\lambda(v) -D\mathcal{E}_\lambda(v)[h] \\
        &\qquad = \frac{1}{4}\int_\OO\int_\OO K(x-y)|v(x)+h(x)-v(y)-h(y)|^2\: \d y\: \d x + \int_\OO F_\lambda(v(x)+h(x))\:\d x \\
        &\qquad \qquad - \frac{1}{4}\int_\OO\int_\OO K(x-y)|v(x)-v(y)|^2\: \d y\: \d x - \int_\OO F_\lambda(v(x))\: \d x \\
        &\qquad \qquad - \frac{1}{2}\int_\OO\int_\OO K(x-y)(v(x)-v(y))(h(x)-h(y))\: \d y\: \d x - \int_\OO F^\prime_\lambda(v(x))h(x)\: \d x.
    \end{align*}
    Using the identity
    \begin{align*}
        |v(x)-v(y)+h(x)-h(y)|^2 = |v(x)-v(y)|^2 + |h(x)-h(y)|^2 + 2(v(x)-v(y))(h(x)-h(y)),
    \end{align*}
    as well as the fundamental theorem of calculus, we obtain
    {\small
    \begin{align*} 
        &\big|\mathcal{E}_\lambda(v+h) - \mathcal{E}_\lambda(v) -D\mathcal{E}_\lambda(v)h\big| \\
        & \leq \frac{1}{4}\int_\OO\int_\OO K(x-y)|h(x)-h(y)|^2\: \d y\: \d x + \left|\int_\OO F_\lambda(v(x)+h(x))\: \d x - \int_\OO F_\lambda(v(x))\: \d x - \int_\OO F^\prime_\lambda(v(x))h(x)\: \d x\right| \\
        &= \frac{1}{4}\int_\OO\int_\OO K(x-y)|h(x)-h(y)|^2\:\d y\:\d x + \left|\int_0^1\int_\OO F^\prime_\lambda(v(x)+th(x))h(x)\: \d x\: \d t - \int_\OO F^\prime_\lambda(v(x))h(x)\: \d x\right| \\
        &= \frac{1}{4}\int_\OO\int_\OO K(x-y)|h(x)-h(y)|^2\: \d y\: \d x + \left|\int_0^1\int_\OO\Big(F^\prime_\lambda(v(x)+th(x)) - F^\prime_\lambda(v(x))\Big)h(x)\: \d x\: \d t\right| \\
        &\leq C\|h\|^2_H + C_\lambda\|h\|_H^2 \\
        &\leq C_\lambda\|h\|_{H_n}^2, 
    \end{align*}}\noindent
    where we used Lemma \ref{lem:nonlocal}, Assumption \ref{hyp:K} and the Lipschitz continuity of $F_\lambda^\prime$. This shows that $\mathcal{E}_\lambda$ is Fréchet differentiable.
    Now, let us show that also $D\mathcal{E}_\lambda$ is Fréchet differentiable with $D^2\mathcal{E}_\lambda: H_n \rightarrow \cL(H_n,H_n^*)$ given by 
    \begin{align*}
        D^2\mathcal{E}_\lambda(v)[h,k] := \frac{1}{2}\Lambda(h,k)+ \int_\OO F^{\prime\prime}_\lambda(v(x))h(x)k(x)\:\d x
    \end{align*}
    for any $v,h,k\in H_n$. In fact, for every $v,h,k\in H_n$, it holds
    \begin{align*}
        &\left|D\mathcal{E}_\lambda(v+k){[h]} - D\mathcal{E}_\lambda(v){[h]} - D^2\mathcal{E}_\lambda(v)[h,k]\right| \\
        &\qquad= \left|\int_\OO F^\prime_\lambda(v(x)+k(x))h(x)\:\d x - \int_\OO F^\prime_\lambda(v(x))h(x)\:\d x - \int_\OO F^{\prime\prime}_\lambda(v(x))h(x)k(x)\:\d x\right| \\
        &\qquad= \left|\int_0^1\int_\OO \big(F^{\prime\prime}_\lambda(v(x)+tk(x))-F^{\prime\prime}_\lambda(v(x))\big)h(x)k(x)\:\d x\:\d t\right| \\
        &\qquad\leq \|h\|_{L^\infty(\OO)}\|k\|_{L^\infty(\OO)}\int_0^1\|F^{\prime\prime}_\lambda(v+tk)-F^{\prime\prime}_\lambda(v)\|_{L^1(\OO)}\:\d t \\
        &\qquad\leq C\|h\|_{H_n}\|k\|_{H_n}\int_0^1\|F^{\prime\prime}_\lambda(v+tk)-F^{\prime\prime}_\lambda(v)\|_{L^1(\OO)}\:\d t,
    \end{align*}
    where we used Hölder's inequality and the continuous embedding $H_n \hookrightarrow L^\infty(\OO)$. Recalling that $F_\lambda^{\prime\prime}$ is continuous and bounded, the dominated convergence theorem yields 
    \begin{align*}
        \int_0^1\|F^{\prime\prime}_\lambda(v+tk)-F^{\prime\prime}_\lambda(v)\|_{L^1(\OO)}\:\d t \rightarrow 0 \qquad \text{ as }\|k\|_{H_n}\rightarrow 0.
    \end{align*}
     Hence, we obtain 
     \begin{align*}
         \sup_{\|h\|_{H_n}\leq 1}\Big|D\mathcal{E}_\lambda(v+k){[h]} - D\mathcal{E}_\lambda(v){[h]} - D^2\mathcal{E}_\lambda(v)[h,k]\Big| = {o}(\|k\|_{H_n})\qquad \text{ as }\|k\|_{H_n}\rightarrow 0.
     \end{align*}
    In particular, this shows that $D\mathcal{E}_\lambda$ is Fréchet-differentiable with derivative $D^2\mathcal{E}_\lambda$ as above. Moreover, we note that by similar arguments, it is easily shown that both $D\mathcal{E}_\lambda$ and $D^2\mathcal{E}_\lambda$ are uniformly continuous and bounded on bounded subsets of $H_n$. Indeed, this follows from Lemma \ref{lem:nonlocal}, as well as the facts that $F^\prime_\lambda$ is Lipschitz-continuous and that $F^{\prime\prime}_\lambda$ is continuous and bounded. 
    \end{proof} \noindent
    Owing to the previous lemma, we are in a position to \review{show the next estimate}.
    \begin{lem} \label{lem:second_estimate}
        \review{It holds that
        \begin{enumerate}[(i)]
            \item the sequence $\{\mu_{\lambda,n}\}_{n \in \mathbb N}$ is uniformly bounded in $L^p(\Omega;L^2(0,T;V))$;
            \item the sequence $\{\Psi_\lambda(\varphi_{\lambda, n})\}_{n \in \mathbb N}$ is uniformly bounded in ${L^\frac p2(\Omega; L^\infty(0,T;L^1(\OO)))}$;
        \end{enumerate}
        for any $p \geq 2$.}
    \end{lem}
    \begin{proof}
        First of all, we note that $D\mathcal{E}_\lambda(\varphi_{\lambda,n}) = \mu_{\lambda,n}$ {due to Assumption \ref{hyp:K}.}  Hence, \review{applying the It\^{o} lemma to the regularized energy functional $\mathcal E_\lambda$, }we obtain
    \begin{multline} \label{eq:estimate21}
        \mathcal{E}_\lambda(\varphi_{\lambda,n}(t)) + \int_0^t\int_\OO|\nabla\mu_{\lambda,n}(s)|^2\:\d x\:\d s 
        = \mathcal{E}_\lambda(\varphi_{0,n}) + \int_0^t\Big(\mu_{\lambda,n}(s),G_n(\varphi_{\lambda,n}(s))\:\d W(s)\Big)_H \\
        + \frac{1}{2}\int_0^t\sum_{k=0}^{+\infty}\left[\frac{1}{2}\Lambda(G_n(\varphi_{\lambda,n}(s))[u_k], G_n(\varphi_{\lambda,n}(s))[u_k])
        +\int_\OO F^{\prime\prime}_\lambda(\varphi_{\lambda,n}(s))|G_n(\varphi_{\lambda,n}(s))[u_k]|^2\:\d x
        \right]\:\d s.
    \end{multline}
    Appealing to Lemma \ref{lem:nonlocal} and Assumption \ref{hyp:G} yields
    \[
    \begin{split}
    \int_0^t\sum_{k=0}^{+\infty}\left|\Lambda(G_n(\varphi_{\lambda,n}(s))[u_k], G_n(\varphi_{\lambda,n}(s))[u_k])\right| \: \d s  & \leq C\int_0^t \|G_n(\varphi_{\lambda, n}(s))\|^2_{\cL^2(U, H)} \: \d s \\
    & \leq C t,
    \end{split}
    \]
    and similarly
    \[
    \begin{split}
        \int_0^t\sum_{k=0}^{+\infty}\int_\OO F^{\prime\prime}_\lambda(\varphi_{\lambda,n}(s))|G_n(\varphi_{\lambda,n}(s))[u_k]|^2\:\d x \: \d s & \leq C_\lambda\int_0^t \|G_n(\varphi_{\lambda, n}(s))\|^2_{\cL^2(U, H)} \: \d s \\
    & \leq C_\lambda t.
    \end{split}
    \]
    As far as the initial energy is concerned, once again by Lemma \ref{lem:nonlocal} and the boundedness of $F''_\lambda$ we have
    \[
    |\mathcal E_\lambda(\varphi_{0,n})| \leq C_\lambda\left(1 + \|\varphi_{0,n}\|^2_H \right) \leq C_\lambda\left(1 + \|\varphi_{0}\|^2_H \right),
    \]
    where we also used the fact that $\|\varphi_{0,n}\|_H \leq \|\varphi_0\|_H$ by construction, $\P$-almost surely. Combining the above estimates, and then taking $\frac{p}{2}$-powers, the supremum over $[0,\christoph{t}]$ and $\P$-expectations in \eqref{eq:estimate21}, we get, recalling the definition of $\mathcal{E}_\lambda$
    \begin{multline}\label{Est_2}
        \E\supt\mathcal E_\lambda^\frac p2(\varphi_{\lambda, n}(\tau)) + \E\left| \int_0^t \|\nabla\mu_{\lambda,n}(\tau)\|_{\b H}^2 \: \d \tau \right|^\frac p2 \\\
        \leq C_\lambda\left(1+\E\|\varphi_{0}\|_H^{p} +\E\supt\left|\int_0^\tau\big(\mu_{\lambda,n}(s),G_n(\varphi_{\lambda,n}(s))\:\d W(s)\big)_H\right|^{\frac{p}{2}} \right),
    \end{multline}
    where the constant $C_\lambda$ is independent of $n$, but depends on $\lambda$ and $T$.
    Next, by the Burkholder–Davis–Gundy inequality and \eqref{eq:calc_G}, we have
    \begin{equation} \label{eq:burkholdermu}
    \begin{split}
        &\E\supt\left|\int_0^\tau \big(\mu_{\lambda,n}(s),G_n(\varphi_{\lambda,n}(s)) \: \d W(s)\big)_H\right|^{\frac{p}{2}} \\
        &\qquad \leq C\E\left[\int_0^t\|\mu_{\lambda,n}(s)\|_H^2\|G_n(\varphi_{\lambda,n}(s))\|_{\cL^2(U,H)}^2 \: \d s\right]^{\frac{p}{4}} \\
        &\qquad \leq C\E\left[\supt \|G_n(\varphi_{\lambda,n}(\tau))\|_{\cL^2(U,H)}^2\int_0^t\|\mu_{\lambda,n}(s)\|_{H}^2 \: \d s\right]^{\frac{p}{4}} \\
        & \qquad \leq C\E\left[\int_0^t\|\mu_{\lambda,n}(s)\|_{H}^2 \: \d s\right]^{\frac{p}{4}} \\
        & \qquad \leq \varepsilon \E \left| \int_0^t \|\mu _{\lambda,n}(s)\|^2_{H} \: \d s \right|^{\frac{p}{2}} + C
    \end{split}
    \end{equation}
    where $\varepsilon > 0$ is arbitrary and $C$ is independent of both $n$ and $\lambda$. In order to handle the last term, we resort to the Poincaré inequality. Choosing $\varepsilon > 0$ suitably and by summing and subtracting the integral average of the approximated chemical potential, we have
    \begin{equation} \label{eq:interp_mu}
        \varepsilon \E \left| \int_0^t \|\mu _{\lambda,n}(s)\|^2_{H} \: \d s \right|^{\frac{p}{2}} \leq \dfrac 12 \E \left| \int_0^t \|\nabla \mu _{\lambda,n}(s)\|^2_{\b H} \: \d s \right|^{\frac{p}{2}} + C\E \left| \int_0^t |\overline{\mu_{\lambda, n}(s)}|^2 \: \d s \right|^{\frac{p}{2}}.
    \end{equation}
    Testing the equation for the approximated chemical potential in \eqref{eq:NCH_galerkin} by 1, we have 
    \begin{align*}
        |\OO|\overline{\mu_{\lambda,n}} = (\mu_{\lambda,n}, 1)_H = (a\varphi_{\lambda,n}, 1)_H + (K * \varphi_{\lambda,n}, 1)_H + (F^\prime_\lambda(\varphi_{\lambda,n}), 1)_H,
    \end{align*}
    implying, recalling that $F'_\lambda$ is linearly bounded,
    \begin{equation} \label{eq:chemicalV}
            \E \left| \int_0^t |\overline{\mu_{\lambda, n}(s)}|^2 \: \d s \right|^{\frac{p}{2}} \leq C_\lambda\left( 1 + \E \left| \int_0^t \|\varphi_{\lambda, n}(s)\|_H^2 \: \d s \right|^{\frac{p}{2}}\right)
    \end{equation}
    via Young's convolutional inequality and the H\"{o}lder inequality. Combining all of the above in \eqref{Est_2}, also recalling \review{Lemma \ref{lem:first_estimate}},  we arrive at
    \begin{equation*}
        \E\supt\mathcal E_\lambda^\frac p2(\varphi_{\lambda, n}(\tau)) + \E\left| \int_0^t \|\nabla\mu_{\lambda,n}(\tau)\|_{\b H}^2 \: \d \tau \right|^\frac p2 
        \leq C_\lambda\left(1+\|\varphi_{0}\|_H^{p} \right).
    \end{equation*}
    The Gronwall lemma and \eqref{eq:chemicalV} then give 
    \begin{align}\label{eq:secondestimate}
        \|\mu_{\lambda,n}\|_{L^p(\Omega;L^2(0,T;V))}\leq C_\lambda.
    \end{align}
    Finally, arguing by comparison observe that
    \[
    \left|\int_\OO \Psi'_\lambda(\varphi_{\lambda, n}) \: \d x\right| \leq  |\mathcal E_\lambda(\varphi_{\lambda, n})| + C\|\varphi_{\lambda, n}\|^2_H + C_\Lambda
    \]
    and hence we conclude that
    \begin{align*}\label{eq:secondestimatefinal}
    \|\mu_{\lambda,n}\|_{L^p(\Omega;L^2(0,T;V))} + \|\Psi_\lambda(\varphi_{\lambda, n})\|_{L^\frac p2(\Omega; L^\infty(0,T;L^1(\OO)))} \leq C_{2, \lambda},
    \end{align*}
    and the proof is complete.
    \end{proof} \noindent
    \review{Concerning the stochastic integrals, we have the following result.}
    \begin{lem} \label{lem:third_estimate}
        \review{It holds that
        \begin{enumerate}[(i)] \itemsep0.3em
            \item the sequence $\{G_n(\varphi_{\lambda,n})\}_{n \in \mathbb N}$ is uniformly bounded in \[L^\infty(\Omega\times(0,T);\cL^2(U,H)) \cap L^p(\Omega;L^2(0,T;\cL^2(U,V)))\]
            for any $p \geq 2$;
            \item the sequence $\left\{\displaystyle\int_0^. G_n(\varphi_{\lambda,n}(s)) \: \d W(s)\right\}_{n \in \mathbb N}$ is uniformly bounded in \[L^q(\Omega;W^{k,q}(0,T;H))\cap L^2(\Omega;H^{k}(0,T;V))\] 
            for every $k \in (0,\frac 12)$ and $q \geq 2$.
        \end{enumerate}
        }
    \end{lem}
    \begin{proof}
         A direct application of Assumption \ref{hyp:G} yields 
    \begin{align*}
        \|G_n(\varphi_{\lambda,n})\|_{\cL^2(U,V)}^2 &\leq \|G(\varphi_{\lambda,n})\|_{\cL^2(U,V)}^2 \leq 2L_G^2(1 + \|\varphi_{\lambda,n}\|^2_V)
    \end{align*}
    and thus, using \review{Lemma \ref{lem:first_estimate}}, it follows that
    \begin{align*}
    \E\left|\int_0^t\|G_n(\varphi_{\lambda,n}(s))\|_{\cL^2(U,V)}^2 \: \d s\right|^\frac p2 \leq \review{C}.
    \end{align*}
    \review{Therefore}, the estimate above and similar computations imply
    \begin{align}
    \|G_n(\varphi_{\lambda,n})\|_{L^\infty(\Omega\times(0,T);\cL^2(U,H))\cap L^p(\Omega;L^2(0,T;\cL^2(U,V)))} \leq C_\lambda,
    \end{align}
    where the constant $C_\lambda$ only depends on $\lambda>0$. 
    By \cite[Lemma 2.1]{fland-gat}, an estimate on the corresponding It\^{o} integrals holds, namely
    \begin{align}\label{eq:unif_est_3}
        \left\|\int_0^. G_n(\varphi_{\lambda,n}(s)) \: \d W(s)\right\|_{L^q(\Omega;W^{k,q}(0,T;H))\cap L^2(\Omega;H^{k}(0,T;V))} \leq C_{\lambda},
    \end{align}
    for every $k\in(0,\frac{1}{2})$ and $q\geq 2$.
    \end{proof} \noindent
    \review{Finally, the last estimate follows by a comparison argument.}
    \begin{lem} \label{lem:fourth_estimate}
        \review{The sequence $\{\varphi_{\lambda, n}\}_{n \in \mathbb N}$ is uniformly bounded in $L^p(\Omega;W^{\alpha, p}(0,T;V^*))$ for any $\alpha \in (0, \frac 12)$ and $p \geq 2$.}
    \end{lem}
    \begin{proof}
        Let us consider the weak formulation of the discretized problem as an equality in $V^*$, i.e.
    \begin{equation*}
        \langle \varphi_{\lambda,n}(t),\psi\rangle_{V^*,V} + \int_0^t\int_\OO\nabla\mu_{\lambda,n}(s)\cdot\nabla\psi\:\d x \, \d s
        = \langle\varphi_{0,n},\psi\rangle_{V^*,V} + \Big(\int_0^tG_n(\varphi_{\lambda,n}(s))\: \d W(s),\psi\Big)_H
    \end{equation*}
    for every $\psi\in V$ such that $\|\psi\|_V=1$, all $t \in [0,T]$ and $\P$-almost surely. Since
    \begin{align*}
        \left|\int_\OO\nabla\mu_{\lambda,n}(s)\cdot\nabla\psi\:\d x\right| \leq \|\nabla\mu_{\lambda,n}(s)\|_\b H
    \end{align*}
    and thus, by \review{Lemma \ref{lem:second_estimate}}, we arrive at
    \begin{align*}
        \left\|\int_0^t\int_\OO\Delta\mu_{\lambda,n}(s)\:\d x \, \d s\right\|_{L^p(\Omega;H^1(0,T;V^*))}\leq C_\lambda,
    \end{align*}
    where $C_\lambda>0$ depends on $\lambda,p$ and $T$. Now, recalling  
    \begin{align*}
        |\langle\varphi_{0,n},\psi\rangle_{V^*,V}| \leq \|\varphi_{0,n}\|_H \leq \|\varphi_{0}\|_H
    \end{align*}
    and estimate \eqref{eq:unif_est_3}, we eventually get by comparison and invoking \cite[Lemma 2.1]{DPGS}
    \begin{align*}\label{eq:unif_est_4}      \|\varphi_{\lambda,n}\|_{L^p(\Omega;W^{\alpha,p}(0,T;V^*))} \leq C_{4, \lambda}
    \end{align*}
    for some $\alpha = \alpha(p) > \frac 1p$ if $p > 2$ and $\alpha \in (0, \frac 12)$ if $p = 2$. Note that the constant may depend on $\lambda,\alpha$ and $T$. \review{The proof is complete.}
    \end{proof}
    \subsection{Passage to the limit as $n \to+\infty$} 
    \label{ssec:lim_n}
    On account of the previous estimates, we now pass to the limit $n\rightarrow+\infty$, keeping $\lambda>0$ fixed, by means of a stochastic compactness argument.
    First of all, let us show that the sequence of laws of $\{\varphi_{\lambda,n}\}_{n\in \enne}$ is tight in a suitable path space. To this end, we recall that the embeddings (see \cite[Theorem 2.1]{fland-gat})
    \begin{align*}
        L^2(0,T;V)\cap H^\alpha(0,T;V^*) & \hookrightarrow L^2(0,T;H) \\
        L^\infty(0,T;H) \cap W^{\alpha, p}(0,T;V^*) &\hookrightarrow C^0([0,T]; H^{-\delta}(\OO))
    \end{align*}
    are compact, with the caveat that the second embedding only holds if $\alpha p > 1$, i.e., only if $p > 2$, for some small $\delta > 0$.
    \begin{lem} \label{lem:tightness}
        If $p > 2$, then the family of laws of $\{\varphi_{\lambda,n}\}_{n\in \enne}$ is tight in the space $\mathcal X := L^2(0,T;H) \cap C^0([0,T]; H^{-\delta}(\OO))$ for any $\delta \in (0, \frac 12)$.
    \end{lem}
    \begin{proof}
        For any $R>0$, let $B_R$ denote the closed ball of radius $R$ in $L^2(0,T;V)\cap H^\alpha(0,T;V^*)$. Then, Markov's inequality as well as \review{Lemmas \ref{lem:first_estimate} and \ref{lem:fourth_estimate}} imply
    \begin{align*}
        \P\{\varphi_{\lambda,n}\in B_R^c\} 
        &= \P\{\|\varphi_{\lambda,n}\|_{L^2(0,T;V)\cap H^\alpha(0,T;V^*)}\geq R\} \\
        &\leq \frac{1}{R^p}\E\|\varphi_{\lambda,n}\|_{L^2(0,T;V)\cap H^\alpha(0,T;V^*)}^p \leq \frac{C}{R^p}.
    \end{align*}
    In particular, this yields
    \begin{align*}
    \lim\limits_{R\rightarrow+\infty}\sup_{n\in\enne}\P\{\varphi_{\lambda,n}\in B_R^c\} = 0.
    \end{align*}
    This shows that the family of laws $(\varphi_{\lambda,n})_{n\in\enne}$ on $L^2(0,T;H)$ is tight. An analogous argument can be carried out for the second embedding.
    \end{proof} \noindent
    Let us now set
    \begin{align*}
       G_n(\varphi_{\lambda,n})\cdot W :=  \int_0^. G_n(\varphi_{\lambda,n}(s))\d W(s).
    \end{align*}
    In the spirit of Lemma \ref{lem:tightness}, using \eqref{eq:unif_est_3} and \cite[Theorem 2.2]{fland-gat} we exploit the compact embedding
    \[
    W^{k, q}(0,T; H) \hookrightarrow C^0([0,T]; V^*)
    \]
    holding whenever $kq > 1$ (here, we are free to choose $q \geq 2$ and $k > \frac 1q$). Hence, we immediately have
    \begin{lem} \label{lem:tightness2}
        The family of laws of $\{G_n(\varphi_{\lambda,n})\cdot W\}_{n \in \mathbb N}$ is tight in the space $C^0([0,T]; V^*)$.
    \end{lem} \noindent
    By identifying $W$ with a constant sequence of random variables $\{W_n\}_{n \in \mathbb N}$ with values in $C^0([0,T];U_0)$, we also get that 
    \begin{lem} \label{lem:tightness3}
        The family of laws of $\{W_n\}_{n \in \mathbb N}$ is tight in the space $C^0([0,T]; U_0)$.
    \end{lem} \noindent
    Finally, the sequence of initial conditions is clearly tight in $V^*$ as $H \hookrightarrow V^*$ compactly. Namely, we have
    \begin{lem} \label{lem:tightness4}
        The family of laws of $\{\varphi_{0,n}\}_{n \in \mathbb N}$ is tight in $V^*$.
    \end{lem} \noindent
    Collecting Lemmas \ref{lem:tightness}-\ref{lem:tightness4}, we obtain that the family of laws 
    \begin{align*}
        \{\varphi_{\lambda,n},\, G_n(\varphi_{\lambda,n})\cdot W,\, W_n,\, \varphi_{0,n}\}_{n\in\enne}
    \end{align*}
    is tight in the space
    \begin{align*}
        \mathcal X \times C^0([0,T]; V^*) \times C^0([0,T];U_0) \times V^*.
    \end{align*}
    Owing to the Prokhorov and Skorokhod theorems (see \cite[Theorem 2.7]{ike-wata} and \cite[Theorem 1.10.4, Addendum 1.10.5]{vaa-well}), there exists a probability space $(\tom, \widetilde{\mathscr{F}},\tP)$ and a sequence of random variables $X_n:  (\tom, \tF)\to(\Omega, \mathscr{F})$ such that the law of $X_n$ is $\P$ for every $n \in \enne$, namely $\tP \circ X_n^{-1} = \P$ (so that composition with $X_n$ preserves laws), and the following convergences hold
    \review{
    \begin{align*}
        \widetilde{\varphi}_{\lambda,n} := \varphi_{\lambda,n}\circ X_n \rightarrow \widetilde{\varphi}_\lambda&\quad\text{ in }\mathcal X, \, \tP\text{-a.s.,}\\
        \widetilde{I}_{\lambda,n} := (G_n(\varphi_{\lambda,n})\cdot W)\circ X_n \rightarrow \widetilde{I}_\lambda&\quad\text{ in }C^0([0,T];V^*), \, \tP\text{-a.s.,}\\
        \widetilde{W}_n:= W\circ X_n \rightarrow \widetilde{W} &\quad\text{ in }C^0([0,T];U_0),\, \tP\text{-a.s.,}\\
        \widetilde{\varphi}_{0,n} := \varphi_{0,n} \circ X_n \rightarrow \widetilde{\varphi}_{0}&\quad\text{ in }V^*,\, \tP\text{-a.s.,}
    \end{align*}
    and, in turn, exploiting the Vitali convergence theorem, as well as the Banach--Alaoglu theorem and reflexivity, we conclude that}
    \begin{align*}
        \widetilde{\varphi}_{\lambda,n} \rightarrow \widetilde{\varphi}_\lambda&\quad\text{ in }L^q(\tom;L^2(0,T;H) \cap C^0([0,T]; H^{-\delta}(\OO))\;\;\forall \: q<p,\\
        \widetilde{\varphi}_{\lambda,n} \overset{*}{\rightharpoonup} \widetilde{\varphi}_\lambda&\quad\text{ in }L^p_w(\tom; L^\infty(0,T;H)) \cap L^p(\tom;H^\alpha(0,T;V^*)),\\
        \widetilde{\varphi}_{\lambda,n} \rightharpoonup \widetilde{\varphi}_\lambda&\quad\text{ in }L^p(\tom;L^2(0,T;V)), \\
        \widetilde{\mu}_{\lambda,n} := \mu_{\lambda,n}\circ X_n \rightharpoonup \widetilde{\mu}_\lambda&\quad\text{ in }L^{p}(\tom;L^2(0,T;V)), \\
        \widetilde{I}_{\lambda,n} \rightarrow \widetilde{I}_\lambda&\quad\text{ in }L^q(\tom;C^0([0,T];V^*))\;\;\forall \: q<p,\\
        \widetilde{W}_n\rightarrow \widetilde{W} &\quad\text{ in }L^q(\tom;C^0([0,T];U_0))\;\;\forall \: q<p, \\
        \widetilde{\varphi}_{0,n}\rightarrow \widetilde{\varphi}_{0}&\quad\text{ in }L^q(\tom;V^*)\;\;\forall \: q<p,
    \end{align*}
    for some limiting measurable processes enjoying the following regularity properties:
    \begin{align*}
        \widetilde{\varphi}_\lambda&\in L^p_w(\tom;L^\infty(0,T;H))\cap L^p(\tom;H^\alpha(0,T;V^*))\cap L^p(\tom;L^2(0,T;V) \cap C^0([0,T]; H^{-\delta}(\OO)),\\
        \widetilde{\mu}_\lambda&\in L^p(\tom;L^2(0,T;V)),\\
        \widetilde{I}_\lambda &\in L^p(\tom;C^0([0,T];V^*)), \\
        \widetilde{W}&\in L^p(\tom;C^0([0,T];U)), \\
        \widetilde{\varphi}_0 & \in L^p(\tom; H).
    \end{align*}
    Observe that from the proof of Skorokhod's theorem, the probability space $(\tom,\widetilde{\mathscr{F}},\tP)$ can be made independent of $\lambda$. 
    \subsection{Identification of the limit solution}
    \label{ssec:limit_n}
    Before moving to proving that the limiting quantities satisfy the regularized problem \eqref{eq:NCH_lambda}, we deduce the behavior of the nonlinearities of the problem. Using the strong convergence of
    \[
\widetilde{\varphi}_{\lambda,n} \rightarrow \widetilde{\varphi}_\lambda\quad\text{ in }L^q(\tom;L^2(0,T;H))
\]
for $q < p$, jointly with the Lipschitz-continuity of $F_\lambda^\prime$ (that is uniform in $n$), we immediately have
    \begin{align*}
        F^\prime_\lambda(\widetilde{\varphi}_{\lambda,n}) \rightarrow F^\prime_\lambda(\widetilde{\varphi}_{\lambda})\quad\text{ in } L^q(\tom;L^2(0,T;H))
    \end{align*}
    for $q < p$.
    Similarly, since $G_n$ is Lipschitz-continuous, and on account of 
    \begin{multline*}
        \|G_n(\widetilde{\varphi}_{\lambda,n}) - G(\widetilde{\varphi}_{\lambda})\|_{L^q(\tom;L^2(0,T;\cL^2(U,H)))} \\
        \leq \|G_n(\widetilde{\varphi}_{\lambda,n}) - G_n(\widetilde{\varphi}_{\lambda})\|_{L^q(\tom;L^2(0,T;\cL^2(U,H)))} + \|G_n(\widetilde{\varphi}_{\lambda}) - G(\widetilde{\varphi}_{\lambda})\|_{L^q(\tom;L^2(0,T;\cL^2(U,H)))} \to 0
    \end{multline*}
    as $n \to +\infty$, it follows that 
    \begin{align*}
        G_n(\widetilde{\varphi}_{\lambda,n}) \rightarrow G(\widetilde{\varphi}_{\lambda})\quad \text{ in }L^q(\tom;L^2(0,T;\cL^2(U,H)))
    \end{align*}
    for all $q<p$. 
    Next, we identify the limit of the stochastic integral sequence. To this end, we employ a standard method illustrated, for instance, in \cite[Section 8.4]{dapratozab}. Consider the family of filtrations on $(\tom, \tF, \tP)$ generated by the Galerkin solutions, that is,
	\[
	\tF_{\lambda,n,t} := \sigma\left\{\widetilde{\varphi}_{\lambda,n}(s),\,  \widetilde{I}_{\lambda,n}(s),\widetilde{W}_{n}(s),\,\widetilde{\varphi}_{0,n}, \, s \in [0,t]\right\},
	\]
	for any $t \geq 0$, $n \in \enne$ and $\lambda \in (0,1)$. By definition, we have that the processes $\widetilde W_n$ are clearly adapted to $\tF_{\lambda, n, t}$. In particular, by preservation of laws and the definitions of Wiener process and stochastic integral, we readily have that $\widetilde W_{n}$ is a $Q_0$-Wiener process on some larger space $U_0 \supset U$ and, correspondingly, 
	\[
	\widetilde{I}_{\lambda, n} = \int_0^t G_{n}(\widetilde{\varphi}_{\lambda,n}(\tau)) \: \mathrm{d}\widetilde{W}_{n}(\tau),
	\]
	is a martingale with values in $V^*$. Similarly, by considering the limit processes, we can define the family of increasing filtrations
	\[
	\tF_{\lambda,t} := \sigma\left\{\widetilde{\varphi}_{\lambda}(s),\, \widetilde{I}_{\lambda}(s), \widetilde{W}(s), \,\widetilde{\varphi}_0,\, s \in [0,t]\right\}.
	\]
	Let now $t > 0$, $s \in [0,t]$ and consider the path spaces
	\begin{align*}
		\mathcal X_{s} &:= L^2(0,s;H) \cap C^0([0,s];H^{-\delta}(\OO)), \\
		\mathcal Z_s &:= \mathcal X_s \times C^0([0,s]; V^*)  \times  C^0([0,s]; U^0) \times V^*.
	\end{align*}
	In what follows, we denote by $\psi$ an arbitrary but fixed element of $C_\text{b}(\mathcal Z_s)$, namely, a continuous and bounded function on $\mathcal Z_s$ with values in $\mathbb R$. By definition of martingale (or Wiener process), we have
	\begin{equation} \label{eq:limit_n0}
		\tE \left[ \left(\widetilde{W}_{n}(t) -\widetilde{W}_{n}(s) \right) \psi\left( \widetilde{\varphi}_{\lambda,n},\, \widetilde{I}_{\lambda,n},\, \widetilde{W}_{n},\,\widetilde{\varphi}_{0,n} \right)  \right] = 0
	\end{equation}
	for any bounded and continuous function $\psi \in C_\text{b}(\mathcal Z_s)$. Of course, the arguments of $\psi$ are intended to be restricted over $[0,s]$. Letting $n \to +\infty$ in \eqref{eq:limit_n0}, an application of the dominated convergence theorem, owing to the proven convergences and the properties of $\psi$, entails
	\begin{equation} \label{eq:limit_n1}
		\tE \left[ \left(\widetilde{W}(t) -\widetilde{W}(s) \right) \psi\left( \widetilde{\varphi}_{\lambda},\, \widetilde{I}_{\lambda},\,\widetilde{W},\,\widetilde{\varphi}_0  \right)  \right] = 0,
	\end{equation}
	which implies that $\widetilde{W}$ is a $U_0$-valued $(\tF_{\lambda,t})_t$-adapted martingale. In order to show that $\widetilde W$ is indeed a Wiener process, we resort to the characterization of $Q$-Wiener processes given by \cite[Theorem 4.6]{dapratozab}. To this end, notice that \eqref{eq:limit_n1} means that, for every $v, w \in \review{U_0}$
	{\footnotesize \begin{equation*}
		\tE \left[ \left(\left( \widetilde{W}_{n}(t), v \right)_{U_0}
		\left( \widetilde{W}_{n}(t), w \right)_{U_0}
		-\left( \widetilde{W}_{n}(s), v \right)_{U_0}
		\left( \widetilde{W}_{n}(s), w \right)_{U_0}
		- (t-s)\left(Q_0v,w\right)_{U_0}\right)
		\psi\left( \widetilde{\varphi}_{\lambda,n},\,
		\widetilde{I}_{\lambda,n},\,
		\widetilde{W}_{n},\,\widetilde{\varphi}_{0,n} \right)  \right] = 0,
	\end{equation*}}
	and using once more the dominated convergence theorem, we get
	{\small \begin{equation*}
		\tE \left[ \left(\left( \widetilde{W}(t), v \right)_{U_0}
		\left( \widetilde{W}(t), w \right)_{U_0}
		-\left( \widetilde{W}(s), v \right)_{U_0}
		\left( \widetilde{W}(s), w \right)_{U_0}
		- (t-s)\left(Q_0v,w\right)_{U_0}\right)
		\psi\left( \widetilde{\varphi}_{\lambda},\,
		\widetilde{I}_{\lambda},\,
		\widetilde{W},\,\widetilde{\varphi}_{0} \right)  \right] = 0,
	\end{equation*}}namely the quadratic variation is exactly
	\begin{equation*}
		\left\llangle \review{\widetilde{W}}\right\rrangle(t) = t\review{Q_0}, \qquad t \in [0,T],
	\end{equation*}
	which is enough to conclude that $\widetilde{W}$ is a $Q_0$-Wiener process, adapted to $(\tF_{\lambda,t})_t$, owing to \cite[Theorem 4.6]{dapratozab}. In a similar way, we address the stochastic integrals. Arguing as in \eqref{eq:limit_n0}-\eqref{eq:limit_n1}, we conclude that $\widetilde{I}_{\lambda}$ is a $V^*$-valued martingale. As far as the quadratic variations are concerned, an application of \cite[Theorem 4.27]{dapratozab} yields
	\begin{equation*}
		\left\llangle \widetilde{I}_{\lambda,n}\right\rrangle(t) = \int_0^t G_{n}(\widetilde{\varphi}_{\lambda,n}(\tau)) \circ G_{n}(\widetilde{\varphi}_{\lambda,n}(\tau))^*\: \mathrm{d}\tau,
	\end{equation*}
	for every $t \in [0,T]$. Let $v, w \in V$ be fixed. Then
	\begin{multline*}
		\tE \left[ \left(\left\langle \widetilde{I}_{\lambda,n}(t), v \right\rangle_{V^*, V}\left\langle \widetilde{I}_{\lambda,n}(t), w \right\rangle_{V^*, V} - \left\langle\widetilde{I}_{\lambda,n}(s), v \right\rangle_{V^*, V}\left\langle \widetilde{I}_{\lambda,n}(s), w \right\rangle_{V^*, V} \right. \right. \\ \left. \left. -\int_0^t (G_{n}(\widetilde{\varphi}_{\lambda,n}(\tau)) \circ G_{n}(\widetilde{\varphi}_{\lambda,n}(\tau))^*v, w)_H \: \mathrm{d}\tau\right) \psi\left( \widetilde{\varphi}_{\lambda,n},\,\widetilde{I}_{\lambda,n},\, \widetilde{W}_{n}
		,\,\widetilde{\varphi}_{0,n}\right)  \right] = 0,
	\end{multline*}
	and, as $n \to +\infty$, the dominated convergence theorem implies that
	\begin{multline*}
		\tE \left[ \left(\left\langle \widetilde{I}_{\lambda}(t), v \right\rangle_{V^*, V}\left\langle \widetilde{I}_{\lambda}(t), w \right\rangle_{V^*, V} - \left\langle\widetilde{I}_{\lambda}(s), v \right\rangle_{V^*, V}\left\langle \widetilde{I}_{\lambda}(s), w \right\rangle_{V^*, V} \right. \right. \\ \left. \left. - \int_0^t (G(\widetilde{\varphi}_{\lambda}(\tau)) \circ G(\widetilde{\varphi}_{\lambda}(\tau))^*v, w)_H \: \mathrm{d}\tau\right) \psi\left( \widetilde{\varphi}_{\lambda},\,\widetilde{I}_{\lambda},\, \widetilde{W}
		,\,\widetilde{\varphi}_{0}\right)  \right] = 0,
	\end{multline*}
	namely, we have shown that the quadratic variation of $\widetilde{I}_{\lambda}$ is
	\[
	\left\llangle \widetilde{I}_{\lambda} \right\rrangle(t) = \int_0^t G(\widetilde{\varphi}_{\lambda}(\tau)) \circ G(\widetilde{\varphi}_{\lambda}(\tau))^* \: \mathrm{d}\tau, \qquad t \in [0,T].
	\]
	In particular, we are left to show that $\widetilde{I}_{\lambda}$ can be identified with the process
	\[
	\widetilde{M}_{\lambda}(t) := \int_0^t \review{G(\widetilde{\varphi}_{\lambda}(\tau))}\: \mathrm{d}\review{\widetilde{W}}(\tau),
	\]
	which is a $V^*$-valued $(\tF_{\lambda,t})_t$-adapted martingale having the same quadratic variation of $ \widetilde{I}_{\lambda}$. By \cite[Theorem 3.2]{Pard}, we can write
	\begin{equation} \label{eq:limit_n2}
		\begin{split}
			\left\llangle \widetilde{M}_{\lambda} - \widetilde{I}_{\lambda} \right\rrangle & = \left\llangle \widetilde{M}_{\lambda} \right\rrangle + \left\llangle \widetilde{I}_{\lambda} \right\rrangle -2 \left\llangle \widetilde{M}_{\lambda}, \widetilde{I}_{\lambda} \right\rrangle \\
			& = 2\int_0^\cdot G(\widetilde{\varphi}_{\lambda}(\tau)) \circ G(\widetilde{\varphi}_{\lambda}(\tau))^*\: \mathrm{d}\tau - 2\int_0^\cdot G(\widetilde{\varphi}_{\lambda}(\tau)) \: \mathrm{d}\left\llangle\widetilde{W}, \widetilde{I}_{\lambda} \right\rrangle (\tau).
		\end{split}
	\end{equation}
	Moreover, still by \cite[Theorem 3.2]{Pard}, we have that the cross quadratic variation satisfies
	\[
	\begin{split}
		\left\llangle\widetilde{I}_{\lambda,n}, \widetilde{W}_{n} \right\rrangle & = \int_0^\cdot G_{n}(\widetilde{\varphi}_{\lambda,n}(\tau)) \circ \iota^{-1}_1 \: \mathrm{d}\left\llangle\widetilde{W}_{n}, \widetilde{W}_{n}\right\rrangle (\tau) \\
		& = \int_0^\cdot G_{n}(\widetilde{\varphi}_{\lambda,n}(\tau)) \circ \iota^{-1} \circ Q_0 \: \mathrm{d}\tau \\
		& = \int_0^\cdot G_{n}(\widetilde{\varphi}_{\lambda,n}(\tau)) \circ \iota^{-1} \circ \iota \circ \iota^* \: \mathrm{d}\tau \\
		& = \int_0^\cdot G_{n}(\widetilde{\varphi}_{\lambda,n}(\tau)) \circ \iota^* \: \mathrm{d}\tau,
	\end{split}
	\]
	where we also used the fact that $Q_0 = \iota_1 \circ \iota^*_1$, where $\iota: U \to U_0$ is the classical Hilbert-Schmidt embedding. This implies that
	\[
	\left\llangle \widetilde{W}_{n}, \widetilde{I}_{\lambda,n} \right\rrangle = \int_0^\cdot \iota \circ G_{n}(\widetilde{\varphi}_{\lambda,n}(\tau))^* \: \mathrm{d}\tau.
	\]
	A further application of the dominated convergence theorem entails that, as $n \to +\infty$,
	\begin{equation} \label{eq:limit_n3}
		\left\llangle \widetilde{W}, \widetilde{I}_{\lambda} \right\rrangle = \int_0^\cdot \iota \circ G(\widetilde{\varphi}_{\lambda}(\tau))^* \: \mathrm{d}\tau.
	\end{equation}
	The identification follows injecting \eqref{eq:limit_n3} in \eqref{eq:limit_n2}. Finally, the proven convergences are enough to pass to the limit in the Galerkin problem and retrieve that the limit processes are indeed a martingale solution of \eqref{eq:NCH_lambda}.
    \subsection{Uniform estimates with respect to $\lambda$}
    Here, we want to show further estimates, now uniform in the Yosida parameter $\lambda > 0$. The symbol $C$ (possibly numbered) denotes a positive constant, always independent of $\lambda$, but possibly dependent on other structural quantities of the problem, and may change from line to line. First of all, \review{from the proof of Lemma \ref{lem:first_estimate}} we observe that the constant $C_1$ in estimate $\eqref{eq:unif_est_1}$ is already independent of $\lambda$. 
    \begin{lem} \label{lem:first_estimate_lam}
        \review{The family $\{\widetilde{\varphi}_\lambda\}_{\lambda > 0}$ is uniformly bounded in $L^p(\tom;L^\infty(0,T;H)) \cap L^p(\tom;L^2(0,T;V))$ for all $p \geq 2$.}
    \end{lem}
    \begin{proof}
        By using Lemma \ref{lem:first_estimate}, lower semicontinuity and preservation of laws of $X_n$, we obtain 
        \begin{align*}
    \|\widetilde{\varphi}_{\lambda}\|_{L^p(\tom;L^\infty(0,T;H))} + \|\widetilde{\varphi}_{\lambda}\|_{L^p(\tom;L^2(0,T;V))} \leq C_1,
    \end{align*}
        and the proof is complete.
    \end{proof} \noindent
    Next, the key energy inequality must be revisited in order to obtain independence from $\lambda$. \review{In particular, we shall deduce the following estimate.
    \begin{lem} \label{lem:second_estimate_lam}
        It holds that
        \begin{enumerate}[(i)]
            \item the family $\{\widetilde{\mu}_\lambda\}_{\lambda > 0}$ is uniformly bounded in $L^\frac p2(\tom;L^2(0,T;V))$;
            \item the family $\{\nabla \widetilde{\mu}_\lambda\}_{\lambda > 0}$ is uniformly bounded in $L^p(\tom;L^2(0,T;\b H))$;
            \item the family $\{\Psi_\lambda(\widetilde \varphi_\lambda)\}_{\lambda >0}$ is uniformly bounded in $L^\frac p2(\tom;L^\infty(0,T;L^1(\OO)))$;
        \end{enumerate}
        for all $p \geq 2$.
    \end{lem} 
    }
    \begin{proof}
    Let us start back from
    \begin{multline} \label{eq:estimate21_bis}
        \mathcal{E}_\lambda(\varphi_{\lambda,n}(t)) + \int_0^t\int_\OO|\nabla\mu_{\lambda,n}(s)|^2\:\d x\:\d s 
        = \mathcal{E}_\lambda(\varphi_{0,n}) + \int_0^t\Big(\mu_{\lambda,n}(s),G_n(\varphi_{\lambda,n}(s))\:\d W(s)\Big)_H \\
        + \frac{1}{2}\int_0^t\sum_{k=0}^{+\infty}\left[\frac{1}{2}\Lambda(G_n(\varphi_{\lambda,n}(s))[u_k], G_n(\varphi_{\lambda,n}(s))[u_k])
        +\int_\OO F^{\prime\prime}_\lambda(\varphi_{\lambda,n}(s))|G_n(\varphi_{\lambda,n}(s))[u_k]|^2\:\d x
        \right]\:\d s.
    \end{multline}
    In order to proceed, it will be much more convenient to recast \eqref{eq:estimate21_bis} with ``lower-order'' terms to the right hand side, namely
    \begin{multline} \label{eq:estimate21_ter}
        \int_\OO \Psi_\lambda(\varphi_{\lambda,n}(t)) \: \d x+ \int_0^t\int_\OO|\nabla\mu_{\lambda,n}(s)|^2\:\d x\:\d s 
        = - \Lambda(\varphi_{\lambda,n}(t), \varphi_{\lambda,n}(t)) + \dfrac{\gamma}{2}\int_\OO \varphi_{\lambda,n}^2(t) \: \d x \\+\mathcal{E}_\lambda(\varphi_{0,n}) + \int_0^t\Big(\mu_{\lambda,n}(s),G_n(\varphi_{\lambda,n}(s))\:\d W(s)\Big)_H \\
        + \frac{1}{2}\int_0^t\sum_{k=0}^{+\infty}\left[\frac{1}{2}\Lambda(G_n(\varphi_{\lambda,n}(s))[u_k], G_n(\varphi_{\lambda,n}(s))[u_k])
        +\int_\OO F^{\prime\prime}_\lambda(\varphi_{\lambda,n}(s))|G_n(\varphi_{\lambda,n}(s))[u_k]|^2\:\d x
        \right]\:\d s.
    \end{multline}
    Firstly, observe that the terms in $\Lambda$ (i.e., the first and fifth terms on the right hand side) can be bounded independently of $\lambda$ through Lemma \ref{lem:nonlocal} (\review{as in the proof of Lemma \ref{lem:second_estimate}}). Then, we consider the initial energy.  It holds
    \[
    \begin{split}
        \mathcal E_\lambda({\varphi_{0,n}}) = \dfrac{1}{4}\Lambda({\varphi_{0,n}}, {\varphi_{0,n}}) + \int_\OO \Psi_\lambda(\varphi_{0,n}(x)) - \dfrac{\gamma}{2}\varphi_{0,n}^2(x)\: \d x + C,
    \end{split}
    \]
    where $C = C_\Lambda + F(0)|\OO|$ is independent of $\lambda$. By Lemma \ref{lem:nonlocal} and the convergence $\varphi_{0,n} \to \varphi_0$ in $H$ with $\|\varphi_{0,n}\|^2_H \leq \|\varphi_0\|^2_H$ for all $n \in \mathbb N$, $\P$-almost surely, we get
    \[
    \begin{split}
        \mathcal E_\lambda({\varphi_{0,n}}) & \leq C\left( 1 + \int_\OO \Psi_\lambda(\varphi_{0,n}(x)) \: \d x + \|\varphi_0\|^2_H\right)
    \end{split}
    \]
    $\P$-almost surely for a constant $C$ independent of $\lambda$ and $n$. Now, collecting all the above controls, raising the result to the power $\frac p2$, taking supremums and $\P$-expectations we arrive at
    \begin{multline*}
        \E\supt\mathcal \|\Psi_\lambda(\varphi_{\lambda, n}(\tau))\|^\frac p2_{L^1(\OO)} + \E\left| \int_0^t \|\nabla\mu_{\lambda,n}(\tau)\|_{\b H}^2 \: \d \tau \right|^\frac p2 
        \\ \leq C\left[1+\E \|\varphi_{0}\|_H^{p} + \review{\E \|\Psi_\lambda(\varphi_{0,n})\|_{L^1(\OO)}^\frac p2}\right] \\
        + C\left[ \E\left|\int_0^t\sum_{k=0}^{+\infty}\int_\OO |F_\lambda^{\prime\prime}(\varphi_{\lambda,n})||G_n(\varphi_{\lambda,n})[u_k]|^2\:\d x \: \d s\right|^{\frac{p}{2}} + \E \supt \left| \int_0^\tau \Big(\mu_{\lambda,n}(s),G_n(\varphi_{\lambda,n}(s))\:\d W(s)\Big)_H\right|^\frac p2 \right].
    \end{multline*}
    Using the Burkholder-Davis-Gundy inequality, as in the \review{proof of Lemma \ref{lem:second_estimate}} (\review{concatenating, in particular \eqref{eq:burkholdermu}, ignoring the last line, and \eqref{eq:interp_mu} replacing, $\frac p2$ with $\frac p4$}), we can bound the last term and obtain an additional term depending on the spatial average of the chemical potential, i.e.,
    \begin{multline} \label{eq:estimate_aux}
        \E\supt\mathcal \|\Psi_\lambda(\varphi_{\lambda, n}(\tau))\|^\frac p2_{L^1(\OO)} + \E\left| \int_0^t \|\nabla\mu_{\lambda,n}(\tau)\|_{\b H}^2 \: \d \tau \right|^\frac p2 
        \leq C\left[1+\E \|\varphi_{0}\|_H^{p} + \review{\E \|\Psi_\lambda(\varphi_{0,n})\|_{L^1(\OO)}^\frac p2}\right] \\
        + C\left[ \E\left|\int_0^t\sum_{k=0}^{+\infty}\int_\OO |F_\lambda^{\prime\prime}(\varphi_{\lambda,n})||G_n(\varphi_{\lambda,n})[u_k]|^2\:\d x \: \d s\right|^{\frac{p}{2}} +\E \left| \int_0^t |\overline{\mu_{\lambda, n}(s)}|^2 \: \d s \right|^{\frac{p}{4}} \right]
    \end{multline}
    Finally, we observe that, given the growth condition Assumption \ref{hyp:potential}
    \[
    \begin{split}
        \left| |\OO|\overline{\mu_{\lambda,n}} \right| & = \left| (\mu_{\lambda,n}, 1)_H \right| \\
        & \leq |(a\varphi_{\lambda,n}, 1)_H| + |(K * \varphi_{\lambda,n}, 1)_H| + |(\Psi^\prime_\lambda(\varphi_{\lambda,n}), 1)_H| + \gamma(\varphi_{\lambda,n}, 1)_H \\
        & \leq C\left[ 1 + \|\varphi_{\lambda,n}\|_H + \int_\OO \Psi(J_\lambda(\varphi_{\lambda,n})) \: \d x\right] \\
        & \leq C\left[ 1 + \|\varphi_{\lambda,n}\|_H + \int_\OO \Psi_\lambda(\varphi_{\lambda,n}) \: \d x\right]
    \end{split}
    \]
    since the general theory (see \cite{barbu-monot}) implies $\Psi \circ J_\lambda \leq \Psi_\lambda$ in $\mathbb R$. Therefore, we arrive at
    \begin{equation} \label{eq:estimate_aux2}
        \E \supt |\overline{\mu_{\lambda,n}}(\tau)|^\frac p2 \leq C\left[1 +\E \supt \| \Psi_\lambda(\varphi_{\lambda,n}(\tau))\|_{L^1(\OO)}^\frac p2  \right].
    \end{equation}
    Hence, \review{multiplying \eqref{eq:estimate_aux2} for a sufficiently small constant, we can sum \eqref{eq:estimate_aux2} to \eqref{eq:estimate_aux} and} conclude, exploiting also the equality of laws preserved by $X_n$,
    \begin{multline} \label{eq:estimate_2}
        \tE\supt\mathcal \|\Psi_\lambda(\widetilde\varphi_{\lambda, n}(\tau))\|^\frac p2_{L^1(\OO)} + \tE \supt |\overline{\widetilde{\mu}_{\lambda,n}(\tau)}|^\frac p2 + \tE\left| \int_0^t \|\nabla\widetilde\mu_{\lambda,n}(\tau)\|_{\b H}^2 \: \d \tau \right|^\frac p2
        \\ \leq C\left[ 1+\tE \|\widetilde\varphi_{0}\|_H^{p} + \tE \|\Psi_\lambda(\review{\widetilde \varphi_{0,n}})\|_{L^1(\OO)}^\frac p2 + \tE \left| \int_0^t |\overline{\widetilde \mu_{\lambda, n}(s)}|^2 \: \d s \right|^{\frac{p}{4}} \right. \\ \left.
        + \tE\left|\int_0^t\sum_{k=0}^{+\infty}\int_\OO |F_\lambda^{\prime\prime}(\widetilde\varphi_{\lambda,n}(s))||G_n(\widetilde\varphi_{\lambda,n}(s))[u_k]|^2\:\d x \: \d s\right|^{\frac{p}{2}} \right].
    \end{multline}
    Now, we would like to let $n\to+\infty$ in \eqref{eq:estimate_2}. By weak lower semicontinuity of the norms, given the weak convergences proven in the previous subsection, it is immediate to pass to the limit to the second and third terms to the left hand side. As $\review{\widetilde\varphi_{0,n}} \to \review{\widetilde\varphi_{0}}$ in $L^\frac p2(\tom; H)$, it follows that $\Psi_\lambda(\review{\widetilde\varphi_{0,n}}) \rightarrow \Psi_\lambda(\review{\widetilde\varphi_{0}})$ in $L^\frac p2(\tom;L^1(\OO))$ as $n\rightarrow+\infty$ by the dominated convergence theorem. Moreover, due to the strong convergences 
    \begin{align*}
        G_n(\widetilde{\varphi}_{\lambda,n}) \rightarrow G(\widetilde{\varphi}_{\lambda})&\quad \text{ in }L^q(\tom;L^2(0,T;\cL^2(U,H))),\\
        \widetilde{\varphi}_{\lambda,n} \rightarrow \widetilde{\varphi}_\lambda&\quad\text{ in }L^q(\tom;L^2(0,T;H))
    \end{align*}
    for all $q<p$, we deduce \review{that, up to subsequences that we do not relabel for readability purposes,}
    \begin{align*}
        |F_\lambda^{\prime\prime}(\widetilde{\varphi}_{\lambda,n})||G_n(\widetilde{\varphi}_{\lambda,n})[u_k]|^2 \rightarrow |F_\lambda^{\prime\prime}(\widetilde{\varphi}_{\lambda})||G(\widetilde{\varphi}_{\lambda})[u_k]|^2\quad \text{a.e. in }\tom\times(0,T)\times \OO
    \end{align*}
    for all $k\in\enne$. Noticing that $|F^{\prime\prime}_\lambda| \leq C_\lambda$ for some constant $C_\lambda>0$, we infer that 
    \begin{align*}
        \sum_{k=0}^{+\infty}\int_\OO |F_\lambda^{\prime\prime}(\widetilde{\varphi}_{\lambda,n})||G_n(\widetilde{\varphi}_{\lambda,n})[u_k]|^2\d x \leq C_\lambda\|G_n(\widetilde{\varphi}_{\lambda,n})\|_{\cL^2(U,H)}^2 \leq C_\lambda
    \end{align*}
    owing to Assumption \ref{hyp:G}. Hence, given \review{Lemma \ref{lem:first_estimate_lam}}, we can apply the dominated convergence theorem to conclude 
    \begin{align*}
        \widetilde{\E}\left(\int_0^t\sum_{k=0}^{+\infty}\int_\OO |F_\lambda^{\prime\prime}(\widetilde{\varphi}_{\lambda,n})||G_n(\widetilde{\varphi}_{\lambda,n})[u_k]|^2\:\d x \: \d s\right)^{\frac{p}{2}} \rightarrow \widetilde{\E}\left(\int_0^t\sum_{k=0}^{+\infty}\int_\OO |F_\lambda^{\prime\prime}(\widetilde{\varphi}_{\lambda})||G(\widetilde{\varphi}_{\lambda})[u_k]|^2\:\d x \: \d s\right)^{\frac{p}{2}}
    \end{align*}
    as $n\rightarrow+\infty$. We are now in a position to pass to the limit as $n \to +\infty$ in \eqref{eq:estimate_2}, yielding
    {\small \begin{multline} \label{eq:estimate_22}
        \tE\supt\mathcal \|\Psi_\lambda(\widetilde\varphi_{\lambda}(\tau))\|^\frac p2_{L^1(\OO)} + \tE \supt |\overline{\widetilde \mu_{\lambda}}(\tau)|^\frac p2 + \tE\left| \int_0^t \|\nabla\widetilde\mu_{\lambda}(\tau)\|_{\b H}^2 \: \d \tau \right|^\frac p2 
        \\ \leq C\left[ 1+\tE \|\widetilde\varphi_{0}\|_H^{p} + \tE \|\Psi_\lambda(\widetilde \varphi_{0})\|_{L^1(\OO)}^\frac p2 + \tE \left| \int_0^t |\overline{\widetilde\mu_{\lambda}(\tau)}|^2 \: \d \tau\right|^\frac p4 \right. \\ \left. + \tE\left|\int_0^t\sum_{k=0}^{+\infty}\int_\OO |F_\lambda^{\prime\prime}(\widetilde\varphi_{\lambda}(s))||G(\widetilde\varphi_{\lambda}(s))[u_k]|^2\:\d x \: \d s\right|^{\frac{p}{2}} \right],
    \end{multline}}where $C > 0$ is independent of $\lambda$. Finally, we observe that since the Yosida approximation is monotonically convergent to $\Psi$ from below, we have
    \[
    \tE \|\Psi_\lambda(\widetilde \varphi_{0})\|_{L^1(\OO)}^\frac p2 \leq \tE \|\Psi(\widetilde \varphi_{0})\|_{L^1(\OO)}^\frac p2 = \E \|\Psi(\varphi_{0})\|_{L^1(\OO)}^\frac p2
    \]
    for any $\lambda > 0$. The last term makes use of the fact that the approximation scheme satisfies the crucial property (thanks to Assumption \ref{hyp:potential})
    \[
    \begin{split}
         |F''_\lambda(r)| & = |\Psi''_\lambda(r)+\gamma| \\
         & = |\Psi''(J_\lambda(r))J'_\lambda(r) + \gamma| \\
         & \leq C\left( 1 + \Psi(J_\lambda(r)) \right) \\
         & \leq C\left( 1 + \Psi_\lambda(r)) \right)
    \end{split}
    \]
    for all $r \in \mathbb R$. The H\"{o}lder inequality then implies
    \[ 
    \begin{split}
        & \tE\left|\int_0^t\sum_{k=0}^{+\infty}\int_\OO |F_\lambda^{\prime\prime}(\widetilde\varphi_{\lambda}(s))||G(\widetilde\varphi_{\lambda}(s))[u_k]|^2\:\d x \: \d s\right|^{\frac{p}{2}} \\
        & \hspace{1cm}\leq C\left[ 1  + \tE\left|\int_0^t\sum_{k=0}^{+\infty}\int_\OO |\Psi_\lambda(\widetilde\varphi_{\lambda}(s))||G(\widetilde\varphi_{\lambda}(s))[u_k]|^2\:\d x \: \d s\right|^{\frac{p}{2}} \right] \\
        & \hspace{1cm}\leq C\left[ 1  + \tE\left|\int_0^t\|\Psi_\lambda(\widetilde\varphi_{\lambda}(s))\|_{L^1(\OO)}\sum_{k=0}^{+\infty} \|G(\widetilde\varphi_{\lambda}(s))[u_k]\|_{L^\infty(\OO)}^2 \: \d s\right|^{\frac{p}{2}} \right] \\
        & \hspace{1cm}\leq C\left[ 1  + \tE\left|\int_0^t\|\Psi_\lambda(\widetilde\varphi_{\lambda}(s))\|_{L^1(\OO)} \: \d s\right|^{\frac{p}{2}} \right].
    \end{split}
    \]
    Injecting in \eqref{eq:estimate_22} gives
    \begin{multline*}
        \tE\supt\mathcal \|\Psi_\lambda(\widetilde\varphi_{\lambda}(\tau))\|^\frac p2_{L^1(\OO)} + \tE \supt |\overline{\widetilde \mu_{\lambda}}(\tau)|^\frac p2 + \tE\left| \int_0^t \|\nabla\widetilde\mu_{\lambda}(\tau)\|_{\b H}^2 \: \d \tau \right|^\frac p2 
        \\ \leq C\left[ 1+\tE \|\widetilde\varphi_{0}\|_H^{p} + \tE \|\Psi_\lambda(\widetilde \varphi_{0})\|_{L^1(\OO)}^\frac p2 + \tE \left| \int_0^t |\overline{\widetilde\mu_{\lambda}(\tau)}|^2 \: \d \tau\right|^\frac p4  + \tE\left|\int_0^t\|\Psi_\lambda(\widetilde\varphi_{\lambda}(s))\|_{L^1(\OO)} \: \d s\right|^{\frac{p}{2}}  \right],
    \end{multline*}
    and the Gronwall lemma gives the second estimate
    \begin{align*}\label{eq:secondestimate2}
        \|\widetilde \mu_{\lambda}\|_{L^\frac p2(\tom;L^2(0,T;V))} + \|\nabla\widetilde \mu_{\lambda}\|_{L^p(\tom;L^2(0,T;\b H))}+\|\Psi_\lambda(\widetilde\varphi_\lambda))\|_{L^\frac p2(\tom;L^\infty(0,T;L^1(\OO)))} \leq C_2,
    \end{align*}
    and the proof is complete.
    \end{proof} \noindent
    The \review{third and fourth} estimates play out very similarly with respect to what is carried out in the \review{proofs of Lemmas \ref{lem:third_estimate} and \ref{lem:fourth_estimate}}. Indeed, they are essentially corollaries of the first two estimates. For the sake of brevity, we do not report the full computations, but only highlight the results \review{in the following lemmas.}
    \begin{lem} \label{lem:third_estimate_lam}
    \review{
        It holds that
        \begin{enumerate}[(i)]
            \item the family $\{G(\widetilde\varphi_{\lambda})\}_{\lambda > 0}$ is uniformly bounded in 
            \[
            L^\infty(\Omega\times(0,T);\cL^2(U,H))\cap L^p(\Omega;L^2(0,T;\cL^2(U,V)))
            \]
            for all $p \geq 2$;
            \item the family $\left\{\displaystyle\int_0^. G(\widetilde \varphi_{\lambda}(s)) \: \d W(s)\right\}_{\lambda > 0}$ is uniformly bounded in
            \[
            L^q(\Omega;W^{k,q}(0,T;H))\cap L^2(\Omega;H^{k}(0,T;V))
            \]
            for all $k \in (0,\frac 12)$ and $q \geq 2$.
        \end{enumerate}
        }
    \end{lem} \noindent
    Exploiting the weak formulation of the problem (passing to the limit of $n\to+\infty$ therein, owing to the proven convergences) we eventually also show that
    \begin{lem} 
        \review{The family $\{\widetilde \varphi_{\lambda}\}_{\lambda > 0}$ is uniformly bounded in $L^p(\tom;W^{\alpha,p}(0,T;V^*))$ for all $\alpha \in (0, \frac 12)$ and $p \geq 2$.}
    \end{lem} \noindent
    \review{Finally,} let us provide an additional estimate, since in order to pass to the limit as $\lambda \to 0^+$ we can no longer rely on the Lipschitz continuity of $\Psi'_\lambda$, as the Lipschitz constants are not uniformly bounded in $\lambda$. 
    \begin{lem} \label{lem:fifth_estimate_lam}
        \review{The family $\{F^\prime_\lambda(\widetilde{\varphi}_\lambda)\}_{\lambda > 0}$ is uniformly bounded in $L^{\frac{p}{2}}(\tom;L^2(0,T;H))$ for all $p \geq 2$.}
    \end{lem}
    \begin{proof}
        Testing the equation for $\widetilde \mu_\lambda$ by $\Psi_\lambda(\widetilde{\varphi}_\lambda) = F_\lambda^\prime(\widetilde{\varphi}_\lambda) + \gamma\widetilde{\varphi}_\lambda$, we have
     \begin{multline*}
         \int_\OO\widetilde{\mu}_\lambda F_\lambda^\prime(\widetilde{\varphi}_\lambda) \: \d x + \gamma\int_\OO\widetilde{\mu}_\lambda \widetilde{\varphi}_\lambda \: \d x \\ =\int_\OO a \review{\widetilde \varphi_\lambda}F_\lambda'(\widetilde \varphi_\lambda)\: \d x + \review{\gamma}\int_\OO a\review{\widetilde \varphi_\lambda^2}\: \d x - \int_\OO (K*\widetilde{\varphi}_\lambda)F^\prime_\lambda(\widetilde{\varphi}_\lambda) \:\d x  - \gamma \int_\OO (K*\widetilde{\varphi}_\lambda)\widetilde{\varphi}_\lambda \:\d x \\+ \int_\OO|F^\prime_\lambda(\widetilde{\varphi}_\lambda)|^2\d x + \gamma\int_\OO \widetilde \varphi_\lambda F'_\lambda(\widetilde \varphi_\lambda) \: \d x.
     \end{multline*}
     Rearranging terms, we can use Young's inequality multiple times to conclude
     \begin{align*}
    \|F_\lambda^\prime(\widetilde{\varphi}_\lambda)\|_{L^2(0,T;H)}^2 \leq C\left[ 1 + \|\widetilde{\mu}_\lambda\|_{L^2(0,T;H)}^2 + \|\widetilde{\varphi}_\lambda\|_{L^2(0,T;H)}^2\right],
     \end{align*}
     for a constant $C >0$ independent of $\lambda$.
     Therefore, thanks to \review{Lemmas \ref{lem:first_estimate_lam} and \ref{lem:second_estimate_lam}}, we infer that
     \begin{align*}\label{Est:lam_5}
    \|F^\prime_\lambda(\widetilde{\varphi}_\lambda)\|_{L^{\frac{p}{2}}(\tom;L^2(0,T;H))}\leq C_5
     \end{align*}
     for some $C_5>0$ independent of $\lambda > 0$.
     \end{proof}
    \subsection{Passage to the limit as $\lambda \to 0^+$}
    \label{ssec:lim_lam}
    The previous estimates allow us to pass to the limit as $\lambda \to 0^+$, by means of a second stochastic compactness argument. From this subsection onwards, we assume to fix once and for all a suitable vanishing sequence $\{\lambda_n\}_{n \in \mathbb N}$ and we understand without relabeling that the passage to the limit is performed along said sequence. The sequence of laws of $\{\varphi_{\lambda}\}_{\lambda > 0}$ is tight in the same path space illustrated in Subsection \ref{ssec:lim_n}. To this end, we still make use of the compact embeddings (see \cite[Theorem 2.1]{fland-gat})
    \begin{align*}
        L^2(0,T;V)\cap H^\alpha(0,T;V^*) & \hookrightarrow L^2(0,T;H) \\
        L^\infty(0,T;H) \cap W^{\alpha, p}(0,T;V^*) &\hookrightarrow C^0([0,T]; H^{-\delta}(\OO))
    \end{align*}
    still assuming that $\alpha p > 1$, i.e., that $p > 2$, and for some small $\delta > 0$. The following tightness lemmas hold; their proofs closely follow the corresponding ones in Lemmas \ref{lem:tightness}-\ref{lem:tightness4}.
    \begin{lem} \label{lem:tightness5}
        If $p > 2$, then the family of laws of $\{\widetilde \varphi_{\lambda}\}_{\lambda > 0}$ is tight in the space $\mathcal X := L^2(0,T;H) \cap C^0([0,T]; H^{-\delta}(\OO))$ for any $\delta \in (0, \frac 12)$.
    \end{lem}
    \begin{lem} \label{lem:tightness6}
        The family of laws of $\{G(\widetilde \varphi_\lambda) \cdot \widetilde W\}_{\lambda > 0}$ is tight in the space $C^0([0,T]; V^*)$.
    \end{lem} \noindent
    By identifying $\widetilde W$ with a constant sequence of random variables $\{\widetilde W_\lambda\}_{\lambda > 0}$ with values in $C^0([0,T];U_0)$, we also get that 
    \begin{lem} \label{lem:tightness7}
        The family of laws of $\{\widetilde W_\lambda\}_{\lambda > 0}$ is tight in the space $C^0([0,T]; U_0)$.
    \end{lem} \noindent
    By the same token, we identify $\widetilde\varphi_0$ with a constant sequence $\{\widetilde \varphi_{0,\lambda}\}_{\lambda >0}$ that is again tight in $V^*$ as $H \hookrightarrow V^*$ compactly. Namely, we have
    \begin{lem} \label{lem:tightness8}
        The family of laws of $\{\widetilde\varphi_{0,\lambda}\}_{\lambda > 0}$ is tight in $V^*$.
    \end{lem} \noindent
    Collecting Lemmas \ref{lem:tightness5}-\ref{lem:tightness8}, we obtain that the family of laws 
    \begin{align*}
        \{\widetilde \varphi_{\lambda},\, G(\widetilde \varphi_{\lambda})\cdot W,\, W_\lambda,\, \varphi_{0,\lambda}\}_{\lambda >0}
    \end{align*}
    is tight in the space
    \begin{align*}
        \mathcal X \times C^0([0,T]; V^*) \times C^0([0,T];U_0) \times V^*.
    \end{align*}
    Owing to the Prokhorov and Skorokhod theorems (see \cite[Theorem 2.7]{ike-wata} and \cite[Theorem 1.10.4, Addendum 1.10.5]{vaa-well}), there exists a probability space $(\hom, \hF,\hP)$ and a family of random variables $Y_\lambda:  (\hom, \hF)\to(\tom, \tF)$ such that the law of $Y_\lambda$ is $\tP$ for every $\lambda \in (0,1)$, namely $\hP \circ Y_\lambda^{-1} = \tP$ (so that composition with $Y_\lambda$ preserves laws), and, \review{thanks to the Vitali and Banach--Alaoglu theorems and reflexivity}, the following convergences hold as $\lambda \to 0^+$:
    \begin{align*}
        \widehat{\varphi}_{\lambda}:= \widetilde{\varphi}_\lambda \circ Y_\lambda \to \widehat{\varphi}
		&\quad\text{ in }L^q(\hom;L^2(0,T;H) \cap C^0([0,T]; H^{-\delta}(\OO))\;\;\forall \: q<p, \\
        \widehat{\varphi}_{\lambda} \overset{*}{\rightharpoonup} \widehat{\varphi}&\quad\text{ in }L^p(\hom;H^\alpha(0,T;V^*)),\\
        \widehat{\varphi}_{\lambda} \rightharpoonup \widehat{\varphi}&\quad\text{ in }L^p(\hom;L^2(0,T;V)), \\
        \widehat{\mu}_{\lambda} := \widetilde
        \mu_{\lambda}\circ Y_\lambda \rightharpoonup \widehat{\mu}&\quad\text{ in }L^{\frac{p}{2}}(\hom;L^2(0,T;V)), \\
        \nabla\widehat{\mu}_{\lambda} \rightharpoonup \nabla\widehat{\mu}&\quad\text{ in }L^p(\hom;L^2(0,T;\b H)), \\
        \widehat{I}_\lambda := (G(\widetilde{\varphi}_\lambda)\cdot \widetilde{W}) \circ Y_\lambda \rightarrow \widehat{I}&\quad\text{ in }L^q(\hom;C^0([0,T];H))\;\;\forall\: q<p,\\
        \widehat{W}_\lambda := \widetilde{W}\circ Y_\lambda \rightarrow \widehat{W}&\quad\text{ in }L^q(\hom;C^0([0,T];U))\;\;\forall\: q<p, 
    \end{align*}
     for some limiting measurable processes enjoying the following regularity properties:
    \begin{align*}
        \widehat{\varphi}&\in L^p_w(\hom;L^\infty(0,T;H))\cap L^p(\hom;H^\alpha(0,T;V^*))\cap L^p(\hom;L^2(0,T;V) \cap C^0([0,T]; H^{-\delta}(\OO)),\\
        \widehat{\mu}&\in L^\frac p2(\hom;L^2(0,T;V)),\\
        \nabla \widehat{\mu}&\in L^\review{p}(\hom;L^2(0,T;\b H)),\\
        \widehat{I} &\in L^p(\hom;C^0([0,T];V^*)), \\
        \widehat{W}&\in L^p(\hom;C^0([0,T];U)), \\
        \widehat{\varphi}_0 & \in L^p(\hom; H).
    \end{align*}
    Finally, we address once again the nonlinearities of the problem. Since the space $L^\frac p2(\hom;L^2(0,T;H))$ is reflexive for $p > 2$, the uniform bound given in \review{Lemma \ref{lem:fifth_estimate_lam}} yields
\[
F'_\lambda(\widehat \varphi_\lambda) \rightharpoonup \xi \quad \text{in } L^\frac p2(\hom;L^2(0,T;H)).
\]
Let us show that indeed $\xi = \Psi'(\widehat \varphi) -\gamma\widehat \varphi = \review{F'}(\widehat \varphi)$. Using the strong convergence of
    \[
\widehat{\varphi}_{\lambda} \rightarrow \widehat{\varphi}\quad\text{ in }L^q(\hom;L^2(0,T;H))
\]
for $q < p$, by comparison it is immediate to deduce that
\[
\Psi'_\lambda(\widehat \varphi_\lambda) \rightharpoonup \xi + \gamma \widehat \varphi \quad \text{in } L^\frac p2(\hom;L^2(0,T;H)).
\]
As by definition it holds $\Psi'_\lambda = \Psi' \circ J_\lambda$, we can investigate the sequence $\{J_\lambda(\widehat \varphi_\lambda)\}_{\lambda > 0}$, namely
\[
\begin{split}
    \|J_\lambda(\widehat \varphi_\lambda(s)) - \widehat \varphi(s)\|^2_H & \leq 2\|J_\lambda(\widehat \varphi_\lambda(s))-\widehat \varphi_\lambda(s)\|^2_H + 2\|\widehat \varphi_\lambda(s)-\widehat \varphi(s)\|^2_H \\
    & = 2\lambda^2\|\Psi'_\lambda(\widehat \varphi_\lambda(s))\|^2_H + 2\|\widehat \varphi_\lambda(s)-\widehat \varphi(s)\|^2_H
\end{split}
\]
implying that
\[
J_\lambda(\hphi_\lambda) \to \hphi \quad \text{in } L^\frac p2(\hom;L^2(0,T;H)),
\]
thanks to \review{Lemma \ref{lem:second_estimate_lam}}. Then, the strong-weak closure of maximal monotone operators (see \cite[Proposition 2.1]{barbu-monot} implies $\xi + \gamma \widehat \varphi= \Psi'(\widehat \varphi)$, as claimed. Given \ref{hyp:G}, we immediately have
    \begin{equation*}
        \|G(\widehat{\varphi}_{\lambda}) - G(\widehat{\varphi})\|_{L^q(\hom;L^2(0,T;\cL^2(U,H)))} \to 0
    \end{equation*}
    as $\lambda \to 0^+$, it follows that 
    \begin{align*}
        G(\widehat{\varphi}_{\lambda}) \rightarrow G(\widehat{\varphi})\quad \text{ in }L^q(\hom;L^2(0,T;\cL^2(U,H)))
    \end{align*}
    for all $q<p$. 
    \subsection{Identification of the limit solution.} This is essentially a \textit{verbatim} iteration of the argument presented in Subsection \ref{ssec:limit_n}, and is therefore omitted. However, it completes the proof of existence by retrieving a martingale solution to \eqref{eq:strongNCH}.
    \section{\review{Pathwise uniqueness and existence of probabilistically-strong solutions}} \label{sec:proof2}
    The present section is devoted to proving Theorem \ref{thm:uniqueness} and Corollary \ref{cor:existenceS}, i.e., to showing a continuous dependence estimate for problem \eqref{eq:strongNCH}, entailing pathwise uniqueness in dimension two and three and, in turn, existence and uniqueness of a probabilistically-strong solution for the stochastic nonlocal Cahn--Hilliard equation. For $i \in \{1,2\}$, let $\varphi_{0i}$ comply with the assumptions listed in Definition \ref{def:mart_sol}. By Theorem \ref{thm:existenceM}, the formal problem
    \begin{equation} \label{eq:CDeq}
        \begin{cases}
        \d \hphi_i - \Delta \hmu_i = G(\hphi_i) \: \d W & \quad \text{in }(0,T) \times \OO, \\
        \hmu_i = a\hphi_i - K * \hphi_i + F'(\hphi_i) & \quad \text{in }(0,T) \times \OO, \\
        \partial_{\bn} \hmu_i = 0 & \quad \text{on }(0,T) \times \partial \OO, \\
        \hphi(0) = \varphi_{0i} & \quad \text{in } \OO,
    \end{cases}
    \end{equation}
    admits at least one martingale solution $\hphi_i$ for $i \in \{1,2\}$, defined on a suitable filtered probability space and satisfying \eqref{eq:CDeq} with respect to some cylindrical Wiener process on $U$. Assuming that $\hphi_1$ and $\hphi_2$ are defined on the same stochastic basis
    \[
    \left\{ \hom,\, \hF,\, (\hF_t)_{t \in [0,T]},\, \hP \right\}
    \]
    and satisfy \eqref{eq:CDeq} with respect to the same cylindrical Wiener process $\hW$ on $U$, we define the differences
    \[
    \hphi:= \hphi_1 -\hphi_2, \qquad \hmu := \hmu_1-\hmu_2, \qquad \hphi_0 := \hphi_{1}(0) - \hphi_{2}(0).
    \]
    Under the additional assumption \ref{hyp:additional}, it is straightforward to see that the process is pathwise mass-conservative, as one can easily see by plugging a constant nonzero function in the variational formulation of the problem (see Definition \ref{def:mart_sol}). In particular, it is immediate to observe that
    \[
    \overline{\varphi_{01}} =\overline{\varphi_{02}} \Rightarrow \overline{\widehat \varphi(t)} = 0
    \]
    for all $t \geq 0$, $\P$-almost surely. Applying the It\^{o} lemma for twice Fréchet-differentiable functionals (see \cite[Theorem 4.32]{dapratozab}) to
    \[
    \Gamma: V^* \to \R, \qquad \Gamma(u) = \dfrac 12 \|\nabla \mathcal{N} (u -|\OO|^{-1}\langle u, 1\rangle_{V^*,V})\|^2_\b H,
    \]
    again, evaluated at $u = \hphi(t)$. This results in
    \begin{multline}
        \label{eq:unique20}
        \dfrac 12 \|\nabla \mathcal{N} \hphi(t)\|^2_\b H + \int_0^t (\hphi(\tau),\,\hmu(\tau))_H\: \d\tau  = \dfrac 12 \|\nabla \mathcal{N} \hphi_0\|^2_\b H \\ + \dfrac 12 \int_0^t \left\|\nabla \mathcal N \left[ G(\hphi_1(\tau))-G(\hphi_2(\tau)) - \left( \overline{G(\hphi_1(\tau))} - \overline{G(\hphi_2(\tau))} \right) \right]\right\|^2_{\cL^2(U, \b H)}  \: \d \tau
        \\ + \int_0^t \left( \hphi(\tau), \, \mathcal N \left[ G(\hphi_1(\tau))-G(\hphi_2(\tau)) - \left( \overline{G(\hphi_1(\tau))} - \overline{G(\hphi_2(\tau))} \right) \right] \: \d \hW \right)_H.
    \end{multline}
    Using Assumption \ref{hyp:potential} and the definition of $\hmu$, we get
    \begin{equation} \label{eq:unique21}
        \begin{split}
             & \int_0^t (\hphi(\tau),\,\hmu(\tau))_H\: \d\tau \\ 
             & \hspace{1cm} =  \int_0^t (\hphi(\tau),\,a\hphi(\tau))_H \: \d \tau - \int_0^t (\hphi(\tau),\,K * \hphi(\tau))_H\: \d\tau \\
             & \hspace{3cm} + \int_0^t (\hphi(\tau),\,F'(\hphi_1(\tau))-F'(\hphi_2(\tau)))_H\: \d\tau \\
             & \hspace{1cm} \geq  C_0\int_0^t\|\hphi(\tau)\|^2_H \: \d \tau + \int_0^t (\nabla \mathcal N\hphi(\tau),\,\nabla K * \hphi(\tau))_\b H\: \d\tau.
        \end{split}
    \end{equation}
    Moreover, we have
    \begin{equation} \label{eq:unique22} 
            \left|\int_0^t (\nabla \mathcal N\hphi(\tau),\,\nabla K * \hphi(\tau))_\bH\: \d\tau \right| 
            \leq \dfrac{C_0}{2}\int_0^t \|\hphi(\tau)\|_H^2 \: \d \tau  + \dfrac{\|\nabla K\|_{\b L^1(B_R)}^2}{2C_0} \int_0^t \|\nabla \mathcal N\hphi(\tau)\|_{\bH}^2 \: \d \tau.
    \end{equation}
    Recalling \cite[Proposition 2.1]{DPGS2}, we observe that
    \[
    \begin{split}
    \|G(\hphi_1(\tau)) - G(\hphi_2(\tau))\|^2_{\cL^2(U,V^*)}  & = \sum_{k=0}^{+\infty} \|G(\hphi_1(\tau))[u_k] - G(\hphi_2(\tau))[u_k]\|^2_{V^*} \\
    & \geq C\sum_{k=0}^{+\infty} \left\|\nabla \mathcal N \left[ G(\hphi_1(\tau))[u_k] - G(\hphi_2(\tau))[u_k] \right]\right\|^2_{\b H}
    \end{split}
    \]
    for a constant $C > 0$ only depending on $\OO$, and therefore, by possibly updating the value of $C$, 
    \begin{equation} \label{eq:unique23}
        \begin{split}
            \int_0^t \left\|\nabla \mathcal N \left[ G(\hphi_1(\tau))-G(\hphi_2(\tau))  \right]\right\|^2_{\cL^2(U, \bH)}  \: \d \tau & \leq C\int_0^t \|G(\hphi_1(\tau)) - G(\hphi_2(\tau))\|^2_{V^*} \: \d \tau \\
            & \leq CL_G^2\int_0^t \|\hphi(\tau)\|^2_{V^*} \: \d\tau
        \end{split}
    \end{equation}
    for some constant $C > 0$ only depending on $|\OO|$. Using \eqref{eq:unique21}-\eqref{eq:unique23} in \eqref{eq:unique20} yields
    \begin{multline}
        \label{eq:unique24}
        \dfrac 12 \|\nabla \mathcal{N} \hphi(t)\|^2_\b H + \dfrac{C_0}{2} \int_0^t \|\hphi(\tau)\|_H^2\: \d\tau  \leq \dfrac 12 \|\nabla \mathcal{N} \hphi_0\|^2_\b H +C\widehat L_G^2\int_0^t \|\hphi(\tau)\|^2_{V^*} \: \d\tau \\+ \dfrac{\|\nabla K\|_{\b L^1(B_R)}^2}{2C_0} \int_0^t \|\nabla \mathcal N\hphi(\tau)\|_{\bH}^2 \: \d \tau
        + \int_0^t \left( \hphi(\tau), \, \mathcal N \left[ G(\hphi_1(\tau))-G(\hphi_2(\tau))\right] \: \d \hW \right)_H.
    \end{multline}
    Raising \eqref{eq:unique24} to the power $\frac p2$, taking supremums in time and $\hP$-expectations leads to
    \begin{multline*}
        \hE \supt \|\hphi(\tau)\|^p_{V^*} + \hE \left| \int_0^t \|\hphi(\tau)\|^2_H \: \d \tau \right|^\frac p2 \leq C\left[ \|\hphi_0\|^p_{V^*} + \hE \left| \int_0^t \|\hphi(\tau)\|^2_{V^*} \: \d \tau \right|^\frac p2 \right. \\
        \left. + \hE \supt \left| \int_0^\tau \left( \hphi(s), \, \mathcal N \left[ G(\hphi_1(s))-G(\hphi_2(s))  \right] \: \d \hW(s) \right)_H \right|^\frac p2 \right],
    \end{multline*}
    where the constant $C > 0$ now also depends on $p$. We are then left with estimating the stochastic integrals on the right hand side of the inequality above. This is done by exploiting the Burkholder--Davis--Gundy inequality. Indeed, recalling \eqref{eq:unique23} (see also \cite[Section 4]{DPGS})
    \begin{equation*} \small
        \begin{split}
        & \hE \supt \left| \int_0^\tau \left( \hphi(s), \, \mathcal N \left[ G(\hphi_1(s))-G(\hphi_2(s))  \right] \: \d \hW(s) \right)_H \right|^\frac p2 \\
        &\hspace{0.3cm} \leq C\hE \left| \int_0^t \|\nabla \mathcal N\hphi(\tau)\|_\bH^2\left\| \nabla \mathcal N \left[ G(\hphi_1(\tau))-G(\hphi_2(\tau))  \right]\right\|^2_{\cL^2(U, \b H)} \: \d\tau\right|^\frac p4 \\
        &\hspace{0.3cm} \leq C\hE \left[ \supt \|\nabla \mathcal N \hphi(\tau)\|_\bH^\frac p2  \left| \int_0^t \| \hphi(\tau)\|^2_{V^*} \: \d\tau\right|^\frac p4 \right] \\
        & \hspace{0.3cm} \leq \dfrac{1}{4}\hE \supt \|\hphi(\tau)\|_{V^*}^p + C\hE \left| \int_0^t \|\hphi(\tau)\|^2_{V^*} \: \d \tau \right|^\frac p2.
        \end{split}
    \end{equation*}
    In light of all of the above, we conclude that
    \begin{equation} \label{eq:unique40}
        \hE \supt \|\hphi(\tau)\|^p_{V^*} + \hE \left| \int_0^t \|\hphi(\tau)\|^2_H \: \d \tau \right|^\frac p2 \leq C\left[ \|\hphi_0\|^p_{V^*} + \hE \left| \int_0^t \|\hphi(\tau)\|^2_{V^*} \: \d \tau \right|^\frac p2 \right],
    \end{equation}
    for a constant $C > 0$ only depending on the structural parameters of the problem, including $p$ and $T$. Therefore, \eqref{eq:unique40} and the Gronwall lemma imply that
    \begin{equation*}
        \|\hphi\|^p_{L^p(\hom; L^\infty(0,T;V^*))} + \|\hphi\|^p_{L^p(\hom; L^2(0,T;H))} \leq C\|\hphi_0\|^p_{L^p(\hom;V^*)},
    \end{equation*}
    and choosing $\varphi_{01} = \varphi_{02}$ in \eqref{eq:CDeq} leads to pathwise uniqueness, and the proof is complete.
    \section{\review{Nonlocal-to-local asymptotics: potentials of polynomial growth}} \label{sec:convergence_unique}
    Hereafter, we show the proof of Theorem \ref{thm:convergence_unique}, i.e., that probabilistically-strong solutions to the stochastic nonlocal Cahn--Hilliard problem \eqref{eq:NCH:eps} converge to the probabilistically-strong solution to the local problem \eqref{eq:localCHstrong} as $\varepsilon \to 0^+$. In particular, for any $\eps > 0$, the process $\varphi_\varepsilon$ solves the nonlocal problem \eqref{eq:NCH:eps}, that for convenience we recall here:
    \begin{equation} \label{eq:base_nonlocal}
    \begin{cases}
        \d \varphi_\eps - \Delta \mu_\eps \: \d t = G(\varphi_\eps) \: \d W & \quad \text{in } (0,T) \times \OO, \\
         \mu_\eps = \mathcal{L}_\eps\varphi_\eps + F'(\varphi_\eps)& \quad \text{in } (0,T) \times \OO, \\
         \partial_{\bn} \mu_\eps = 0 & \quad \text{on } (0,T) \times \partial\OO, \\
         \varphi_\eps(0) = \varphi_{0,\eps} & \quad \text{in } \OO.
    \end{cases}
    \end{equation}
    On the other hand, the process $\varphi$ solves the local Cahn--Hilliard system \eqref{eq:localCHstrong}, reading
    \begin{equation} \label{eq:base_local}
        \begin{cases} 
        \d \varphi - \Delta \mu \: \d t = G(\varphi) \: \d W & \quad \text{in } (0,T) \times \OO, \\
         \mu = -\Delta\varphi+ F'(\varphi)& \quad \text{in } (0,T) \times \OO, \\
         \partial_{\bn} \mu = \partial_\bn \varphi=0 & \quad \text{on } (0,T) \times \partial\OO, \\
         \varphi(0) = \varphi_0 & \quad \text{in } \OO.
    \end{cases}
    \end{equation}
    Since both problems admit a unique probabilistically-strong solution, we fix once and for all a filtered probability space 
    \[
    \left\{ \Omega,\, \cF,\, (\cF_t)_{t \in [0,T]},\, \P \right\}
    \]
    and, accordingly, a cylindrical Wiener process $W$ on $U$ to formulate \eqref{eq:base_nonlocal} and \eqref{eq:base_local}. As in the previous section, we define the differences
    \[
    \xi:= \varphi_\eps - \varphi, \qquad \eta := \mu_\eps-\mu, \qquad \xi_0 := \varphi_{0,\eps}-\varphi_{0}.
    \]
    Once again, Assumption \ref{hyp:additional} yields the key conservation property
    \[
    \overline{\xi(t)} = 0
    \]
    for all $t \geq 0$, $\P$-almost surely and for all $\varepsilon > 0$. Applying the It\^{o} lemma for twice Fréchet-differentiable functionals (see \cite[Theorem 4.32]{dapratozab}) to
    \[
    \Gamma: V^* \to \R, \qquad \Gamma(u) = \dfrac 12 \|\nabla \mathcal{N} (u -|\OO|^{-1}\langle u, 1\rangle_{V^*,V})\|^2_\b H,
    \]
    evaluated at $u = \xi(t)$ results in
    \begin{multline}
        \label{eq:unique50}
        \dfrac 12 \|\nabla \mathcal{N} \xi(t)\|^2_\b H + \int_0^t (\xi(\tau),\,\eta(\tau))_H\: \d\tau  = \dfrac 12 \|\nabla \mathcal{N} \xi_0\|^2_\b H \\ + \dfrac 12 \int_0^t \left\|\nabla \mathcal N \left[ G(\varphi_\eps(\tau))-G(\varphi(\tau)) \right]\right\|^2_{\cL^2(U, \b H)}  \: \d \tau
         + \int_0^t \left( \xi(\tau), \, \mathcal  N \left[ G(\varphi_\eps(\tau))-G(\varphi(\tau)) \right] \: \d W \right)_H.
    \end{multline}
    Expanding the scalar product to the left hand side, by summing and subtracting $\mathcal L_\eps \varphi$ we get
    \begin{multline} \label{eq:unique51}
             \int_0^t (\xi(\tau),\,\eta(\tau))_H\: \d\tau  =  \int_0^t (\xi(\tau),\,\mathcal L_\varepsilon\varphi(\tau)+\Delta  \varphi(\tau))_H \: \d \tau + \int_0^t (\xi(\tau),\,\mathcal L_\varepsilon\xi(\tau)) \: \d \tau\\
              +\int_0^t (\xi(\tau),\,F'(\varphi_\varepsilon(\tau))-F'(\varphi(\tau)))_H\: \d\tau. 
    \end{multline}
    On the right hand side of \eqref{eq:unique51}, we control the three terms as follows. We exploit \cite[Corollary 4.2]{abels2023strong} and the Young inequality to infer that
    \begin{equation} \label{eq:unique52}
        \begin{split}
            \left|\int_0^t (\xi(\tau),\,\mathcal L_\varepsilon\varphi(\tau)+\Delta  \varphi(\tau))_H \: \d \tau\right| \leq \frac 12 \int_0^t \|\xi(\tau)\|^2_H \: \d \tau +C\varepsilon\int_0^t\|\varphi(\tau)\|^2_{H^3(\OO)} \: \d \tau,
        \end{split}
    \end{equation}
    while, given the positivity of the kernel $K_\eps$, we can apply the nonlocal interpolation inequality \cite[Lemma 3.4]{Davoli1}
    \begin{equation} \label{eq:unique53}
        \int_0^t (\xi(\tau),\,\mathcal L_\varepsilon\xi(\tau)) \: \d \tau \geq (C_F +1)\int_0^t \|\xi(\tau)\|^2_H \: \d\tau - C\int_0^t \|\xi(\tau)\|^2_{V^*} \:\d\tau,
    \end{equation}
    and finally by Assumption \ref{hyp:potential}
    \begin{equation} \label{eq:unique54}
        \int_0^t (\xi(\tau),\,F'(\varphi_\varepsilon(\tau))-F'(\varphi(\tau)))_H\: \d\tau \geq - C_F\int_0^t \|\xi(\tau)\|^2_H \: \d \tau.
    \end{equation}
    The stochastic diffusion can be handled as usual, i.e., recalling  again \cite[Proposition 2.1]{DPGS2}, we have
    \begin{equation} \label{eq:unique55}
        \begin{split}
            \int_0^t \left\|\nabla \mathcal N \left[ G(\varphi_\eps(\tau))-G(\varphi(\tau))  \right]\right\|^2_{\cL^2(U, \bH)}  \: \d \tau & \leq C\int_0^t \|G(\varphi_\eps(\tau)) - G(\varphi(\tau))\|^2_{V^*} \: \d \tau \\
            & \leq CL_G^2\int_0^t \|\xi(\tau)\|^2_{V^*} \: \d\tau
        \end{split}
    \end{equation}
    for some constant $C > 0$ only depending on $|\OO|$. Collecting \eqref{eq:unique51}-\eqref{eq:unique55} in \eqref{eq:unique50} we arrive at
    \begin{multline}
        \label{eq:unique57}
        \|\nabla \mathcal{N} \xi(t)\|^2_\b H + \int_0^t \|\xi(\tau)\|_H^2\: \d\tau  \leq C\left[ \|\nabla \mathcal{N} \xi_0\|^2_\b H + \int_0^t \|\xi(\tau)\|^2_{V^*} \: \d\tau \right. \\ \left.
        + \int_0^t \left( \hphi(\tau), \, \mathcal N \left[ G(\hphi_1(\tau))-G(\hphi_2(\tau))\right] \: \d \hW \right)_H + \eps\int_0^t \|\varphi(\tau)\|^2_{H^3(\OO)} \: \d \tau\right].
    \end{multline}
    Raising \eqref{eq:unique57} to the power $\frac p8$, taking supremums in time and $\P$-expectations leads to
    \begin{multline*}
        \E \supt \|\xi(\tau)\|^\frac p4_{V^*} + \E \left| \int_0^t \|\xi(\tau)\|^2_H \: \d \tau \right|^\frac p8 \leq C\left[ \E \|\xi_0\|^\frac p4_{V^*} + \E \left| \int_0^t \|\xi(\tau)\|^2_{V^*} \: \d \tau \right|^\frac p8 \right. \\
        \left. + \E \supt \left| \int_0^\tau \left( \xi(s), \, \mathcal N \left[ G(\varphi_\eps(s))-G(\varphi(s))  \right] \: \d W(s) \right)_H \right|^\frac p8 + \varepsilon^\frac p8\E\left|\int_0^t \|\varphi(\tau)\|^2_{H^3(\OO)} \: \d \tau\right|^\frac p8\right],
    \end{multline*}
    where the constant $C > 0$ depends on $p$. The stochastic integral is controlled using the Burkholder--Davis--Gundy inequality, namely
    \begin{equation*} 
        \begin{split}
        & \E \supt \left| \int_0^\tau \left( \xi(s), \, \mathcal N \left[ G(\varphi_\eps(s))-G(\varphi(s))  \right] \: \d W(s) \right)_H \right|^\frac p8  \\
        &\hspace{0.3cm} \leq C\E \left| \int_0^t \|\nabla \mathcal N\xi(\tau)\|_\bH^2\left\| \nabla \mathcal N \left[ G(\varphi_\eps(s))-G(\varphi(s))  \right] \right\|^2_{\cL^2(U, \b H)} \: \d\tau\right|^\frac p{16} \\
        &\hspace{0.3cm} \leq C\E \left[ \supt \|\nabla \mathcal N \xi(\tau)\|_\bH^\frac p8  \left| \int_0^t \| \xi(\tau)\|^2_{V^*} \: \d\tau\right|^\frac p{16} \right] \\
        & \hspace{0.3cm} \leq \dfrac{1}{4}\E \supt \|\xi(\tau)\|_{V^*}^\review{\frac p4} + C\E \left| \int_0^t \|\xi(\tau)\|^2_{V^*} \: \d \tau \right|^\frac p8.
        \end{split}
    \end{equation*}
    All of the estimates above then imply
    \begin{multline} \label{eq:unique5}
        \E \supt \|\xi(\tau)\|^\frac p4_{V^*} + \E \left| \int_0^t \|\xi(\tau)\|^2_H \: \d \tau \right|^\frac p8 \\ \leq C\left[ \E \|\xi_0\|^\frac p4_{V^*} + \E \left| \int_0^t \|\xi(\tau)\|^2_{V^*} \: \d \tau \right|^\frac p8 +  \varepsilon^\frac p8\E\left|\int_0^t \|\varphi(\tau)\|^2_{H^3(\OO)} \: \d \tau\right|^\frac p8\right],
    \end{multline}
    for a constant $C > 0$ only depending on the structural parameters of the problem, including $p$ and $T$. Therefore, \eqref{eq:unique5} and the Gronwall lemma imply that
    \begin{equation*}
        \|\hphi\|^\frac p4_{L^\frac p4(\Omega; L^\infty(0,T;V^*))} + \|\hphi\|^\frac p4_{L^\frac p4(\Omega; L^2(0,T;H))} \leq C\left[ \|\xi_0\|^\frac p4_{L^\frac p4(\Omega;V^*)} +  \varepsilon^\frac p8\E\left|\int_0^t \|\varphi(\tau)\|^2_{H^3(\OO)} \: \d \tau\right|^\frac p8\right]
    \end{equation*}
    and using the assumption on the initial condition the claim is proven.
    \section{\review{Nonlocal-to-local asymptotics: general regular potentials}} \label{sec:proof3}
    Finally, this section is devoted to show Theorem \ref{thm:convergence}. To establish the nonlocal-to-local asymptotics for the system in the case of general regular potentials, we devise uniform estimates and perform again a compactness argument. As we shall see, some additional effort is needed to obtain $\varepsilon$-uniform estimates.
    \subsection{Uniform estimates with respect to $\varepsilon$}
    In this subsection, we prove that the solutions to \eqref{eq:NCH:eps} satisfy uniform estimates independent of $\varepsilon$. The dependence of the estimate for the $H$-norm of $\varphi_\varepsilon$ on the norm of $\nabla K$ makes it challenging to retrieve a uniform estimate with respect to $\varepsilon$. Therefore, we resort to a weaker norm, in the same spirit of the computations shown in Section \ref{sec:convergence_unique}. More precisely, \review{we show the following lemma.}
    \begin{lem} \label{lem:first_estimate_eps}
        \review{The estimate
        \begin{multline*} 
        \E \supt \|\varphi_\eps(\tau)\|^p_{V^*} + \E\left| \int_0^t \Lambda_\eps(\varphi_\eps(\tau),\varphi_\eps(\tau)) \: \d \tau \right|^\frac p2 \\
        \leq C\left[ \E \|\varphi_{0,\eps}\|^p_{V^*} + \E \left| \int_0^t \|\varphi_\eps(\tau)\|^2_{V^*} \: \d \tau \right|^\frac p2 + \E\left|\int_0^t \|F(\varphi_\eps(\tau))\|_{L^1(\OO)}\: \d \tau \right|^\frac p2\right]
    \end{multline*}
    holds for all $p \geq 2$, all $\eps > 0$ and for a constant $C >0$ independent of $\eps$.}
    \end{lem}
    \begin{proof}
    We apply the It\^{o} formula to the functional
    \[
    \Gamma: V^*_0 \to \R, \qquad \Gamma(u) = \dfrac 12 \|\nabla \mathcal{N} u \|^2_\b H
    \]
    evaluated at $u=\varphi_\eps(t)-\overline{\varphi_{0,\varepsilon}}$. This is well defined owing to Assumption \ref{hyp:additional}, and results in 
    \begin{multline}
        \label{eq:conv1}
        \dfrac 12 \|\nabla \mathcal{N} (\varphi_\eps(t)-\overline{\varphi_{0,\varepsilon}})\|^2_\b H + \int_0^t (\varphi_\eps(\tau)-\overline{\varphi_{0,\varepsilon}},\,\mu_\eps(\tau))_H\: \d\tau \\ 
        = \dfrac 12 \|\nabla \mathcal{N} (\varphi_{0,\eps}-\overline{\varphi_{0,\varepsilon}})\|^2_\b H  + \dfrac 12 \int_0^t \left\|\nabla \mathcal N \left[ G(\varphi_\eps(\tau)) \right]\right\|^2_{\cL^2(U, \b H)}  \: \d \tau
         \\+ \int_0^t \left( \varphi_\eps(\tau)-\overline{\varphi_{0,\varepsilon}}, \, \mathcal N \left[ G(\varphi_\eps(\tau))\right] \: \d W \right)_H.
    \end{multline}
    Recalling the definition of $\mu_\eps$, it follows that
    \begin{align*}
        &\int_0^t (\varphi_\eps(\tau)-\overline{\varphi_{0,\varepsilon}},\,\mu_\eps(\tau))_H\: \d\tau \\
        & \hspace{1cm}= \int_0^t \big(\varphi_\eps(\tau)-\overline{\varphi_{0,\varepsilon}},\,\mathcal{L}_\eps\varphi_\eps(\tau) + F^\prime(\varphi_\eps(\tau))\big)_H\: \d\tau \\
        &\hspace{1cm}= 2\int_0^t\Lambda_\eps(\varphi_\eps(\tau),\, \varphi_\eps(\tau)) \: \d\tau + \int_0^t \big(\varphi_\eps(\tau), F^\prime(\varphi_\eps(\tau))\big)_H\: \d\tau - \int_0^t \overline{\varphi_{0,\varepsilon}} \, \overline{F'(\varphi_\eps(\tau))}|\OO| \: \d \tau,
    \end{align*}
    owing to the fact that $\overline{\mathcal L_\eps \varphi_\eps} = 0$. Here, we also defined the positive nonlocal energy contribution
    \begin{align*}
        \Lambda_\eps: H\times H\rightarrow \R^+, \quad
        \Lambda_\eps(h,k) := \dfrac 12\int_\OO\int_\OO K_\eps(x-y)(h(x)-h(y))(k(x)-k(y))\:\d y\: \d x.
    \end{align*}
    In particular, Assumption \ref{hyp:potential} and Remark \ref{rem:assumptions} yield
    \begin{align*}
        \int_0^t \big(\varphi_\eps(\tau), F^\prime(\varphi_\eps(\tau))\big)_H\: \d\tau 
        \geq -\alpha\int_0^t \|\varphi_\eps(\tau)\|_H^2\: \d\tau,
    \end{align*}
    while, proceeding similarly as in \eqref{eq:unique23}, Assumption \ref{hyp:additional} gives 
    \begin{align*}
        \int_0^t \left\|\nabla \mathcal N \left[ G(\varphi_\eps(\tau)) \right]\right\|^2_{\cL^2(U, \b H)}  \: \d \tau \leq C\int_0^t\|\varphi_\eps(\tau)\|_{V^*}^2\d\tau.
    \end{align*}
    Finally, observe that, for $\varepsilon$ sufficiently small, 
    \[
    \begin{split}
      \int_0^t |\overline{\varphi_{0,\varepsilon}}| \, |\overline{F'(\varphi_\eps(\tau))}||\OO|\:\christoph{\d\tau} & \leq C(1 + |\overline{\varphi_0}|)\left(t+\int_0^t \|F(\varphi_\eps(\tau))\|_{L^1(\OO)}\: \d \tau\right) \\
      & = C\left(t+\int_0^t \|F(\varphi_\eps(\tau))\|_{L^1(\OO)} \: \d \tau\right).
    \end{split}
    \]
    In light of all of the above, recalling the embedding $H \embed V^*$, \eqref{eq:conv1} reads as
    \begin{multline}
        \label{eq:conv2}
        \|\nabla \mathcal{N} (\varphi_\eps(t)-\overline{\varphi_{0,\varepsilon}})\|^2_\b H + \int_0^t\Lambda_\eps(\varphi_\eps(\tau), \varphi_\eps(\tau)) \: \d\tau  \\ 
        \leq C\left[ 1+ \|\nabla \mathcal{N} (\varphi_{0,\eps}-\overline{\varphi_{0,\varepsilon}})\|^2_\b H+ \int_0^t \|\varphi_\eps(\tau)\|_H^2\: \d\tau 
         + \int_0^t \|F(\varphi_\eps(\tau))\|_{L^1(\OO)}\: \d \tau \right. \\ \left. + \int_0^t \left( \varphi_\eps(\tau), \, \mathcal N \left[ G(\varphi_\eps(\tau))\right] \: \d W \right)_H \right].
    \end{multline}
    One crucial result enabled by the peculiar structure of the nonlocal kernels $K_\varepsilon$ is given by the so-called nonlocal Ehrling inequality (see \cite[Lemma 3.4]{Davoli1}), enabling to control the second term to the right hand side in terms of $\Lambda_\varepsilon$ and the $V^*$-norm. This then leads from \eqref{eq:conv2} to
    \begin{multline}
        \label{eq:conv3}
        \|\nabla \mathcal{N} (\varphi_\eps(t)-\overline{\varphi_{0,\varepsilon}})\|^2_\b H + \int_0^t\Lambda_\eps(\varphi_\eps(\tau), \varphi_\eps(\tau)) \: \d\tau  \\ 
        \leq C\left[ 1 + \|\nabla \mathcal{N} (\varphi_{0,\eps}-\overline{\varphi_{0,\varepsilon}})\|^2_\b H+ \int_0^t \|\varphi_\eps(\tau)\|_{V^*}^2\: \d\tau+\int_0^t \|F(\varphi_\eps(\tau))\|_{L^1(\OO)}\: \d \tau 
         \right. \\ \left. + \int_0^t \left( \varphi_\eps(\tau), \, \mathcal N \left[ G(\varphi_\eps(\tau))\right] \: \d W \right)_H \right].
    \end{multline}
    Summing to both sides of \eqref{eq:conv3} the average $|\overline{\varphi_{0,\varepsilon}}|^2$, then raising the result to the power $\tfrac{p}{2}$, taking supremums in time and $\P$-expectations, yields
    \begin{multline*}
        \E \supt \|\varphi_\eps(\tau)\|^p_{V^*} + \E\left| \int_0^t \Lambda_\eps(\varphi_\eps(\tau),\varphi_\eps(\tau)) \: \d \tau \right|^\frac p2 \\
        \leq C\left[ \E \|\varphi_{0,\eps}\|^p_{V^*} + \E \left| \int_0^t \|\varphi_\eps(\tau)\|^2_{V^*} \: \d \tau \right|^\frac p2 + \E\left|\int_0^t \|F(\varphi_\eps(\tau))\|_{L^1(\OO)}\: \d \tau \right|^\frac p2\right. \\
        \left. + \E \supt \left| \int_0^\tau \left( \varphi_\eps(s), \, \mathcal N \left[ G(\varphi_\eps(s))\right] \: \d W(s) \right)_H \right|^\frac p2 \right],
    \end{multline*}
    where the constant $C>0$ depends on $p$ but is independent of $\eps$. Exploiting the Burkholder-Davis-Gundy inequality as in the previous section, it holds
    \begin{align*}
        \E \supt \left| \int_0^\tau \left( \varphi_\eps(s), \, \mathcal N \left[ G(\varphi_\eps(s))\right] \: \d W(s) \right)_H \right|^\frac p2 \leq  \dfrac{1}{2}\E \supt \|\varphi_\eps(\tau)\|_{V^*}^p + C\E \left| \int_0^t \|\varphi_\eps(\tau)\|^2_{V^*} \: \d \tau \right|^\frac p2.
    \end{align*}
    Altogether, we conclude
    \begin{multline*} 
        \E \supt \|\varphi_\eps(\tau)\|^p_{V^*} + \E\left| \int_0^t \Lambda_\eps(\varphi_\eps(\tau),\varphi_\eps(\tau)) \: \d \tau \right|^\frac p2 \\
        \leq C\left[ \E \|\varphi_{0,\eps}\|^p_{V^*} + \E \left| \int_0^t \|\varphi_\eps(\tau)\|^2_{V^*} \: \d \tau \right|^\frac p2 + \E\left|\int_0^t \|F(\varphi_\eps(\tau))\|_{L^1(\OO)}\: \d \tau \right|^\frac p2\right]
    \end{multline*}
    where the constant $C>0$ is independent of $\eps$.
    \end{proof} \noindent
    In order to be able to close the estimate in Lemma \ref{lem:first_estimate_eps}, we need to refine the energy estimate. \review{Precisely, we show the following result.
    \begin{lem} \label{lem:second_estimate_eps}
        It holds that
        \begin{enumerate}[(i)]
            \item the family $\{\varphi_\eps\}_{\eps > 0}$ is uniformly bounded in $L^p(\Omega;L^\infty(0,T;V^*))$;
            \item the family $\{\Lambda_\eps(\varphi_\eps,\varphi_\eps)\}_{\eps > 0}$ is uniformly bounded in $L^\frac p2(\Omega;L^\infty(0,T))$;
            \item the family $\{\mu_\eps\}_{\eps>0}$ is uniformly bounded in $L^\frac p2(\Omega;L^2(0,T;V))$;
            \item the family $\{\nabla \mu_\eps\}_{\eps>0}$ is uniformly bounded in $L^p(\Omega;L^2(0,T;\b H))$;
            \item the family $\{\Psi(\varphi_\eps)\}_{\eps>0}$ is uniformly bounded in $L^\frac p2(\Omega;L^\infty(0,T;L^1(\OO)))$.
        \end{enumerate}
    \end{lem}
    }
    \begin{proof}
    We apply the It\^{o} formula to the full nonlocal energy functional
    \begin{align*}
        \mathcal E_\eps: H\rightarrow\R,\qquad \mathcal E_\eps(u) := \frac{1}{4}\int_\OO\int_\OO K_\eps(x-y)|u(x) - u(y)|^2 \: \d y \: \d x + \int_\OO F(u(x)) \: \d x
    \end{align*}
    evaluated at $u=\varphi_\eps(t)$.
    This yields
    \begin{multline} \label{eq:conv4}
        \mathcal E_\eps(\varphi_\eps(t)) + \int_0^t\int_\OO|\nabla\mu_{\eps}(\tau)|^2\:\d x\:\d \tau
        = \mathcal E_\eps(\varphi_{0, \eps}) + \int_0^t\Big(\mu_{\eps}(s),G(\varphi_\eps(s))\:\d W(s)\Big)_H \\
        + \frac{1}{2}\int_0^t\sum_{k=0}^{+\infty}\left[\frac{1}{2}\Lambda_\eps(G(\varphi_\eps(\tau))[u_k], G(\varphi_\eps(\tau))[u_k])
        +\int_\OO F^{\prime\prime}(\varphi_\eps(\tau))|G(\varphi_\eps(\tau))[u_k]|^2\:\d x
        \right]\:\d \tau.
    \end{multline}
    By means of the nonlocal Poincaré (see \cite[Theorem 1.1]{Ponce}) and Ehrling inequalities, as well as Assumption \ref{hyp:G2}, we can control the third term on the right hand side as
    \begin{align*}       \int_0^t\sum_{k=0}^{+\infty}\Lambda_\eps(G(\varphi_\eps(\tau))[u_k], G(\varphi_\eps(\tau))[u_k])\:\d \tau &\leq C\int_0^t\sum_{k=0}^{+\infty}\|G(\varphi_\eps(\tau))[u_k]\|_V^2\:\d \tau \\
        &=C\int_0^t\|G(\varphi_\eps(\tau))\|_{\cL^2(U,V)}^2\:\d \tau \\
        &\leq C\left( 1+ \int_0^t\|\varphi_\eps(\tau)\|_H^2\:\d \tau\right)\\
        &\leq C\left( 1+\int_0^t \Lambda_\eps(\varphi_\eps(\tau), \varphi_\eps(\tau)) \: \d \tau + \int_0^t\|\varphi_\eps(\tau)\|_{V^*}^2\:\d \tau\right).
    \end{align*}
    As, by assumption, the initial energies $\mathcal E_\eps(\varphi_{0,\eps})$ are suitably uniformly bounded, taking into account the above inequality in \eqref{eq:conv4}, raising to the power $\tfrac{p}{2}$ and taking the supremum in time as well as $\P$-expectations, we have
    \begin{multline*}
    \E \supt |\Lambda_\varepsilon(\varphi_\eps(t), \varphi_\eps(\tau))|^\frac p2 + \E\supt\|F(\varphi_\eps(t))\|_{L^1(\OO)}^{\tfrac{p}{2}} + \E\left|\int_0^t\|\nabla\mu_\eps(\tau)\|_\b H^2\:\d \tau\right|^{\tfrac{p}{2}} \\
        \leq C\Bigg(1+\E \left| \int_0^t \Lambda_\eps(\varphi_\eps(\tau), \varphi_\eps(\tau)) \: \d \tau\right|^\frac p2 +\E\left|\int_0^t\|\varphi_\eps(\tau)\|_{V^*}^2\:\d \tau\right|^\frac p2 \\  +\E\left|\int_0^t\Big(\mu_{\eps}(s),G(\varphi_\eps(s))\:\d W(s)\Big)_H\right|^{\tfrac{p}{2}} 
        + \E\left|\int_0^t\sum_{k=0}^{+\infty}\int_\OO F^{\prime\prime}(\varphi_\eps(\tau))|G(\varphi_\eps(\tau)[u_k]|^2\:\d x
        \:\d \tau\right|^{\tfrac{p}{2}}\Bigg).
    \end{multline*} 
    Using the Burkholder-Davis-Gundy inequality (as in \eqref{eq:burkholdermu}), we can bound the stochastic integral on the right-hand side and get an additional term involving the mean of $\mu_\eps$. For the remaining term, we invoke Assumptions \ref{hyp:potential}-\ref{hyp:growth}, \ref{hyp:G}-\ref{hyp:Gii} and Hölder's inequality to get 
    \[ 
    \begin{split}
        & \E\left|\int_0^t\sum_{k=0}^{+\infty}\int_\OO |F^{\prime\prime}(\varphi_{\eps}(\tau))||G(\varphi_{\eps}(\tau)[u_k]|^2\:\d x \: \d s\right|^{\frac{p}{2}} \\
        & \hspace{1cm}\leq C\left[ 1  + \E\left|\int_0^t\sum_{k=0}^{+\infty}\int_\OO |F(\varphi_{\eps}(\tau))||G(\varphi_{\eps}(\tau)[u_k]|^2\:\d x \: \d s\right|^{\frac{p}{2}} \right] \\
        & \hspace{1cm}\leq C\left[ 1  + \E\left|\int_0^t\|F(\varphi_{\eps}(\tau))\|_{L^1(\OO)}\sum_{k=0}^{+\infty} \|G(\varphi_{\eps}(\tau)[u_k]\|_{L^\infty(\OO)}^2 \: \d s\right|^{\frac{p}{2}} \right] \\
        & \hspace{1cm}\leq C\left[ 1  + \E\left|\int_0^t\|F(\varphi_{\eps}(\tau))\|_{L^1(\OO)} \: \d s\right|^{\frac{p}{2}} \right].
    \end{split}
    \]
    Finally, it is left to control the mean value of $\mu_\eps$. Testing the equation for the chemical potential by $1$ and recalling that $\overline{\mathcal{L}_\eps\varphi_\eps}=0$, we arrive at
    \[
    \begin{split}
        \left| |\OO|\:\overline{\mu_{\eps}} \right|  = \left| (\mu_{\eps}, 1)_H \right| 
        & \leq |(\mathcal{L}_\eps\varphi_{\eps}, 1)_H| + |(F^\prime(\varphi_{\eps}), 1)_H| \\
        & \leq C\left( 1 + \| F(\varphi_{\eps})\|_{L^1(\OO)}\right),
    \end{split}
    \]
    implying
    \begin{equation} \label{eq:averagemu}
        \E \supt  |\overline{\mu_\varepsilon(\tau)}|^\frac p2 \le C\left( 1 + \E \supt \|F(\varphi_\eps(\tau))\|_{L^1(\OO)}^\frac p2 \right).
    \end{equation}
    Eventually, collecting all the terms above, we arrive at
    \begin{multline} \label{eq:convergence2}
    \E \supt |\Lambda_\varepsilon(\varphi_\eps(\tau), \varphi_\eps(\tau))|^\frac p2+\E\supt\|F(\varphi_\eps(\tau))\|_{L^1(\OO)}^{\frac{p}{2}} + \E\left|\int_0^t\|\nabla\mu_\eps(\tau)\|_H^2\:\d \tau\right|^{\frac{p}{2}} \\
        \leq C\Bigg(1+\E \left| \int_0^t \Lambda_\eps(\varphi_\eps(\tau), \varphi_\eps(\tau)) \: \d \tau\right|^\frac p2+ \E\left|\int_0^t\|\varphi_\eps(\tau)\|_{V^*}^2\:\d \tau\right|^\frac p2  +\E \left|\int_0^t \|F(\varphi_\eps(\tau))\|_{L^1(\OO)} \: \d \tau \right|^\frac p2\Bigg).
    \end{multline}
    Summing \review{the estimate in Lemma \ref{lem:first_estimate_eps}},  \eqref{eq:convergence2} as well as \eqref{eq:averagemu} multiplied by a sufficiently small constant yields
    \begin{multline}
        \E \supt \|\varphi_\eps(\tau)\|^p_{V^*} + \E \supt |\Lambda_\varepsilon(\varphi_\eps(\tau), \varphi_\eps(\tau))|^\frac p2+\E\supt\|F(\varphi_\eps(\tau))\|_{L^1(\OO)}^{\frac{p}{2}} \\+ \E \supt  |\overline{\mu_\varepsilon(\tau)}|^\frac p2+\E\left|\int_0^t\|\nabla\mu_\eps(\tau)\|_H^2\:\d \tau\right|^{\frac{p}{2}}  \\
        \leq C\left[ 1 +\E \|\varphi_{0,\eps}\|^p_{V^*} + \E \left| \int_0^t \|\varphi_\eps(\tau)\|^2_{V^*} \: \d \tau \right|^\frac p2 + \E\left|\int_0^t \|F(\varphi_\eps(\tau))\|_{L^1(\OO)}\: \d \tau \right|^\frac p2\right.\\
        \left. + \E \left| \int_0^t \Lambda_\eps(\varphi_\eps(\tau), \varphi_\eps(\tau)) \: \d \tau\right|^\frac p2\right]
    \end{multline}
    and the Gronwall lemma gives the estimates
    \begin{align*}
    \|\varphi_\eps\|_{L^p(\Omega;L^\infty(0,T;V^*))}+\|\Lambda_\eps(\varphi_\eps, \varphi_\eps)\|_{L^\frac p2(\Omega;L^\infty(0,T))} & \leq C,\\\label{eq:conv52}
        \| \mu_{\eps}\|_{L^\frac p2(\Omega;L^2(0,T;V))} + \|\nabla \mu_{\eps}\|_{L^p(\Omega;L^2(0,T;\b H))}+\|\Psi(\varphi_\eps)\|_{L^\frac p2(\Omega;L^\infty(0,T;L^1(\OO)))} &\leq C,
    \end{align*}
    and the proof is complete.
    \end{proof} \noindent
    \review{For brevity, we collect the next three estimates in the following result.
    \begin{lem} \label{lem:third_fourth_estimate_eps}
        It holds that
        \begin{enumerate}[(i)] \itemsep0.3em
            \item the family $\{\varphi_\eps\}_{\eps >0}$ is uniformly bounded in $L^p(\Omega;L^\infty(0,T;H))$ for all $p \geq 2$;
            \item the family $\{G(\varphi_\eps)\}_{\eps >0}$ is uniformly bounded in 
            \[
            L^\infty(\Omega\times(0,T);\cL^2(U,H))\cap L^p(\Omega;L^\infty(0,T;\cL^2(U,V)))
            \]
            for all $p \geq 2$;
            \item the family $\left\{\displaystyle\int_0^. G(\varphi_{\eps}(s)) \: \d W(s)\right\}_{\eps>0}$ is uniformly bounded in 
            \[
            L^q(\Omega;W^{k,q}(0,T;H))\cap L^p(\Omega;W^{k,p}(0,T;V))
            \]
            for all $p \geq 2$, $q \geq 2$ and $k \in (0, \frac 12)$.
        \end{enumerate}
    \end{lem}
    }
    \begin{proof}
    Once again, the nonlocal Ehrling inequality proves to be a crucial result to carry out the nonlocal-to-local asymptotics. Indeed, for $\varepsilon$ sufficiently small, it holds
    \begin{align*}
        \|\varphi_\eps(t)\|_H^2 \leq \delta\Lambda_\eps(\varphi_\eps(t),\varphi_\eps(t)) + C_\delta\|\varphi_\eps(t)\|_{V^*}^2
    \end{align*}
    for some fixed $\delta > 0$, and therefore
    \[
    \|\varphi_\eps\|_{L^p(\Omega;L^\infty(0,T;H))} \leq C\left(\|\Lambda_\eps(\varphi_\eps,\varphi_\eps)\|_{L^\frac p2(\Omega;L^\infty(0,T))} + \|\varphi_\eps\|_{L^p(\Omega;L^\infty(0,T;V^*))} \right).
    \]
    Owing to Lemma \ref{lem:second_estimate_eps}, we then conclude that 
    \begin{align} \label{eq:conv6}
        \|\varphi_\eps\|_{L^p(\Omega;L^\infty(0,T;H))}\leq C
    \end{align}
    for a constant $C > 0$ independent of $\eps$.
    Arguing as in Subsection \ref{ssec:unif_n}, exploiting \cite[Lemma 2.1]{fland-gat} jointly with Assumptions \ref{hyp:G} and \ref{hyp:G2}, we get
    \begin{align*}
        \|G(\varphi_\eps)\|_{L^\infty(\Omega\times(0,T);\cL^2(U,H))\cap L^p(\Omega;L^\infty(0,T;\cL^2(U,V)))} &\leq C, \\
        \left\|\int_0^. G(\varphi_{\eps}(s)) \: \d W(s)\right\|_{L^q(\Omega;W^{k,q}(0,T;H))\cap L^p(\Omega;W^{k,p}(0,T;V))} &\leq C
    \end{align*}
    for every $k\in(0,\frac{1}{2})$ and $q\geq 2$. Similarly as in \review{the proof of Lemma \ref{lem:fourth_estimate}}, it follows that
    \begin{align*}     
    \|\varphi_{\eps}\|_{L^p(\Omega;W^{\alpha,p}(0,T;V^*))} \leq C
    \end{align*}
    for some $\alpha = \alpha(p) > \frac 1p$ if $p > 2$.
    \end{proof} \noindent
    \review{Finally, we show the last estimate.
    \begin{lem} \label{lem:fifth_estimate_eps}
        The families $\{F^\prime(\varphi_\eps)\}_{\eps>0}$ and $\{\mathcal L_\eps \varphi_\eps\}_{\eps>0}$ are uniformly bounded in $L^{\frac{p}{2}}(\Omega;L^2(0,T;H))$.
    \end{lem}}
    \begin{proof}
    Using Assumption \ref{hyp:potential} and recalling Remark \ref{rem:assumptions}, we can multiply the equation for $\mu_\eps$ by $F^\prime(\varphi_\eps)+\gamma\varphi_\eps$, where $\gamma > \alpha$. Hence, we get
    \begin{align*}
        \int_\OO\mu_\eps F^\prime(\varphi_\eps)\;\d x + \gamma\int_\OO \varphi_\eps\mu_\eps\;\d x = \int_\OO \mathcal{L}_\eps\varphi_\eps (F^\prime(\varphi_\eps)+\gamma\varphi_\eps)\;\d x +\int_\OO|F^\prime(\varphi_\eps)|^2\;\d x + \int_\OO \varphi_\eps F'(\varphi_\eps) \: \d x.
    \end{align*}
    Letting $\Psi'(s) := F^\prime(s)+\gamma s$ for any $s \in \mathbb R$, we observe that it is monotone by construction. Therefore, since the kernel $K_\eps$ is nonnegative, we have
    \begin{align*}
        \int_\OO \mathcal{L}_\eps\varphi_\eps\Psi^\prime(\varphi_\eps)\;\d x = \frac{1}{2}\int_\OO\int_\OO K_\eps(x-y)(\varphi_\eps(x)-\varphi_\eps(y))(\Psi^\prime(\varphi_\eps(x))-\Psi^\prime(\varphi_\eps(y)))\;\d y\;\d x \geq 0.
    \end{align*}
    Eventually, invoking Young`s inequality twice, it follows that 
    \begin{align*}
        \|F^\prime(\varphi_\eps)\|_{H}^2 \leq C(\|\mu_\eps\|_H^2 + \|\varphi_\eps\|_H^2)
    \end{align*}
    and therefore, owing to the proven estimates,
    \begin{align*}
        \|F^\prime(\varphi_\eps)\|_{L^\frac{p}{2}(\Omega;L^2(0,T;H))} \leq C,
    \end{align*}
    where the constant $C$ does not depend on $\eps$. As a byproduct, by comparison in the equation for $\mu_\eps$, we obtain
    \begin{align*}
        \|\mathcal{L}_\eps\varphi_\eps\|_{L^\frac{p}{2}(\Omega;L^2(0,T;H))} \leq C,
    \end{align*}
    and the proof is complete.
    \end{proof}
    \subsection{Passage to the limit as $\eps\rightarrow0^+$.}
    In this section, we pass to the limit as $\eps\rightarrow0^+$ in the nonlocal problems \eqref{eq:NCH:eps}, exploiting the previously shown uniform estimates. Accordingly, we tacitly fix a decreasing vanishing sequence $\{\eps_n\}_{n \in \mathbb N}$ without relabeling. In particular, following the same logic of Subsections \ref{ssec:lim_n} and \ref{ssec:lim_lam}, we investigate the tightness properties of the sequence of approximate solutions.
    Once again, we recall that the embedding
    \begin{align*}
        L^\infty(0,T;H) \cap W^{\alpha, p}(0,T;V^*) &\hookrightarrow C^0([0,T]; H^{-\delta}(\OO))
    \end{align*}
    is compact whenever $\alpha p > 1$ for some small $\delta > 0$. To recover compactness in $L^2(0,T;H)$, a different route must be followed, as a uniform estimate in $L^p(\Omega;L^2(0,T;V))$ is typically out of reach, due to the singular behavior of $\nabla K_\eps$ in the limit $\eps\to0^+$. Here, we shall proceed a posteriori. The following result holds: its proof is entirely analogous to the one of Lemma \ref{lem:tightness} and is therefore omitted.
    \begin{lem} \label{lem:tightnonlocal}
        If $p > 2$, then the family of laws of $\{\varphi_{\eps}\}_{\eps>0}$ is tight in the space $C^0([0,T]; H^{-\delta}(\OO))$ for any $\delta \in (0, \frac 12)$.
    \end{lem}
    \noindent
    Similarly, setting
    \begin{align*}
       G(\varphi_{\eps})\cdot W :=  \int_0^. G(\varphi_{\eps}(s))\:\d W(s),
    \end{align*}
    we immediately have the counterpart of Lemma \ref{lem:tightness2}, i.e., 
    \begin{lem} \label{lem:tightnonlocal2}
        The family of laws of $\{G(\varphi_\eps) \cdot W\}_{\eps>0}$ is tight in the space $C^0([0,T]; V^*)$.
    \end{lem} \noindent
    It even holds a stronger statement, due to Assumption \ref{hyp:G2}: the sequence is even tight in $C^0([0,T];H)$ if $p > 2$. Nonetheless, Lemma \ref{lem:tightnonlocal2} is sufficient for our purposes. As for the last two tightness lemmas, arguing similarly to Lemmas \ref{lem:tightness3} and \ref{lem:tightness4}, i.e., identifying the Wiener process $W$ with a constant sequence of random variables, we first have
    \begin{lem} \label{lem:tightnonlocal3}
        The family of laws of $\{W_\eps\}_{\eps>0}$ is tight in the space $C^0([0,T]; U_0)$.
    \end{lem} \noindent
    Then, for the initial conditions, given that $\varphi_{0,\eps} \to \varphi_0$ in $L^p(\Omega;H)$, we have
    \begin{lem} \label{lem:tightnonlocal4}
        The family of laws of $\{\varphi_{0,\eps}\}_{\eps>0}$ is tight in $V^*$.
    \end{lem} \noindent
    Therefore, by Lemmas \ref{lem:tightnonlocal}-\ref{lem:tightnonlocal4}, the family of laws of
    \begin{align*}
        \{\varphi_{\eps},\, G(\varphi_{\eps})\cdot W,\, W_\eps,\, \varphi_{0,\eps}\}_{\eps>0}
    \end{align*}
    is tight in the space
    \begin{align*}
        C^0([0,T];H^{-\delta}(\OO)) \times C^0([0,T]; V^*) \times C^0([0,T];U_0) \times V^*.
    \end{align*}
    In order to pass to the limit, we can exploit the approach established by Gy\"{o}ngy and Krylov in \cite [Lemma 1.1]{gyo-kry}, by showing convergence in probability. Fix two arbitrary subsequences of $\{\varphi_\eps\}_{\eps>0}$, that we denote with
        \[
        \{\varphi_{\eps_j}\}_{j \in \enne}, \qquad 
        \{\varphi_{\eps_k}\}_{k \in \enne},
        \]
        having values in the Polish space $\mathcal X := C^0([0,T];H^{-\delta}(\OO))$. The sequence of joint laws of $(\varphi_{\eps_j}, \varphi_{\eps_k})$ is clearly tight in $\mathcal X^2$ by Lemma \ref{lem:tightnonlocal}. Therefore, by the Prokhorov theorem, the sequence is relatively compact in the space of probability measures over $\mathcal X^2$ (and its Borel $\sigma$-algebra), i.e., it admits a weakly convergent subsequence, namely
        \[
        \operatorname{Law}_\P (\varphi_{\eps_{j_i}}, \varphi_{\eps_{k_i}}) \rightharpoonup \varsigma
        \]
        for suitable choices of $j_i$ and $k_i$ as $i \in \mathbb N$ and for some limit probability measure $\christoph{\varsigma}$. Using the Skorokhod theorem, as in the proof of existence of martingale solutions, we can frame this convergence in a new probability space $(\Omega', \cF', \P')$ through law-preserving maps $Z_i:(\Omega', \cF')\to(\Omega, \cF)$ (i.e., such that $\P' \circ Z_i^{-1} = \P$ for all $i \in \mathbb N$). In particular, we obtain
        \[
        (\varphi_{\eps_{j_i}}', \varphi_{\eps_{k_i}}') :=(\varphi_{\eps_{j_i}}, \varphi_{\eps_{k_i}}) \circ Z_i \to (\varphi_1,\varphi_2) \qquad \text{in }\mathcal X^2, \quad \P\text{-a.s.}
        \]
        for some limit process
        \[
        (\varphi_1,\varphi_2): (\Omega', \cF') \to (\mathcal X^2, \mathscr B(\mathcal X^2)).
        \]
        By the same token, the remaining tightness results yield the convergences
        \begin{align*}
            (I'_{\ji}, I'_{\ki}) := (G(\varphi_{\ji}) \cdot W,\, G(\varphi_{\ki}) \cdot W) \circ Z_i \to (I_1, I_2) & \qquad \text{in }\mathcal X^2, \quad \P\text{-a.s.} \\
            (W'_{\ji}, W'_{\ki}) := (W_{\ji}, W_{\ki}) \circ Z_i \to (W_1, W_2) & \qquad \text{in } [C^0([0,T];U_0]^2, \quad \P\text{-a.s.}\\
            (\varphi_{0, \ji}', \varphi_{0, \ki}') := (\varphi_{0, \ji}, \varphi_{0, \ki}) \circ Z_i \to (\varphi_{01}, \varphi_{02}) & \qquad \text{in }V^* \times V^*, \quad \P\text{-a.s.}\qquad 
        \end{align*}
        where the limit processes belong to the specified spaces. Since the sequence of Wiener processes is constant and given the known convergence of the initial conditions, as laws are preserved it is immediate to identify
        \[
        W_1 = W_2 = W, \qquad \varphi_{01}=\varphi_{02}= \varphi_0.
        \]
        On account of the uniform estimates proven in the previous subsection, by reflexivity and compactness we can establish the following convergences:
       \begin{align*}
        (\varphi_{\eps_{j_i}}', \varphi_{\eps_{k_i}}') \to (\varphi_1, \varphi_2) & \qquad \text{ in } L^q(\Omega', \mathcal X^2) \quad \forall \: q < p, \\
       (\varphi_{\eps_{j_i}}', \varphi_{\eps_{k_i}}') \overset{*}{\rightharpoonup}(\varphi_1, \varphi_2) &\qquad\text{ in }L^p_w(\pom; L^\infty(0,T;H)) \cap L^p(\pom;W^{\alpha, p}(0,T;V^*)),\\
        (\mu_{\eps_{j_i}}', \mu_{\eps_{k_i}}')  := (\mu_{\eps_{j_i}}, \mu_{\eps_{k_i}})\circ Z_i \rightharpoonup \
        {(\mu_1, \mu_2)}&\qquad\text{ in }L^{\frac p2}(\pom;L^2(0,T;V)). \\
        (F'(\varphi_{\ji}'), F'(\varphi_{\ki}')) \rightharpoonup {(\eta_1,\eta_2)} &\qquad\text{ in }L^{\frac p2}(\pom;L^2(0,T;H)). \\
        (\mathcal L_{\ji}(\varphi_{\ji}'), \mathcal L_{\ki}(\varphi_{\ki}')) \rightharpoonup {(\xi_1,\xi_2)} &\qquad\text{ in }L^{\frac p2}(\pom;L^2(0,T;H)).
        \end{align*}
        Fix then any two sufficiently large values $\ell$ and $m \in \mathbb N$. Exploiting the nonlocal version of Ehrling`s inequality, for any $\rho>0$ there exists $C_\rho>0$ such that
    \begin{align*}
        \E\:\hspace{-1.4mm}'\|\varphi_{\eps_{j_\ell}}' - \varphi_{\eps_{j_m}}'\|_{L^\infty(0,T;H)}^q 
        &\leq \rho\pE\|\Lambda_{\eps_{j_\ell}}(\varphi_{\eps_{j_\ell}}') + \Lambda_{{j_m}}(\varphi_{\eps_{j_m}}')\|_{L^\infty(0,T)}^\frac q2 + C_\rho\pE\|\varphi_{\eps_{j_\ell}}' - \varphi_{\eps_{j_m}}'\|_{C^0(0,T;V^*)}^q \\
        &\leq C\rho + C_\rho\pE\|\varphi_{\eps_{j_\ell}}' - \varphi_{\eps_{j_m}}'\|_{C^0(0,T;V^*)}^q,
    \end{align*}
    where $q < p$ and we exploited the \review{second uniform estimate in Lemma \ref{lem:second_estimate_eps}}. Given the strong convergence in $\mathcal X$ and that $\rho>0$ can be chosen arbitrarily small, this yields that $\{\varphi_{\ji}'\}_{i\in \mathbb N}$ is a Cauchy sequence in $L^q(\pom;L^\infty(0,T;H))$, and therefore
    \begin{align}\label{sconv:phi}
        (\varphi_{\ji}', \varphi_{\ki}') \rightarrow (\varphi_1, \varphi_2) & \quad \text{ in }L^q(\pom;[L^\infty(0,T;H)]^2) \quad  \forall \: q < p.
    \end{align}
    As a byproduct, recalling Assumption \ref{hyp:G} on $G$, we derive that
    \begin{align*}
        \|G(\varphi_{\ji}') - G(\varphi_1)\|_{L^q(\pom;L^\infty(0,T;\cL^2(U,H)))} \rightarrow 0
    \end{align*}
    as $i\to+\infty$, and therefore,
    \begin{align*}
        (G(\varphi_{\ji}'), G(\varphi_{\ki}')) \rightarrow (G(\varphi_1),G(\varphi_2)) \quad \text{ in }L^q(\Omega';[L^\infty(0,T;\cL^2(U,H))]^2) \quad \forall \: q <p.
    \end{align*}
    for all $q<p$. Next, setting as usual $\Psi'(s) = F'(s) + \gamma s$ for some $\gamma > \alpha$ and all $s \in \erre$, let us show that $\eta_m = \Psi^\prime(\varphi_m) -\gamma\varphi_m = F^\prime(\varphi_m)$ for $m \in \{1,2\}$. Using the strong convergence \eqref{sconv:phi}, by comparison we infer that
    \begin{align*}
        (\Psi^\prime(\varphi_{\ji}'), \Psi^\prime(\varphi_{\ki}')) \rightharpoonup (\eta_1 + \gamma\varphi_1, \eta_2+\gamma\varphi_2) \quad \text{ in }L^{\frac{p}{2}}(\Omega';[L^2(0,T;H)]^2).
    \end{align*}
    The identification is then a consequence of the weak-strong closure of maximal monotone operators.
    Finally, we show that $\varphi_m\in L^p(\pom;L^\infty(0,T;V))\cap L^{\frac{p}{2}}(\pom;L^2(0,T;H^2(\OO)))$ and $\xi_m = -\Delta\varphi_m$ for $m \in \{1,2\}$. For instance, since \eqref{sconv:phi} holds in $L^q(\Omega' \times(0,T) ; H)$ as well, it holds that
    \[
    \varphi_{\ji}'(t) \to \varphi_1(t) \quad \text{in }H
    \]
    in the $\P' \otimes \d t$-almost sure sense, up to a non-relabeled subsequence. We can exploit pathwise \cite[Lemma 3.3]{Davoli1} to get
    \begin{align*}
        \Lambda_0(\varphi_1(t)) \leq \liminf_{i\to+\infty} \Lambda_\eps(\varphi_{\ji}'(t), \varphi_{\ji}'(t)),
    \end{align*}
    holding $\P' \otimes \d t$-almost surely, where we set
    \begin{align*}
        \Lambda_0(u) := \frac{1}{2}\int_\OO|\nabla u(x)|^2\;\d x\qquad \forall \:u\in V.
    \end{align*}
    Eventually, taking powers $\frac p2$, supremums in time and $\P'$-expectations yields, owing also to Fatou's lemma,
    \[
    \pE\sup_{t \in [0,T]}|\Lambda_0(\varphi_1(t))|^\frac p2 \leq \liminf_{i\to+\infty} \pE\sup_{t \in [0,T]}|\Lambda_\eps(\varphi_{\ji}'(t), \varphi_{\ji}'(t)|^\frac p2 \leq C,
    \]
    ensuring that $\varphi_1\in L^p(\pom;L^\infty(0,T;V))$. Furthermore, due to \cite[Lemma 3.2]{Davoli1}, it holds
    \begin{multline*}
        \int_0^T\Lambda_{\ji}(\varphi_{ji}'(\tau), \varphi_{ji}'(\tau))\;\d\tau + \int_0^T\int_\OO ((K_{\ji}*1)\varphi_{\ji}' - K_\eps*\varphi_{\ji}')(\tau,x)(\zeta-\varphi_{\ji}')(\tau,x)\;\d x\:\d\tau \\ \leq \int_0^T\Lambda_{\ji}(\zeta(\tau),\zeta(\tau))\;\d\tau
    \end{multline*}
    for all $\zeta\in L^2(0,T;V)$. On account of \eqref{sconv:phi}, \cite[Lemma 3.3]{Davoli1} as well as Fatou`s lemma, it follows that 
    \begin{align*}
        \frac{1}{2}\int_0^T\int_\OO|\nabla\varphi_1(\tau,x)|^2\d x\:\d\tau + \int_0^T\int_\OO \xi_1(\tau,x)(\zeta-\varphi_1)(\tau,x)\;\d x\:\d\tau \leq \frac{1}{2}\int_0^T\int_\OO|\nabla\zeta(\tau,x)|^2\;\d x\:\d\tau
    \end{align*}
    for all $\zeta\in L^2(0,T;V)$, $\P'$-almost surely. Indeed, this qualifies $\xi_1$ as an element of the subdifferential $\partial \widetilde \Lambda_0$ at the point $u=\varphi_1$, where
    \[
    \widetilde{\Lambda}_0(u) = \dfrac 12\int_0^T\int_\OO |\nabla u(x,t)|^2 \: \d x \: \d t \qquad \forall \: u \in L^2(0,T;V),
    \]
    extended with value $+\infty$ to the whole $L^2(0,T;H)$. Therefore, we have
    \begin{align*}
        \int_0^T\int_\OO\xi_1(x, t)\zeta(x, t)\;\d x \: \d t = \int_0^T\int_\OO\nabla\varphi_1(x,t)\cdot\nabla\zeta(x,t)\;\d x\: \d t\qquad \forall \: \zeta\in L^2(0,T;V),
    \end{align*}
    and, in turn,
    \begin{align*}
        \int_\OO\xi_1(x, t)\zeta(x)\;\d x = \int_\OO\nabla\varphi_1(x, t)\cdot\nabla\zeta(x)\;\d x\qquad \forall \: \zeta\in V.
    \end{align*}
    Thus, by elliptic regularity theory, we conclude that $\varphi_1\in L^2(0,T;H^2(\OO)))$ $\P'$-almost surely and that $\xi_1(t) = -\Delta\varphi_1(t)$ for almost all times and $\P'$-almost surely, as well as $\partial_\bn\varphi_1=0$ almost everywhere on $\partial\OO \times (0,T)$, and $\P'$-almost surely. Furthermore, given
    \[
    \|\varphi_1\|_{H^2(\OO)}^2 \leq C\left(\|\varphi_1\|^2_H + \|\Delta \varphi_1\|_H^2 \right)
    \]
    we eventually deduce that $\varphi_1 \in L^\frac p2(\pom, L^2(0,T;H^2(\OO))$, matching the regularity provided in \cite{scarpa21}. The proven convergences are now enough to identify the limit of the problem. After carrying out a procedure analogous to the one illustrated in Subsection \ref{ssec:limit_n} to identify the stochastic integrals, it is straightforward to show that the limit processes $(\varphi_1, \varphi_2)$ satisfy the local Cahn--Hilliard system
    \[
    \begin{cases}
        \d \varphi_m -\Delta \mu_m \: \d t = G(\varphi_m) \: \d W & \quad \text{in } \OO \times (0,T), \\
        \mu_m = -\Delta\varphi_m + F'(\varphi_m)& \quad \text{in } \OO \times (0,T), \\
        \partial_\bn \mu_m = \partial_\bn\varphi_m = 0 & \quad \text{on } \partial\OO \times (0,T), \\
        \varphi_m(0) = \varphi_0 & \quad \text{in } \OO, 
    \end{cases}
    \]
    where $m \in \{1,2\}$. The pathwise uniqueness result given in \cite{scar-SCH} for the system, however, yields immediately that $\varphi_1 = \varphi_2$, $\P'$-almost surely. Thanks to preservation of laws, therefore, the limit law $\varsigma$ must be supported on the diagonal of $\mathcal X^2$, i.e., 
    \[
    \varsigma((v_1, v_2) \in \mathcal X^2: v_1 = v_2) = 1.
    \]
    Then, \cite[Lemma 1.1]{gyo-kry} yields convergence in probability with respect to the topology of $\mathcal X$ for the original sequence $\varphi_\varepsilon$, on the original probability space $(\Omega, \cF, \P)$. Iterating precisely the previous argument on the sequence $\{\varphi_\eps\}_{\eps > 0}$ finally completes the proof of Theorem \ref{thm:convergence}.
    \section*{Acknowledgments}
     The authors warmly thank the anonymous reviewer for their useful remarks, which helped to improve the quality of the manuscript significantly. The first author is a member of Gruppo Nazionale per l'Analisi Matematica, la Probabilità e le loro Applicazioni (GNAMPA), Istituto Nazionale di Alta Matematica  (INdAM). The research of the first author is funded by the European Union (ERC, NoisyFluid, No. 101053472). The second author was partially supported by the RTG 2339 ``Interfaces, Complex Structures, and Singular Limits'' of the Deutsche Forschungsgemeinschaft (DFG, German Research Foundation).
    \printbibliography
\end{document}
Altogether, we conclude that 
    \begin{align}
        \label{eq:conv:4}
        \|\varphi_\eps\|_{L^p(\Omega;L^\infty(0,T;V^*))}^p + \|\varphi_\eps\|_{L^p(\Omega;L^2(0,T;H))}^p \leq C,
    \end{align}
 Therefore, letting $n\rightarrow\infty$ in \eqref{eq:estimate_2} we arrive at
    \begin{align*}
        &\frac{1}{2}\widetilde{\E}\sups\mathcal{E}_{1,\lambda}(\widetilde{\varphi}_{\lambda}(s))^{\frac{p}{2}} + \widetilde{\E}\sups\|F_\lambda(\widetilde{\varphi}_{\lambda}(s))\|_{L^1(\OO)}^{\frac{p}{2}} + \widetilde{\E}\|\nabla\widetilde{\mu}_{\lambda}\|_{L^2(0,t;H)}^p + \widetilde{\E}\sups|(\widetilde{\mu}_{\lambda}(s))_\OO|^\frac{p}{2} \\
        &\leq  C\big(1+ \widetilde{\E}\|(\widetilde{\mu}_{\lambda})_\OO\|_{L^2(0,t)}^\frac{p}{2} + \widetilde{\E}\|\nabla\widetilde{\varphi}_{\lambda}\|_{L^2(0,t;H)}^p\big) + \|F_\lambda(\varphi_{0})\|_{L^1(\OO)}^\frac{p}{2} \\ 
        &+ \widetilde{\E}\left(\int_0^t\sum_{k=0}^\infty\int_\OO |F_\lambda^{\prime\prime}(\widetilde{\varphi}_{\lambda})||G(\widetilde{\varphi}_{\lambda})u_k|^2\d x\d s\right)^{\frac{p}{2}},
    \end{align*}
    where the constant $C$ is independent of $\lambda$. It remains to control the last two terms on the right-hand side uniformly in $\lambda$. Notice first that
    \begin{align*}
        \|F_\lambda(\varphi_0)\|_{L^1(\OO)} \leq \|F(\varphi_0)\|_{L^1(\OO)}
    \end{align*}
    by definition of $F_\lambda$. Next, invoking Hölder's inequality and recalling \ref{hyp:G}, we infer that 
    \begin{align*}
        \int_0^t\sum_{k=0}^\infty\int_\OO |F_\lambda^{\prime\prime}(\widetilde{\varphi}_{\lambda})||G(\widetilde{\varphi}_{\lambda})u_k|^2\d x\d s 
        &\leq \int_0^t\|F_\lambda^{\prime\prime}(\widetilde{\varphi}_\lambda(s))\|_{L^1(\OO)}\d s\sum_{k=0}^\infty\|g_k(\widetilde{\varphi}_\lambda)\|_{L^\infty(Q)}^2 \\
        &\leq |Q|L_G^2\int_0^t\|F_\lambda^{\prime\prime}(\widetilde{\varphi}_\lambda(s))\|_{L^1(\OO)}\d s.
    \end{align*}
    On account of these estimates as well as the growth assumption on $F^{\prime\prime}$, cf. \ref{hyp:potential}, we obtain
    \begin{align*}
        &\frac{1}{2}\widetilde{\E}\sups\mathcal{E}_{1,\lambda}(\widetilde{\varphi}_{\lambda}(s))^{\frac{p}{2}} + \widetilde{\E}\sups\|F_\lambda(\widetilde{\varphi}_{\lambda}(s))\|_{L^1(\OO)}^{\frac{p}{2}} + \widetilde{\E}\|\nabla\widetilde{\mu}_{\lambda}\|_{L^2(0,t;H)}^p + \widetilde{\E}\sups|(\widetilde{\mu}_{\lambda}(s))_\OO|^\frac{p}{2} \\
        &\leq C\big(1+ \widetilde{\E}\|(\widetilde{\mu}_{\lambda})_\OO\|_{L^2(0,t)}^\frac{p}{2} + \widetilde{\E}\|\nabla\widetilde{\varphi}_{\lambda}\|_{L^2(0,t;H)}^p + \widetilde{\E}\|F_\lambda(\widetilde{\varphi}_{\lambda})\|_{L^1(0,t;L^1(\OO))}^{\frac{p}{2}}\big) + \|F(\varphi_{0})\|_{L^1(\OO)}^\frac{p}{2},
    \end{align*}
    where the constant $C$ is independent of $\lambda$. Finally, Gronwall's lemma yields
    \begin{align}\label{Est:lam_2}
        \|\widetilde{\mu}_\lambda\|_{L^{\frac{p}{2}}(\tom;L^2(0,T;V))} + \|\nabla \widetilde{\mu}_\lambda\|_{L^p(\tom;L^2(0,T;H))} \leq C.
    \end{align}
    Thanks to \ref{hyp:G}, we also have
    \begin{align*}
        \|G(\widetilde{\varphi}_{\lambda})\|_{L^\infty(\tom\times(0,T);\cL^2(U,H))\cap L^p(\tom;L^2(0,T;\cL^2(U,V)))} \leq C,
    \end{align*}
    and thus, by \cite[Lemma 2.1]{fland-gat},
    \begin{align}\label{Est:lam_3}
        \left\|\int_0^. G(\widetilde{\varphi}_{\lambda}(s))\d \widetilde{W}(s)\right\|_{L^p(\tom;W^{k,p}(0,T;H))\cap L^2(\tom;H^{k}(0,T;V))} \leq C
    \end{align}
    for every $k\in(0,\frac{1}{2})$ and $p\geq 1$.

    Since $\gamma_\lambda$ is monotone and since $J$ is non-negative, it holds
     \begin{align*}
         \int_\OO\mathcal{L}\widetilde{\varphi}_\lambda\,\gamma_\lambda(\widetilde{\varphi}_\lambda)\d x = \frac{1}{2} \int_\OO\int_\OO J(x-y)(\widetilde{\varphi}_\lambda(x) - \widetilde{\varphi}_\lambda(y))(\gamma_\lambda(\widetilde{\varphi}_\lambda(x))-\gamma_\lambda(\widetilde{\varphi}_\lambda(y)))\d y\d x \geq 0.
     \end{align*}